\documentclass[10pt]{elsarticle}

\usepackage{amsmath}
\usepackage{amsthm}
\usepackage{amssymb}
\usepackage[margin=3cm]{geometry}
\usepackage{graphicx}
\usepackage{epstopdf}
\usepackage[dvipsnames]{xcolor}
\usepackage[hidelinks]{hyperref}
\usepackage[normalem]{ulem}
\usepackage[section]{algorithm}
\usepackage{algorithmic}
\usepackage{caption}
\usepackage{xcolor}
\newcommand{\bm}[1]{{\boldsymbol #1}} % bold math
\newcommand{\ii}{\bm{i}}
\newcommand{\pp}{\bm{p}}

\DeclareMathOperator{\Span}{span}
\DeclareMathOperator{\polylog}{polylog}

\usepackage{xspace}
\newcommand{\cossiga}{\textsc{CossIGA}\xspace} % \xspace adds an empty space after \cossiga, if necessary. 

\newtheorem{remark}{Remark}[section]

\newcounter{lnote}

\begin{document}

\begin{frontmatter}

  \title{Compressive Isogeometric Analysis} 

% \author{S. Brugiapaglia$^1$, L.  Tamellini$^2$,  M. Tani$^{2,3}$ \\[6pt]
%   {\footnotesize $^1$ SFU and/or CONCORDIA affiliation????} \\
%   {\footnotesize $^2$ Consiglio Nazionale delle Ricerche, Istituto di Matematica Applicata e Tecnologie Informatiche (CNR-IMATI)} \\
%   {\footnotesize $^3$ Universit\`a degli Studi di Pavia, Dipartimento di Matematica????} \\
% }

\author[label1]{Simone Brugiapaglia\corref{cor1}}
\ead{simone.brugiapaglia@concordia.ca}
\author[label2]{Lorenzo Tamellini}
\ead{tamellini@imati.cnr.it}
\author[label2]{Mattia Tani}
\ead{mattia.tani@imati.cnr.it}
\cortext[cor1]{Corresponding author}
\address[label1]{Department of Mathematics and Statistics, Concordia University. Montreal, Canada.}
\address[label2]{Consiglio Nazionale delle Ricerche, Istituto di Matematica Applicata e Tecnologie Informatiche ``E. Magenes'' (CNR-IMATI). Pavia, Italy.}

\begin{abstract}
This work is motivated by the difficulty in assembling the Galerkin matrix when solving 
Partial Differential Equations (PDEs) with Isogeometric Analysis (IGA) using B-splines of moderate-to-high polynomial degree. To mitigate this problem, we
propose a novel methodology named \cossiga (COmpreSSive IsoGeometric Analysis), which combines the IGA principle with CORSING, a recently introduced sparse recovery approach
for PDEs based on compressive sensing.
\cossiga assembles only a small portion of a suitable IGA Petrov-Galerkin
discretization and is effective whenever the PDE solution is
sufficiently sparse or compressible, i.e., when most of its
coefficients are zero or negligible.
The sparsity of the solution is promoted by employing a multilevel
dictionary of B-splines as opposed to a basis. 
Thanks to sparsity and the fact that only a
fraction of the full discretization matrix is assembled, the proposed technique has
the potential to lead to significant computational savings.
We show the effectiveness of \cossiga for the solution of the 2D and 3D Poisson
equation over nontrivial geometries by means of an extensive numerical investigation.

\end{abstract}


\end{frontmatter}

% \pagenumbering{arabic}
% \setcounter{page}{1}

\section{Introduction}

% In this paper we focus on a method which is designed to compute numerical approximation of solutions of PDE which
% have a certain sparsity, in the sense that if expanded on a suitable basis (actually, frame), most of the coefficients
% are zero. Whether or not it is reasonable to expect such a feature in PDE solution is therefore the first question
% to be addressed. 

% We propose a hybrid discretization scheme, obtained by combining the Compressed Solving method with the Iso-Geometric Analysis (IGA), named CossIGA (COmpreSSive IGA). 

Isogeometric Analysis (IGA) is an alternative to standard Finite Element Analysis (FEA) that 
has attracted considerable attention by researchers in computational science end engineering communities
since the seminal paper \cite{Hughes:2005}, published in 2005. IGA methodologies are quite similar to standard FEA,
with the main difference that the basis functions used for representing the domain and
the solution to the Partial Differential Equation (PDE) 
considered are splines rather than finite element basis functions.
This apparently simple change generates many interesting features.
These include the possibility of representing the geometry of the domain exactly,
a more flexible choice of the polynomial degree and of the regularity
of the basis used for approximating the solution, and a more effective error vs.\ degrees-of-freedom ratio than standard FEA.
We refer the interested reader to \cite{BBSV14} for a detailed mathematical analysis of the IGA method.

As all ``young methods'', many aspects of IGA are still the subject of scientific investigation.
The method proposed in this manuscript addresses one of these aspects,
i.e.\ the fact that the assembly and resolution of IGA Galerkin linear systems is usually
very expensive from the computational point of view, especially for moderate-to-high polynomial degrees of the spline basis functions. This aspect has been already tackled in several ways.
Proposed approaches include switching from the usual element-based quadrature to a function-based quadrature \cite{calabro2016fast,sangalli2018matrix},
more efficient matrix computation algorithms \cite{antolin2015efficient}, low rank and sparse grids techniques
\cite{beck.eal:sparse-IGA,Hofreither:2018,mantzaflaris2016low}, and efficient preconditioning
\cite{Collier2013,da2014isogeometric,donatelli2017symbol,hofreither2015robust,sangalli2016isogeometric,Tielen2019}.

In this paper, we propose instead a method in which the computational savings are potentially obtained by assembling
only a small fraction of the full Galerkin matrix. The crucial underlying assumption is the sparsity (or compressibility) of the solution. Namely, if expanded with respect to a suitable basis, most of the solution coefficients should be zero (or negligible). Whether or not it is reasonable to expect such a feature in the solution to a PDE is therefore the first question that should be addressed. As will be clearer later, one notable example is the case of PDEs whose solutions exhibit multiscale features. 

We mention in passing that another approach based on computing only a fraction of the Galerkin matrix has been recently proposed in \cite{drzisga2020surrogate}. However, in \cite{drzisga2020surrogate} the authors propose to compute a portion of the matrix exactly and then use this information to estimate the remaining entries of the matrix, so that, in the end, an approximation of the entire matrix is available. Instead, we propose to completely neglect some part of the matrix.

The method proposed in this paper is named \cossiga (COmpreSSive IsoGeometric Analysis).
This name refers to the fact that the method is an adaptation of the recently introduced CORSING (COmpRessed SolvING) method
\cite{brugiapaglia2016compressed,brugiapaglia2015compressed,brugiapaglia2018theoretical}
to the IGA framework.
CORSING combines the Petrov-Galerkin method with compressive sensing \cite{candes2006robust, donoho2006compressed}.
It assembles only a fraction of the discretization matrix and approximates the PDE
solution via sparse recovery, using techniques such as $\ell^1$ minimization or greedy algorithms.
In this work, we extend the previous CORSING works,
which were restricted to piecewise multilinear basis functions and rectangular domains, by using the
IGA B-splines-based machinery to represent both the domain and the solution.

One notable difference with previous CORSING works is that here we do not use a hierarchical basis as set of trial functions (although this could be a viable approach given the well-established theory
of hierarchical B-splines and their use in adaptive methods for PDEs - see \cite{bracco2019} and references therein).
In fact, we consider as trial functions a \emph{dictionary} of splines obtained as a union of spline bases at different refinement levels (we use the term dictionary to refer to a generic, possibly redundant, system).  The task of determining which splines should be activated is left to a sparse recovery algorithm. Thus, our approach can be seen as an alternative to standard adaptivity using hierarchical B-splines, where one does not need to implement any hierarchical basis, nor any marking/refining/derefining algorithms. Although redundant systems are widely employed in signal processing, their application in numerical analysis and scientific computing is a largely unexplored, yet very promising direction (see also \cite{adcock2019frames}).

A first important disclaimer is that the goal of the manuscript is to show the potential of \cossiga, but our implementation is not yet computationally effective. We actually prepared this initial manuscript with the idea of testing whether \cossiga has \emph{enough potential} to be worth implementing in a cost-effective way (which, by now, we believe it has). All numerical tests shown here have been implemented in GeoPDEs \cite{VAZQUEZ2016523}, available at \url{http://rafavzqz.github.io/geopdes/download/}.
The fact that we do not have an efficient implementation of \cossiga is the reason why in the numerical
tests we do not discuss computational times and only show abstract indicators of computational cost such as the spline refinement level and the number of computed coefficients.

A second disclaimer is that upon completion of the manuscript we became aware of the
work \cite{KANG201978}, which bears some similarities with \cossiga. Our work shares
with \cite{KANG201978} the idea that using a spline dictionary instead of a basis could promote sparsity/compressibility
of the PDE solution, and that a sparse version of the solution can then be recovered by suitable $\ell^1$ minimization
algorithms. However, in \cite{KANG201978} randomized selection of the rows is not present and the compressive sensing paradigm is therefore not fully exploited. In fact, the
entire Galerkin matrix is assembled, and not just a fraction of it as proposed here (in our
opinion, this is where most of the computational gain can be potentially obtained).
Furthermore, the numerical tests discussed in \cite{KANG201978} only include 2D square domains,
while here we take into account more general (2D and 3D) geometries. Finally, in \cite{KANG201978} sparse recovery is performed via $\ell^1$ minimization, whereas we employ the greedy algorithm orthogonal matching pursuit.

The rest of this paper is organized as follows. The methodology is explained in Section \ref{section:method}. In particular, a brief recap on the basics of IGA is given in Section \ref{subsection:IGA},
and the construction of the multilevel spline dictionary is detailed in Section \ref{subsection:frameandsparsity}.
%Moreover, the proper formal setting for the \cossiga method is actually not a plain Galerkin formulation, but a Petrov-Galerkin formulation, as explained in Section \ref{subsection:petrovG}. 
The Petrov-Galerkin formulation, which is the proper formal setting of the \cossiga method, is introduced in Section \ref{subsection:petrovG}.
Finally,
the two Sections \ref{subsection:cossiga} and \ref{subsection:effective_CossIGA} give a detailed
explanation of the specifics of \cossiga (the approach is summarized in Algorithm~\ref{alg:CossIGA}).
An extensive numerical investigation is then carried out in Section \ref{section:tests}:
in particular, we consider two test cases in 2D and one in 3D;
a test comparing the effectiveness of $C^0$ vs.\ $C^{p-1}$ splines ($p$ being their polynomial degree) is also provided.
Finally, Section \ref{section:conclusions} gathers conclusive remarks and future work directions.

\section{The \cossiga method}\label{section:method}

\subsection{Problem setting}\label{subsection:problem}
Let us consider the weak formulation of the homogeneous Poisson equation over a domain $\Omega \subseteq \mathbb{R}^d$, with $d = 2,3$
\begin{equation}
\label{eq:weak_problem}
\mathrm{find } \quad u \in H_0^1(\Omega): \quad a(u,v) = \int_\Omega f(\bm{x})v(\bm{x}) \mathrm{d}\bm{x}, \quad \forall v \in H_0^1(\Omega),
\end{equation}
where $a(\cdot,\cdot):H_0^1(\Omega) \times H_0^1(\Omega) \to \mathbb{R}$ is the bilinear form defined by
\begin{equation}
a(u,v) = \int_\Omega \nabla u(\bm{x}) \cdot \nabla v(\bm{x}) \mathrm{d}\bm{x},
\end{equation}
and where  $f \in H^{-1}(\Omega)$ is a forcing term. 
The method can be easily generalized to more general weak problems in Hilbert spaces, such as advection-diffusion-reaction equations. 

In this section we will formally introduce all the technical elements needed to define the CossIGA approach.
The main ideas employed are based on the CORSING method \cite{brugiapaglia2018theoretical,brugiapaglia2015compressed}.
In short, we will consider a Petrov-Galerkin discretization of \eqref{eq:weak_problem},
where the test space is randomly subsampled according the so-called local $a$-coherence
and sparse recovery of an approximate solution to \eqref{eq:weak_problem} is performed via Orthogonal Matching Pursuit (OMP).

\subsection{B-splines and the Isogeometric Analysis principle}\label{subsection:IGA}

%We start by briefly recalling the fundamentals of IGA and refer to \cite{BBSV14} for a more thorough discussion. 
Given two natural numbers $n,p \in \mathbb{N}$, we define the \emph{knot vector}
over the unit interval  $\widehat{I}:=[0,1]$  as $\bm{\Xi} = [\xi_1,\xi_2,\ldots,\xi_{n+p+1}]$
with nondecreasing and possibly repeated entries, such that $\xi_1 = 0$ and $\xi_{n+p+1} = 1$. Each
$\xi_i$ is a \emph{knot} and any interval $(\xi_i,\xi_{i+1})$ having nonzero length is an \emph{element}.
Let us further denote the total number of elements as $N_{\mathrm{el}}$.
In this paper, the elements will have the same length, called \emph{mesh size} and denoted by $h$.
Moreover, we assume the knot vector $\bm{\Xi}$ to be \emph{open},
i.e.\ we let its first and last knots have multiplicity $p + 1$ (i.e., they are repeated $p+1$ times).
Observe that also internal knots could have multiplicity greater than one.
Finally, we define the nondecreasing vector $\bm{Z} = [\zeta_1, \ldots , \zeta_{N_{\mathrm{el}} +1}]$
as the vector of knots of $\bm{\Xi}$ without repetitions,
and let $m_i$ be the multiplicity of $\zeta_i$ in $\bm{\Xi}$, so that $\sum_{i = 1}^{N_{\mathrm{el}}+1} m_i = n + p + 1$.

Given the knot vector $\bm{\Xi}$ thus built, we define the B-splines by means of the Cox--De Boor recursive formula. 
% \begin{align}
% \widehat{B}_{i,0}(\xi) 
% & = 
% \begin{cases}
% 1, & \xi_i \leq \xi < \xi_{i+1},\\
% 0, & \mathrm{otherwise},
% \end{cases}\\
% \widehat{B}_{i,p}(\xi) 
% & = \begin{cases}
% \displaystyle\frac{\xi - \xi_i}{\xi_{i+p}-\xi_i} \widehat{B}_{i,p-1}(\xi) + \frac{\xi_{i+p+1} - \xi}{\xi_{i+p+1}-\xi_{i+1}} \widehat{B}_{i+1,p-1}(\xi),
% & \xi_i \leq \xi < \xi_{i+p+1},\\
% 0, & \mathrm{otherwise},
% \end{cases}
% \end{align}
% where we adopt the convention $0/0 = 0$; note that the basis corresponding to an open knot vector will be interpolatory in the first and last knots.
We start with piecewise constant splines
\[
\widehat{B}_{i,0}(\xi)= 
\begin{cases}
  1, &  \xi_{i}\leq \xi<\xi_{i+1}, \\ 
  0, &  \textrm{otherwise,} \\
\end{cases} \qquad \qquad  \mbox{for } i = 1,\ldots,n+p.
\]
Then, for $\tilde{p} = 1,\ldots,p$, we have the recursive step 
\[
  \widehat{B}_{i,\tilde{p}}(\xi)=
  \begin{cases}
    \dfrac{\xi-\xi_{i}}{\xi_{i+\tilde{p}}-\xi_{i}}\widehat{B}_{i,\tilde{p}-1}(\xi)+\dfrac{\xi_{i+\tilde{p}+1}-\xi}{\xi_{i+\tilde{p}+1}-\xi_{i+1}}\widehat{B}_{i+1,\tilde{p}-1}(\xi),
       & \xi_{i}\leq \xi<\xi_{i+\tilde{p}+1}, \\
       0, & \textrm{otherwise},
     \end{cases}
     \;\; \mbox{ for } i = 1,\ldots,n+p-\tilde{p},
\]  
with the convention that $0/0=0$. Note that if the knot vector $\bm{\Xi}$ is open, the corresponding basis
is interpolatory in the first and last knots.
  The B-splines just defined form a basis for the space $S_p(\bm{\Xi},\widehat{I})$ of
  spline-polynomials, i.e., of piecewise polynomials of degree $p$ and regularity $C^ {p-m_i}$ at $\zeta_i$,
$$
S_p(\bm{\Xi},\widehat{I}) 
= \mathrm{span}\{\widehat{B}_{i,p}: i=1,\ldots,n\}.
$$
In particular, the maximal regularity of a spline at the knots  is $C^{p-1}$.  
In the following, we will be interested in the cases where
all the internal knots of $\bm{\Xi}_l$ are either repeated once or $p-1$ times:
in the former case, we talk about spline of maximal regularity or $C^{p-1}$ splines,
while in the latter we talk about $C^{0}$ splines. 
Moreover, since we are considering homogeneous boundary conditions, we consider the set of B-splines that vanish at the boundary: 
$$
S_p^{\mathrm{int}}(\bm{\Xi},\widehat{I}) 
= \mathrm{span}\{\widehat{B}_{i,p}: i=2,\ldots,n-1\}.
$$

For $d = 2$ we define the parametric domain
$\widehat{\Omega} = \widehat{I} \times \widehat{I}$ (extension to the
case $d > 2$ is analogous). We consider two open knot vectors
$\bm{\Xi}_1, \bm{\Xi}_2$ with $n_1 + p_1 + 1$ and $n_2 + p_2 + 1$ knots
respectively, the corresponding knots without repetitions $Z_1$,
$Z_2$, and the tensor products $\bm{\Xi} = \bm{\Xi}_1 \times \bm{\Xi}_2$,
$\bm{Z} = Z_1 \times Z_2$; in particular, $\bm{Z}$ generates a cartesian mesh
over $\widehat{\Omega}$ composed of $N_{\mathrm{el},1} N_{\mathrm{el},2}$
rectangular elements. Taking tensor products of the univariate
B-splines over $\bm{\Xi}_1$ and $\bm{\Xi}_2$ we obtain a basis for the space of
bivariate spline polynomials and the corresponding basis of B-splines satisfying
homogeneous boundary conditions. To this end, we introduce the multi-indices
  $\ii = (i_1,i_2), \pp = (p_1, p_2)$ and let $\widehat{B}_{\ii,\pp}(\xi_1,\xi_2)=\widehat{B}_{i_1,p_1} (\xi_1)\widehat{B}_{i_2,p_2} (\xi_2)$,
  so that
\begin{align*}
  S_\pp(\bm{\Xi} , \widehat{\Omega} )
  & = \mathrm{span}\{\widehat{B}_{\ii,\pp}:  1 \leq i_1 \leq n_1, 1 \leq i_2 \leq n_2 \},
   % = S_p(\bm{\Xi}_1,\hat{I}) \otimes S_p(\bm{\Xi}_2,\hat{I}), 
   \\
  S_\pp^{\mathrm{int}}(\bm{\Xi} , \widehat{\Omega} )
  & = \mathrm{span}\{\widehat{B}_{\ii,\pp}:  2 \leq i_1 \leq n_1-1, 2 \leq i_2 \leq n_2-1 \}.
   % = S_p^{\mathrm{int}}(\bm{\Xi}_1,\hat{I}) \otimes S_p^{\mathrm{int}}(\bm{\Xi}_2,\hat{I}).
\end{align*}
In the following, we will assume that $p_1=p_2=p$,
  so that we can drop the bold notation $\pp$ and write $p$ instead.
  Moreover, we enumerate the B-splines with a single index $i$ ranging from $1$ to $n = n_1 n_2$
  for $S_p(\bm{\Xi} , \widehat{\Omega})$, and from $1$ to $n^{\mathrm{int}} = (n_1-1) (n_2-1)$
  for $S_p^{\mathrm{int}}(\bm{\Xi} , \widehat{\Omega} )$, i.e.,
  $$
  S_p(\bm{\Xi} , \widehat{\Omega} ) = \mathrm{span}\{\widehat{B}_{i,p}:  1 \leq i \leq n \},
  \quad 
  S_p^{\mathrm{int}}(\bm{\Xi} , \widehat{\Omega} ) = \mathrm{span}\{\widehat{B}_{i,p}: 1 \leq i \leq n^{\mathrm{int}} \}.
  $$

We assume that the computational domain $\Omega$ can be parameterized by  an invertible mapping 
$F : \widehat{\Omega} \to \Omega$, obtained as a
linear combination of B-splines with given control points $\mathbf{P}_1,\ldots,\mathbf{P}_n \in \mathbb{R}^2$, i.e.,
$$
\bm{x} \in \Omega \Longleftrightarrow 
\bm{x}=F(\bm{\xi}) := \sum_{1 \leq i \leq n} \widehat{B}_{i,p}(\bm{\xi}) \mathbf{P}_i , \quad \text{ for some } \bm{\xi}  \in \widehat{\Omega}.
$$
We mention in passing that many geometries $\Omega$ of practical interest, such as circles and ellipses,
cannot be represented exactly by B-splines. However, 
\emph{nonuniform rational B-splines} (NURBS) can be employed for this purpose (see \cite{BBSV14,Hughes:2005}
for details). As the name suggests, NURBS are ratios of B-splines and retain most of the properties of B-splines,
so in the rest of this manuscript we use ``splines'' as a comprehensive term for both B-splines and NURBS.

According to the IGA principle, splines are also employed to approximate the solution $u$ to the weak problem \eqref{eq:weak_problem}.
To this end, we introduce the splines on the physical domain, defined by
$$
B_{i,p} := \widehat{B}_{i,p} \circ F^{-1},
$$
and the spline space on the physical domain $\Omega$ as follows:
$$
S_p(\bm{\Xi},\Omega)=\mathrm{span}\{B_{i,p}: 1 \leq i \leq n\}, \quad S_p^{\mathrm{int}}(\bm{\Xi},\Omega)=\mathrm{span}\{B_{i,p}: 1 \leq i \leq n^{\mathrm{int}}\}.
$$ 
Finally, we define a basis of $S_p^{\mathrm{int}}(\bm{\Xi},\Omega)$ normalized with respect to the $H^1(\Omega)$-seminorm as follows:
$$
\mathcal{B}_p^{\mathrm{int}}(\bm{\Xi},\Omega) = \left\{\frac{B_{i,p}}{|B_{i,p}|_{H^1(\Omega)}} : 1 \leq i \leq n^{\mathrm{int}} \right\}.
$$

\subsection{Multilevel dictionary of B-splines and the sparsity assumption}\label{subsection:frameandsparsity}

\begin{figure}[t]
  \centering
  \includegraphics[width=\linewidth]{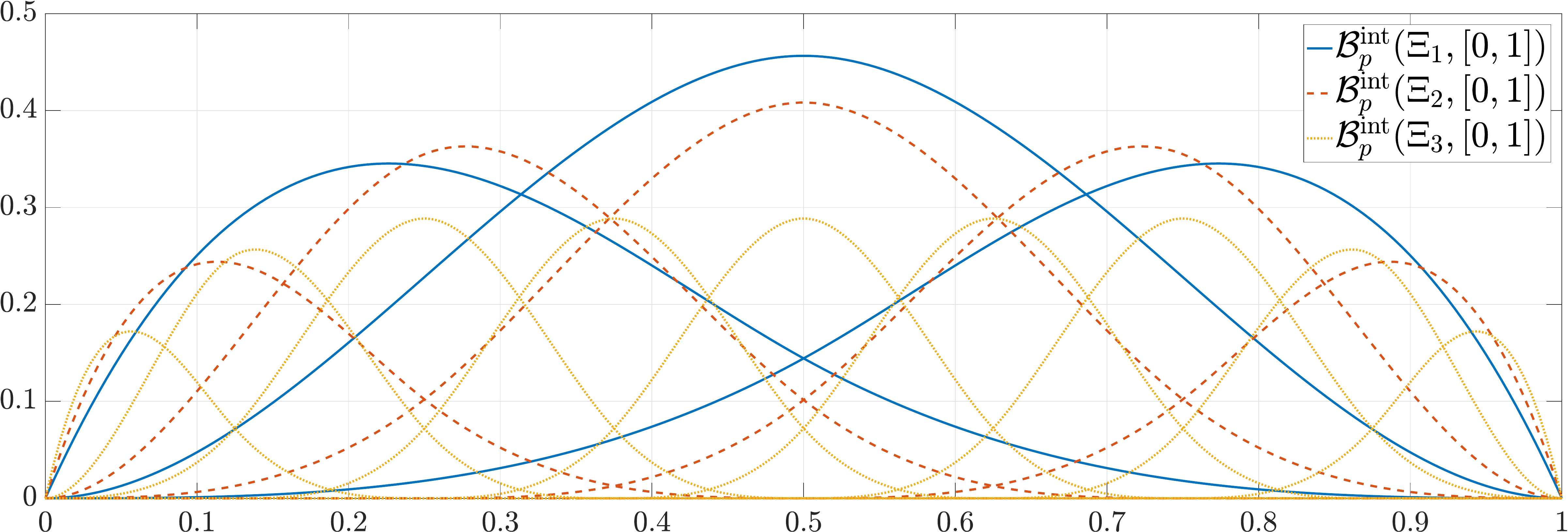}
  \caption{The dictionary $\Psi_{p,l_0,L}$, for $p=3, l_0=1, L=3$, and splines of regularity $C^2$.}
  \label{fig:frame}
\end{figure}

To apply the compressive sensing  principle, we need to generate a sparse (or compressible) representation of the solution $u$ to \eqref{eq:weak_problem}.
Namely, we need to identify a basis or, more in general, a dictionary such that most of the coefficients of the corresponding expansion of $u$ are zero (or negligible).
With this aim, we resort to a multiscale decomposition that is able to enhance compressibility
in solutions with, e.g., local features or sharp transitions. Given $l_0,L \in \mathbb{N}$ such that $1 \leq l_0 < L$, we consider the multilevel dictionary of B-splines
\begin{equation}
\label{eq:dictionary}
\Psi_{p,l_0,L} 
:= \bigcup_{l = l_0}^{L} \mathcal{B}^{\mathrm{int}}_p(\bm{\Xi}_l,\Omega)
= \{\psi_j\}_{j \in [N_{\mathrm{dict}}]},
\end{equation}
where $\bm{\Xi}_{l_0} \subseteq \cdots \subseteq \bm{\Xi}_L$ is a nested sequence of knot vectors
such that $\bm{\Xi}_l$ corresponds to a grid of meshsize $h_l = 2^{-l}$ associated with B-splines of degree $p$
and where we adopted the notation $[k] = \{1,\ldots,k\}$, for every $k \in \mathbb{N}$. The dictionary $\{\psi_j\}_{j \in [N_{\mathrm{dict}}]}$ is assumed to be ordered lexicographically with respect to the multi-index $(l,i)$: the level $l$ and the index $i$ of each (normalized) spline $B_{i,p}$ in $\mathcal{B}^{\mathrm{int}}_p(\bm{\Xi}_l,\Omega)$.
A plot of the dictionary $\Psi_{3,1,3}$ in 1D is shown in Figure \ref{fig:frame}.
The intuition is the following: all (or most of) the splines in the lowest levels are  activated
to approximate the coarse component of the solution and only a few splines in the high-resolution
levels are activated to capture local features or sharp transitions.
The dictionary $\Psi_{p,l_0,L}$ has been also considered in \cite{KANG201978}.

In standard hierarchical approaches, only a linearly independent subset of $\Psi_{p,l_0,L}$ is selected.
This is typically done by starting from the coarsest basis, then marking a region where the error is concentrated, refining in that region while preserving linear independence, and repeating this process until the solution is accurate enough.
Instead, in our approach the splines in the dictionary to be activated are identified by solving a particular minimization problem, as we will discuss in the following. However, before doing this, we need to make a short digression about the cardinality and the number of degrees of freedom (i.e., the dimension of the span) of the dictionary $\Psi_{p,l_0,L}$.

In our numerical tests (Section~\ref{section:tests}) we will use either $C^{p-1}$ splines or $C^{0}$ splines. 
Hence, we recall explicit formulas for the cardinality of $\mathcal{B}_p^{\mathrm{int}}(\bm{\Xi}_l,\Omega)$ in these two cases:
%Note that the cardinality of the generic $\mathcal{B}_p^{\mathrm{int}}(\bm{\Xi}_l,\Omega)$ is
\begin{equation}
\label{eq:cardinality_B_p^int}
|\mathcal{B}_p^{\mathrm{int}}(\bm{\Xi}_l,\Omega)| =
\begin{cases}
\left(2^l+p-2\right)^d, & \text{ for $C^{p-1}$ splines,}\\  
\left(2^l \, p - 1\right)^d, & \text{ for $C^{0}$ splines.}
\end{cases}
\end{equation}
Clearly, the cardinality of the dictionary is
\begin{equation}
\label{eq:cardinality_dictionary}
N_{\mathrm{dict}} 
= |\Psi_{p,l_0,L}|
= \sum_{l = l_0}^L |\mathcal{B}^{\mathrm{int}}_p(\bm{\Xi}_l,\Omega)|.
% = \sum_{l = l_0}^L (2^l+p-2)^2.
\end{equation}
%At each hierarchical level, the basis $\mathcal{B}^{\mathrm{int}}_p(\bm{\Xi}_l,\Omega)$ contains $N_l = (2^l +p-2)^2$ trial functions.
Notice that the span of the last hierarchical level coincides with the span of the whole dictionary,
due to the fact that B-splines in $\mathcal{B}^{\mathrm{int}}_p(\bm{\Xi}_l,\Omega)$ are linear combinations
of B-splines in $\mathcal{B}^{\mathrm{int}}_p(\bm{\Xi}_L,\Omega)$, for every $l \leq L$. Namely,
\begin{equation}
\label{eq:span_dictionary=span_level_L}
\mathrm{span}(\Psi_{p,l_0,L}) = S_p^{\mathrm{int}}(\bm{\Xi}_L,\Omega).
\end{equation}
Therefore, combining \eqref{eq:cardinality_B_p^int} and \eqref{eq:span_dictionary=span_level_L}, the number of degrees of freedom (dof)
of the dictionary $\Psi_{p,l_0,L}$ coincides with the cardinality of the basis at level $L$,  i.e., 
\begin{equation}
\label{eq:NDOFs}
N_{\mathrm{dof}} = 
\begin{cases}
(2^L+p-2)^d & \text{for $C^{p-1}$ splines,}\\  
(2^L \, p - 1)^d & \text{for $C^{0}$ splines.}\\
\end{cases}
\end{equation}
The size of the dictionary is in general comparable to the size of the basis at level $L$. In fact, it is not difficult to show that, for any $d$ and $p$,  $N_{\mathrm{dict}} \leq 3 N_{\mathrm{dof}}$ for $C^{p-1}$ splines (for $L$ large enough)  and $N_{\mathrm{dict}} \leq 2 N_{\mathrm{dof}}$ for $C^0$ splines (for any $L$).\footnote{For $C^{p-1}$ splines, $$N_{\mathrm{dict}}
= N_{\mathrm{dof}} + \sum_{l = l_0}^{L-1}(2^l+p-2)^d 
\leq N_{\mathrm{dof}} + \left(\sum_{l = l_0}^{L-1}(2^l+p-2)\right)^d 
= N_{\mathrm{dof}} + \left(2^L - 2^{l_0} + (L-l_0)(p-2)\right)^d.$$ Therefore, if $p \geq 2$ and $L$ is large enough to have $(L-l_0-2^{1/d})(p-2) \leq 2^L (2^{1/d}-1) +2^{l_0}$, we obtain $N_{\mathrm{dict}} \leq 3 N_{\mathrm{dof}}$. The computation is similar for $p=1$ or for $C^0$ splines.}

We can now come back to the main topic of introducing the minimization approach to  select the B-splines from the dictionary
  to be used to represent the solution $u$. As already said, we considered a multilevel dictionary hoping that
  only a few B-splines will be needed in order to well approximate $u$: most of those in the lowest levels (which are not many) will be used to approximate the coarse component
  of the solution and only a few splines in the higher levels will be used to capture local features or sharp transitions. In other words,
we aim at computing a sparse approximation $\widetilde{u}$ to $u$, i.e., a function of the form
\begin{equation}
\label{eq:def_u_tilde}
\widetilde{u} = \sum_{j =1}^{N_{\mathrm{dict}}} \widetilde{u}_j \psi_j, \quad \text{with } \|\widetilde{\bm{u}}\|_0  \ll N_{\mathrm{dict}},
\end{equation}
where, for every $\bm{v} \in \mathbb{R}^k$, $\|\bm{v}\|_0 := |\{v_j \neq 0\}|$. More specifically, we say that $\widetilde{u}$ is $s$-sparse if $\|\widetilde{\bm{u}}\|_0 \leq s$. Given a budget of $s$ coefficients, the goal is to compute an $s$-sparse approximation such that $\|u-\widetilde{u}\|_{H^1(\Omega)}$ is as close as possible to the best $s$-term approximation error of $u$ with respect to $\Psi_{p,l_0,L}$, defined by
\begin{equation}
\label{eq:best_s-term}
\sigma_{s}(u)_{H^1(\Omega)} = \inf_{\|\bm{z}\|_0 \leq s} \bigg\|u - \sum_{j \in [N_{\mathrm{dict}}]} z_j \psi_j\bigg\|_{H^1(\Omega)}.
\end{equation}
If $\sigma_{s}(u)_{H^1(\Omega)}$ has a fast decay with respect to $s$ (e.g., $\sigma_{s}(u)_{H^1(\Omega)} \leq C s^{-\alpha}$ for some $C,\alpha>0$), $u$ is informally said to be \emph{compressible} with respect to $\Psi_{p,l_0,L}$.  

\subsection{Petrov-Galerkin: B-splines vs.\ sine functions}\label{subsection:petrovG}

Together with the sparsity-promoting dictionary just introduced, we consider
a Petrov-Galerkin (PG) discretization of \eqref{eq:weak_problem}. We use the functions in the dictionary $\Psi_{p,l_0,L}$ as trial functions of the PG formulation, i.e.\ we approximate $u$ as
a linear combination of functions in $\Psi_{p,l_0,L}$.
We choose the test functions according to a principle that lies at the core of compressive sensing and also employed in CORSING. Namely, since the trial functions are localized in the space domain, it is convenient to choose test functions localized in the frequency (or Fourier) domain.
The underlying intuition is that functions that are sparse in the space domain cannot be too sparse in the frequency domain
(this is the so-called uncertainty principle \cite{donoho1989uncertainty}).
In our setting, we employ test functions of Fourier type to measure the solution in the frequency domain. Now, thanks to sparsity the amount of information intrinsically needed to represent the solution is very small; yet, the information in the frequency domain is spread over the whole spectrum due to the uncertainty principle. Therefore, Fourier measurements of signals that are sparse in space are highly redundant.  In order to get rid of this redundancy, the idea of compressive sensing is to select only a few of them in a randomized way.

For this reason, we consider the sine functions over $\widehat{I}^d$, defined by 
\begin{equation}\label{eq:sin_basis_fun_with_vector_index_rr}
% s_r(\xi) := \sin(r_1 \pi \xi_1) \sin(r_2 \pi \xi_2), \quad \forall r \in \mathbb{N}^2
\sin_{\bm{r}}(\bm{\xi}) 
:= \prod_{i=1}^d \sin(r_i \pi \xi_i), \quad \forall \bm{\xi} \in \widehat{I}^d,\; \forall \, \bm{r} \in \mathbb{N}^d.  
\end{equation}
Given a maximum frequency $R \in \mathbb{N}$, the corresponding basis of test functions defined over $\widehat{\Omega}$ is
\begin{equation}\label{eq:sin_basis_fun_with_scalar_index_q}
\Phi_{R}
:= \left\{\frac{\sin_{\bm{r}} \circ F^{-1}}{\left|\sin_{\bm{r}} \circ F^{-1}\right|_{H^1(\Omega)}} : \bm{r} \in [R]^d\right\} 
= \{\varphi_q\}_{q \in[N_{\mathrm{test}}]},  
\end{equation}
where the definition of the set $\{\varphi_q\}_{q \in [N_{\mathrm{test}}]}$ implicitly depends on the ordering used over the multi-index set $[R]^d$ (e.g., the lexicographic ordering) and where
%where $N_{\mathrm{test}} = R^2$. 
$\left|\Phi_{R}\right| =: N_{\mathrm{test}} = R^d $.
The resulting PG discretization of \eqref{eq:weak_problem} with respect to the trial and test functions in $\Psi_{p,l_0,L}$ and $\Phi_{R}$, respectively, is
\begin{equation}
\label{eq:PG_system}
B \bm{z} = \bm{c},
\end{equation}
where $B \in \mathbb{R}^{N_{\mathrm{test}} \times N_{\mathrm{dict}}}$ and $\bm{c} \in \mathbb{R}^{N_{\mathrm{test}}}$ are defined as
\begin{equation}
\label{eq:def_B_c}
B_{qj} := a(\psi_j, \varphi_q), \quad c_q := \int_\Omega f \varphi_q, \quad \forall j \in [N_{\mathrm{dict}}], \; \forall q \in [N_{\mathrm{test}}].
\end{equation}
A sufficient requirement to have a well-posed PG formulation is the following discrete inf-sup condition (see, e.g., \cite[Theorem 5.3.1]{quarteroni2008numerical}):
\begin{equation}
\label{eq:discrete_infsup}
\inf_{u \in \Span(\Psi_{p,l_0,L})} \sup_{v \in \Span(\Phi_R)} \frac{a(u,v)}{\|u\|_{H^1(\Omega)} \|v\|_{H^1(\Omega)}} \geq \alpha > 0.
\end{equation}
Note that \eqref{eq:discrete_infsup} is a condition on the vector spaces spanned by the trial and the test functions.
% In view of  \eqref{eq:span_dictionary=span_level_L} and \eqref{eq:NDOFs},
A necessary condition to have $\alpha >0$ is 
\begin{equation}
\label{eq:N_test_inf-sup}
N_{\mathrm{test}} \geq N_{\mathrm{dof}}. %= (2^L +p -2)^2.
\end{equation}
Moreover, $\alpha$ is nondecreasing with respect to $R$ or, equivalently, to $N_{\mathrm{test}}$. In practice,
in view of \eqref{eq:cardinality_B_p^int}, we make the heuristic choice
\begin{equation}
\label{eq:choice_R}
R =
\begin{cases}
1.5 \lceil 2^L+p-2 \rceil & \text{for $C^{p-1}$ splines,}\\  
1.5 \lceil 2^L \times p - 1 \rceil & \text{for $C^{0}$ splines.}
\end{cases}
\end{equation}
where $\lceil x \rceil$ rounds a real number $x$ to the closest integer greater than or equal to $x$. 
The factor 1.5 in Equation \eqref{eq:choice_R},
which in particular implies condition \eqref{eq:N_test_inf-sup}, has been empirically chosen based on numerical experimentation.
Studying the relation between $\alpha$ and $R$ from the theoretical standpoint is an open problem. Some theoretical insights on this issue are given by the so-called restricted inf-sup property analysis, introduced in \cite{brugiapaglia2018theoretical}.

\begin{figure}[t!]
  \centering
  \includegraphics[width=0.3\linewidth]{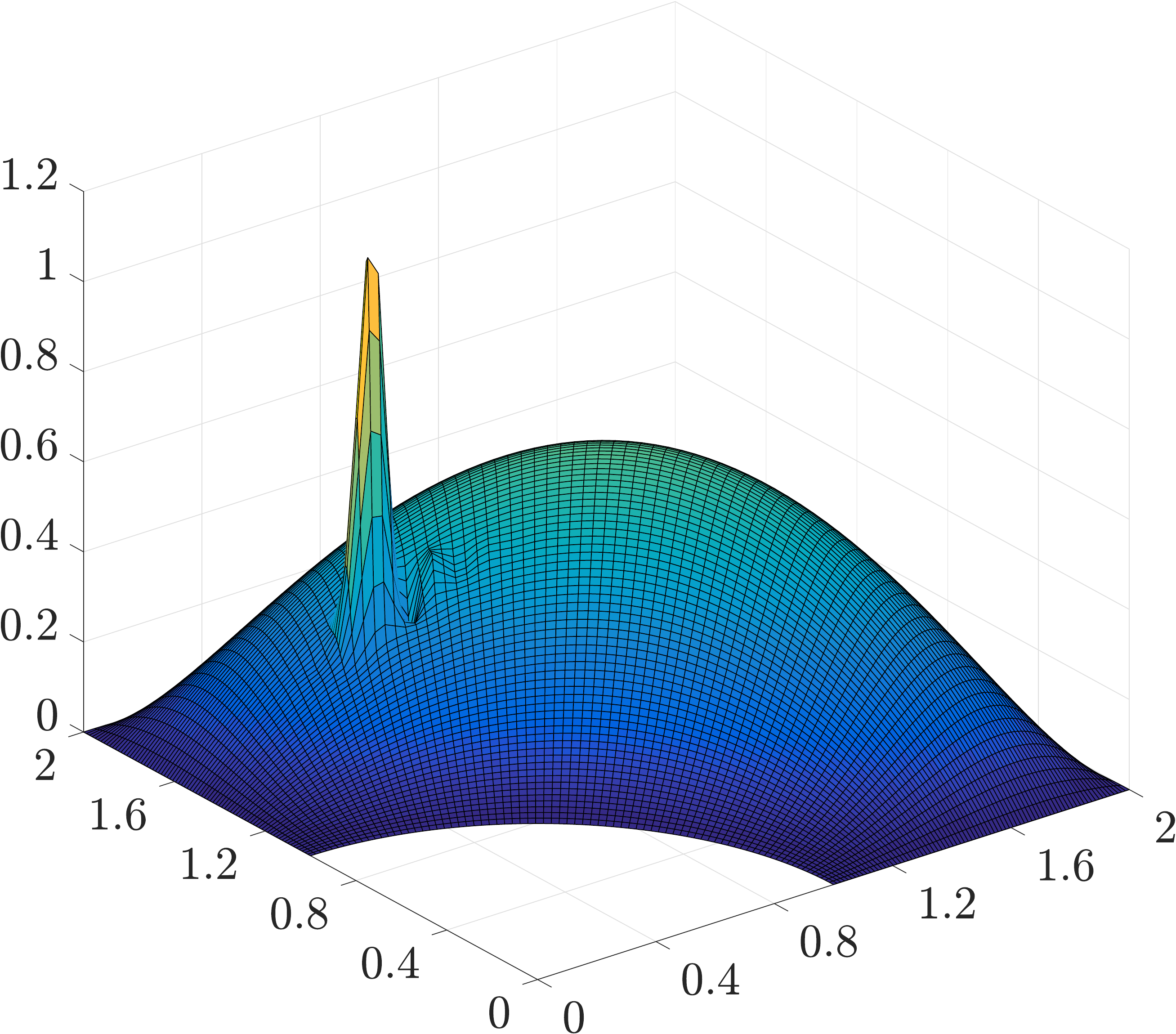} \quad 
  \includegraphics[width=0.3\linewidth]{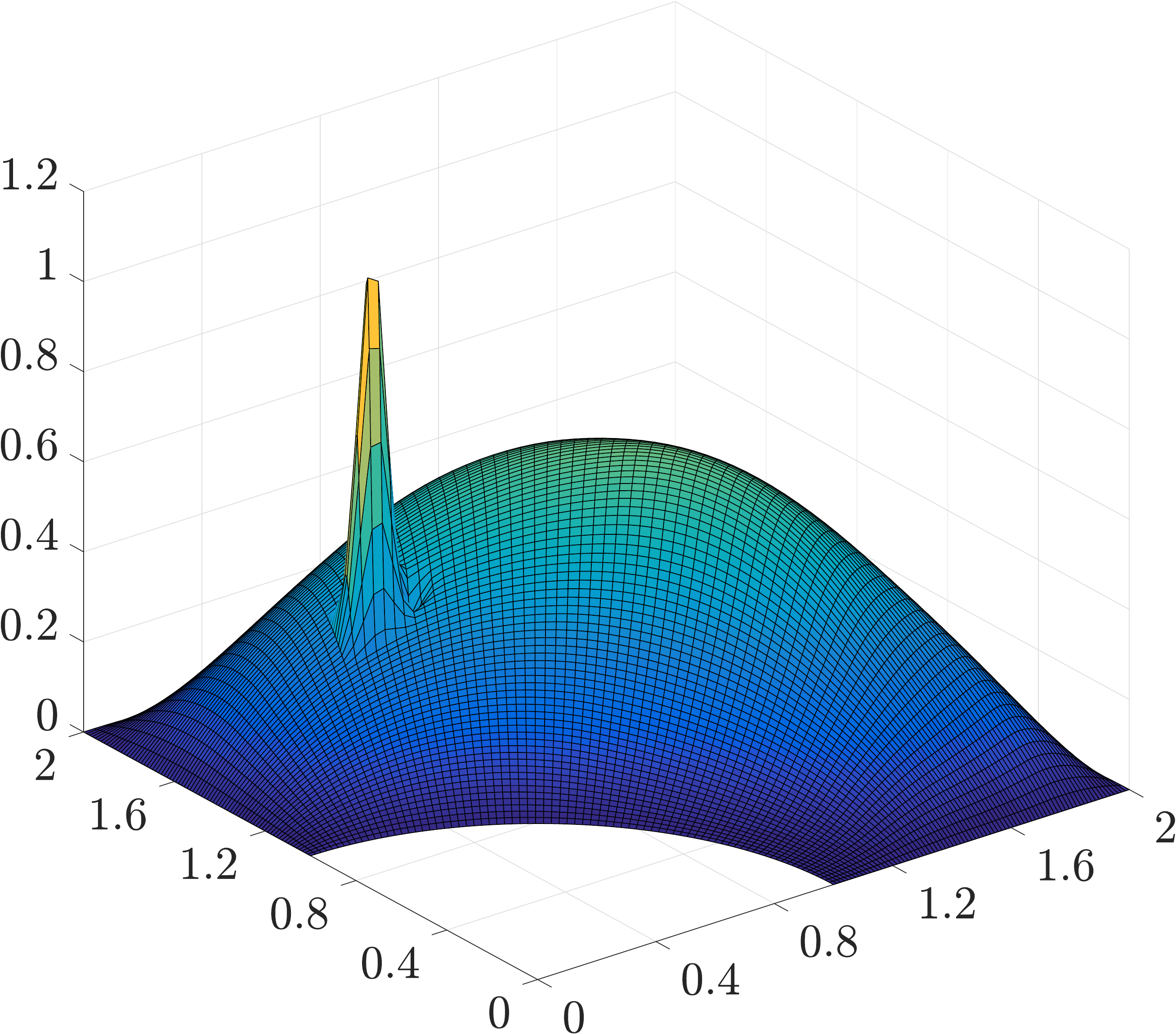} \quad
  \includegraphics[width=0.3\linewidth]{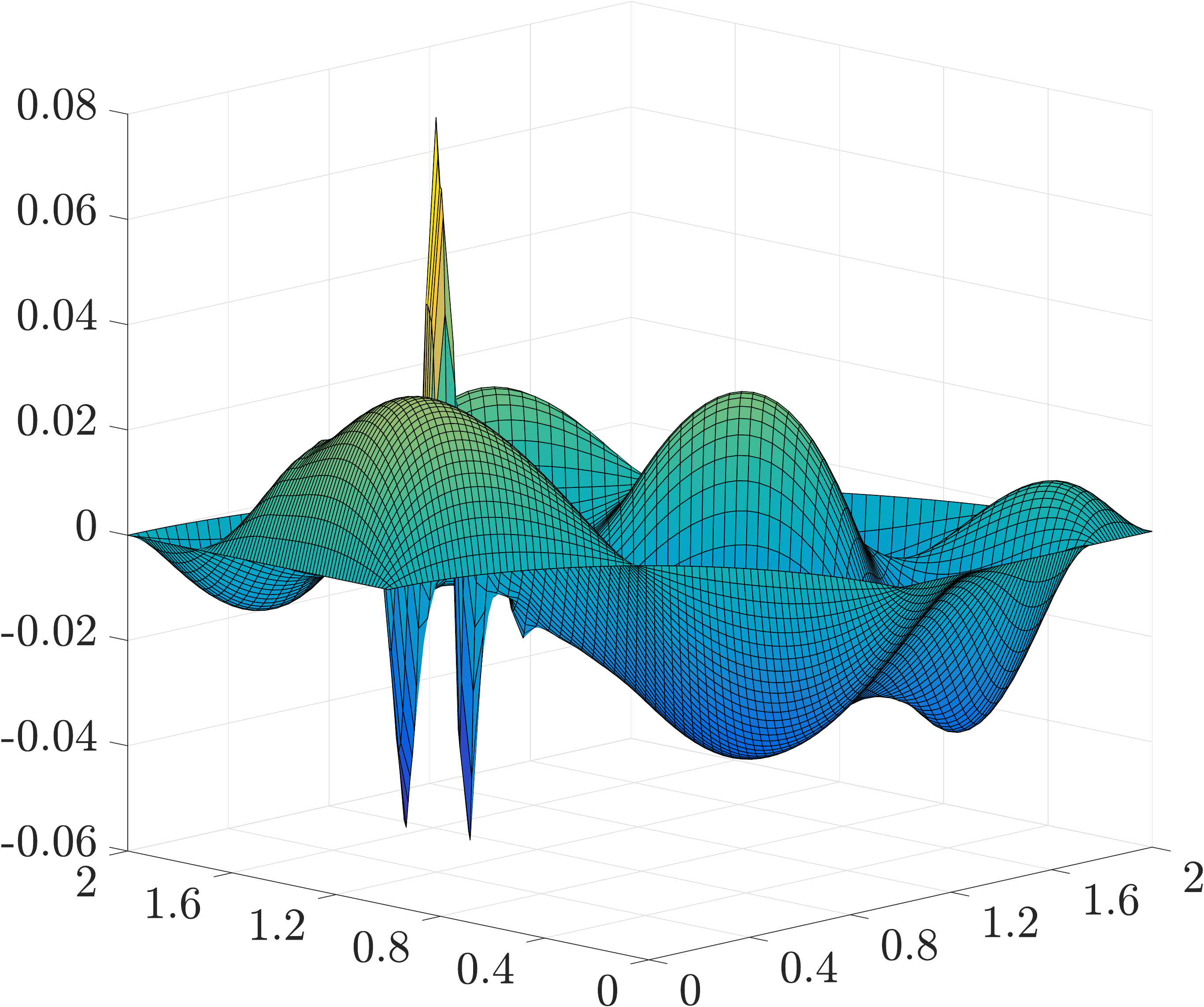}\\[16pt]
  \includegraphics[width=0.4\linewidth]{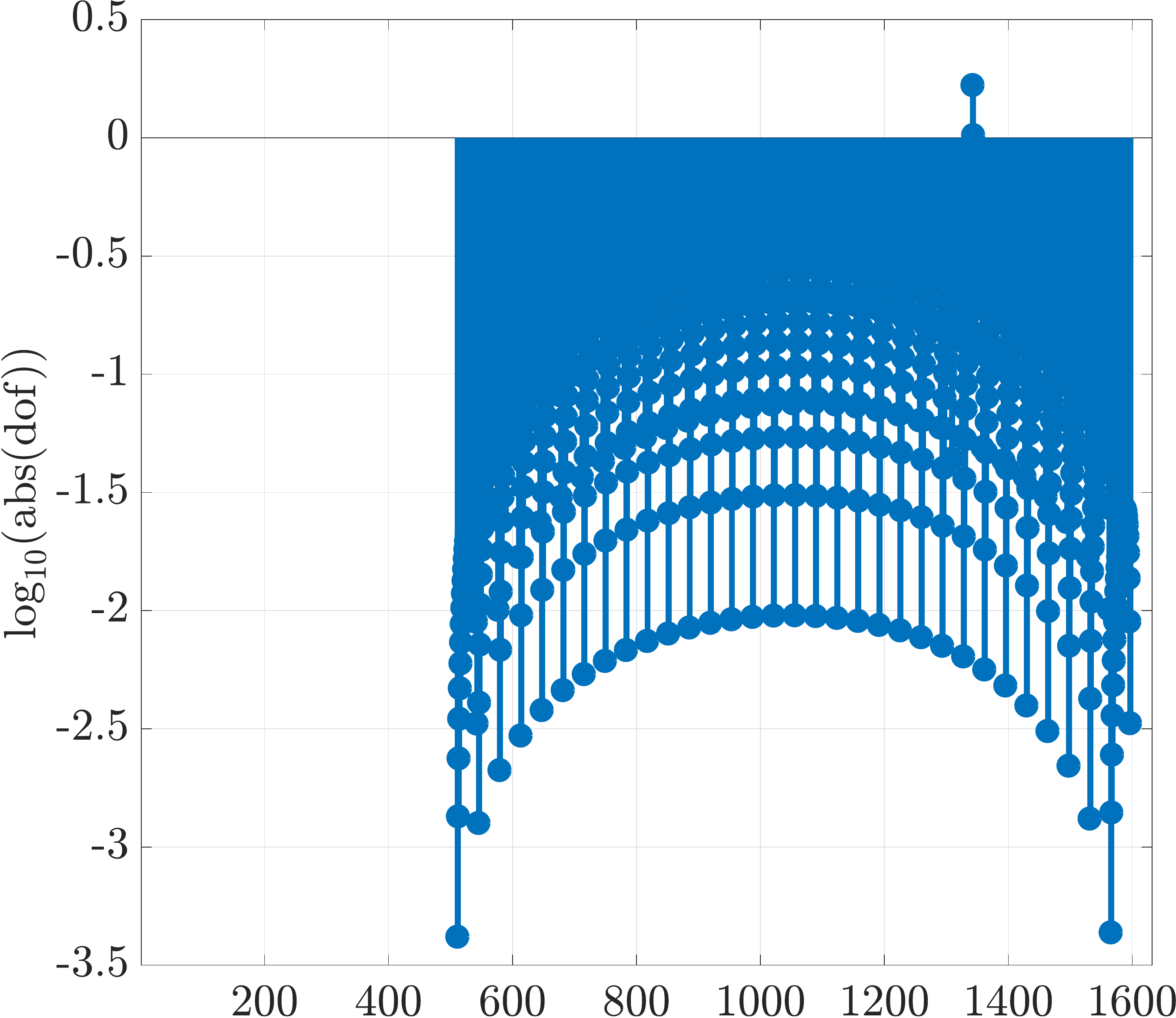} \quad \quad
  \includegraphics[width=0.4\linewidth]{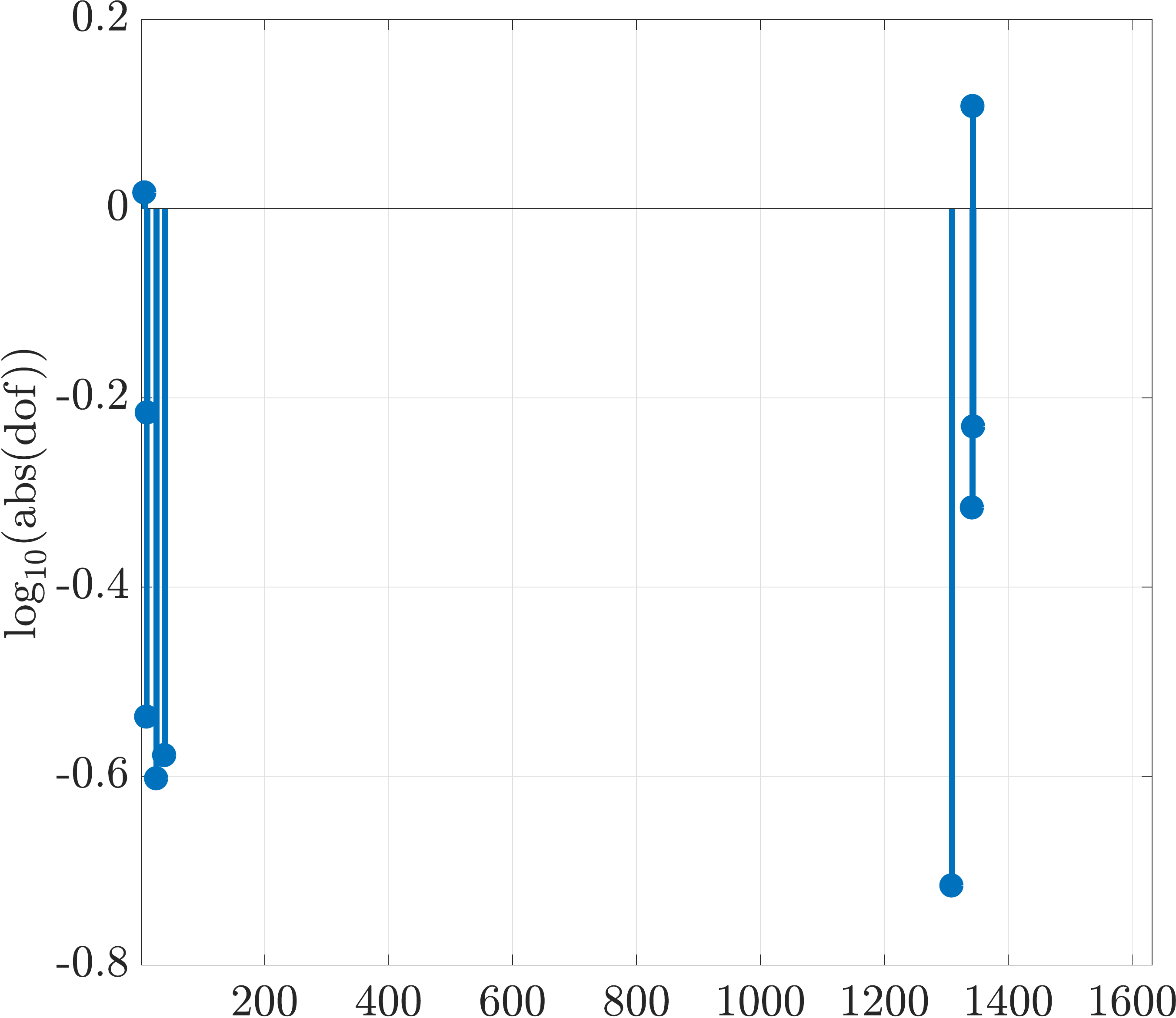}\\
  \caption{In this picture, we compare the standard IGA-Galerkin solution of the problem discussed in Section~\ref{subsection:polygauss}
    (i.e., the solution computed by solving the square linear system obtained using
    as trial and test functions the B-splines at the finest discretization level $\bm{\Xi}_L$)
    with the solution computed using OMP to solve the PG system (\ref{eq:PG_system})
    in the sense of equation (\ref{eq:sparse_recovery}); here we have enforced OMP to compute $s=9$ coefficients.
    The top row shows the computed solutions with both approaches
    (IGA: left panel; PG-OMP: central panel), and their difference in the right panel.
      The vertical axis of right panel allows to appreciate the fact that the difference is substantially smaller in size than the solution itself.
    The bottom row shows the size of the coefficients in log scale (so-called ``stem plots''): the left panel is for IGA and the right panel for PG-OMP.
    The horizontal axis in both  stem plots
    shows the lexicographically-ordered indices of the B-splines in the dictionary
    $\Psi_{p,l_0,L}$ with $p=2,l_0=1,L=5$ ($N_{\mathrm{dict}} = 1632$ and $N_{\mathrm{dof}} = 1024$),
    with the understanding that the IGA-Galerkin solution only uses
    the basis at level $L=5$, i.e., the last of the bases that compose the dictionary $\Psi_{p,l_0,L}$ (this is why the IGA stem plot is
    shifted rightward). For the standard IGA-Galerkin, we need to activate all the coefficients;
    conversely, the redundancy of the spline dictionary significantly promotes sparsity, so that we only need $s=9$ coefficients
    to recover a decent approximation of the solution by solving the PG system with OMP:
    5 for the coarse part of the solution and 4 for the localized feature.}
  \label{fig:IGAvsOMP}
\end{figure}

Since we are assuming $u$ to be well approximated by a sparse function $\widetilde{u}$ of the form \eqref{eq:def_u_tilde}, we look for an $s$-sparse approximate solution to \eqref{eq:PG_system}, obtained by solving
\begin{equation}\label{eq:sparse_recovery}
\min_{\bm{z} \in \mathbb{R}^{N_{\mathrm{dict}}}} \|B \bm{z} - \bm{c} \|_2 \quad \text{s.t.} \quad \|\bm{z}\|_0 \leq s,
\end{equation}
for a suitable small value of $s \in \mathbb{N}$ (such that $s \ll N_{\mathrm{dofs}} \leq N_{\mathrm{dict}}$) chosen by the user.
This problem is actually NP-hard \cite{natarajan1995sparse} but it can be approximately solved by sparse recovery approaches
such as Orthogonal Matching Pursuit (OMP) (see, e.g., \cite[Section 3.2]{foucart2013mathematical}).
Of course, OMP is not the only option to compute sparse solutions to \eqref{eq:PG_system}.
Other choices include $\ell^1$ minimization and thresholding algorithms (see \cite[Section 3]{foucart2013mathematical}).
In this context, we choose OMP thanks to its ability to easily control the number of iterations given an estimate $s$
of the sparsity level and its computational efficiency for small values of $s$
(see \cite[Section 5]{brugiapaglia2015compressed} for a numerical comparison between $\ell^1$ minimization and OMP for sparse numerical approximation of PDEs).
As an example, Figure \ref{fig:IGAvsOMP} clearly shows the effectiveness of using a multilevel dictionary
in a PG setting and then resorting to OMP to compute a sparse approximate solution to the corresponding linear system. 

\subsection{\cossiga}\label{subsection:cossiga}

The final step is to reduce the dimensionality of the linear system \eqref{eq:PG_system} via randomized subsampling.
In other words, we aim at computing a sparse approximation $\widetilde{u}$ to $u$ of the form \eqref{eq:def_u_tilde}
without assembling the full PG matrix $B$ (that is in general densely populated),
but only a small submatrix of it composed by a randomized selection of its rows.
This is possible thanks to the choice of the trial and test functions (localized in space and frequency, respectively) and to the
previously mentioned uncertainty principle~\cite{donoho1989uncertainty}.

To begin with, we draw $m \ll N_{\mathrm{dof}}$ test indices $\tau_1, \ldots, \tau_m \in [N_{\mathrm{test}}]$ i.i.d.\ at random according to
a suitable discrete probability distribution $\bm{\pi} \in \mathbb{R}^{N_{\mathrm{test}}}$ over $[N_{\mathrm{test}}]$, i.e.\
$$
\mathbb{P}\{\tau_i = q \} = \pi_q, \quad \forall \, q \in [N_{\mathrm{test}}],\; \forall \, i \in [m].
$$
Next, we consider the $m \times N_{\mathrm{dict}}$ \cossiga discretization
\begin{equation}
\label{eq:Ax=b}
A \bm{z} = \bm{b},
\end{equation}
where 
\begin{equation}
\label{eq:def_A_b}
A_{ij} := a(\psi_j,\varphi_{\tau_i}), 
\quad b_i := \int_{\Omega} f \varphi_{\tau_i}, \quad \forall \, j \in [N_{\mathrm{dict}}], \; \forall \, i \in [m].
\end{equation}
% A sparse solution to the highly underdetermined  be now computed is found by sparse recovery techniques. In particular,
The \cossiga solution is then found by using OMP in order to compute an approximate solution $\widetilde{\bm{x}}$ to
\begin{equation}
\label{eq:sparse_recovery_cossiga}
\min_{\bm{z} \in \mathbb{R}^{N_{\mathrm{dict}}}} \|E(A \bm{z} - \bm{b})\|_2 \text{ s.t. } \|\bm{z}\|_0 \leq s.
\end{equation}
The diagonal scaling $E \in \mathbb{R}^{m \times m}$ is defined as
\begin{equation}
\label{eq:def_E}
E_{ik} = \frac{\delta_{ik} }{\sqrt{m\pi_{\tau_i}}}, \quad \forall \, i,k \in [m],
\end{equation}
and accounts for the effect of the nonuniform sampling and it is chosen such that $\mathbb{E}[(EA)^*(EA)] = B^*B$, where $\mathbb{E}[\,\cdot\,]$ denotes the expected value (see \cite{krahmer2014stable, rauhut2012sparse}).
% We note in passing that OMP is not the only option to compute sparse solutions to \eqref{eq:def_A_b}. Other choices include $\ell^2$ minimization and thresholding algorithms (see \cite[Section 3]{foucart2013mathematical}). In this context, we choose OMP thanks to its ability to easily control the number of iterations given an estimate $s$ of the sparsity level and thanks to its computational efficiency for small values of $s$ (see \cite{brugiapaglia2015compressed} for  a numerical comparison between $\ell^1$ minimization and OMP for sparse numerical approximation of PDEs)
A very important quantity in this context is the so-called subsampling rate, i.e, the ratio $m/N_{\mathrm{dof}}$:
  a successful application of the \cossiga method will deliver a good approximation of the true solution with a very small
  subsampling rate, i.e. with $m \ll N_{\mathrm{dof}}$. 

Of course, the choice of the sampling probability distribution $\bm{\pi}$ is crucial for the effectiveness of the method.
Following ideas from \cite{krahmer2014stable,rauhut2012sparse} and the theoretical recipe in \cite{brugiapaglia2018theoretical},
we define $\bm{\pi}$ as a normalized upper bound to the so-called local $a$-coherence $\bm{\mu}$ of $\Psi_{p,l_0,L}$ with respect to $\Phi_R$, which is defined as
$$
\mu_q := \max_{j \in[N_{\mathrm{dict}}]} (a(\psi_j,\varphi_q))^2, \quad q \in [N_{\mathrm{test}}]. 
$$
%If an upper bound $\bm{\nu}$ to the local $a$-coherence $\bm{\mu}$ is available, i.e.\ if
In practice, the exact local $a$-coherence $\bm{\mu}$ is replaced with a suitable upper bound $\bm{\nu}$. Namely, if  
$$
\mu_q \leq \nu_q, \quad \forall q \in [N_{\mathrm{test}}],
$$
then, we let 
\begin{equation}
\label{eq:prob_density}
\pi_q = \frac{\nu_q}{\|\bm{\nu}\|_1}, \quad \forall q \in [N_{\mathrm{test}}].
\end{equation}
In order to estimate $\bm{\nu}$, we employ the theoretical results in \cite{brugiapaglia2018wavelet}.
In particular, we employ the following upper bound, corresponding to \cite[Equation (4.27)]{brugiapaglia2018wavelet}
(simplified by observing that $\|\bm{r}\|_0 \leq d$):
\begin{equation}
\label{eq:def_nu}
\nu_{q(\bm{r})} = \min\left\{ \frac{2^{(3d-2)L}\|\bm{r}\|_2^2}{\prod_{k = 1}^d r_k^4} ,  \frac{\|\bm{r}\|_2^2}{\|\bm{r}\|_\infty^2 \prod_{k = 1}^d r_k} \right\},
\quad \forall \bm{r} \in[R]^d,
\end{equation}
where $q: [R]^d \to [N_{\mathrm{test}}]$ corresponds to the ordering on $[R]^d$ used for the test functions.

\begin{figure}[t]
  \centering
  \includegraphics[width=0.34\linewidth]{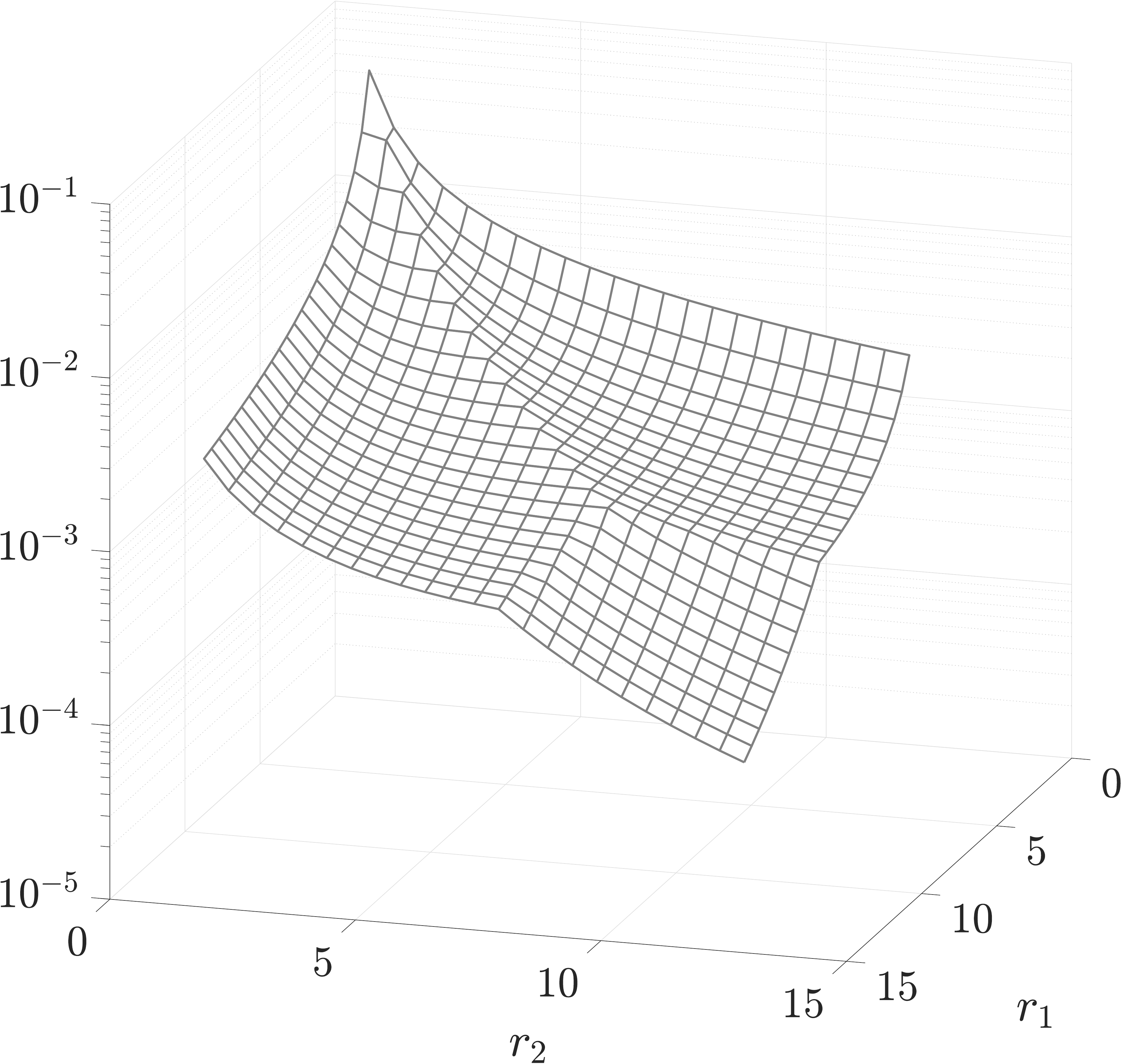}\hspace{0.5cm}
  \includegraphics[width=0.60\linewidth]{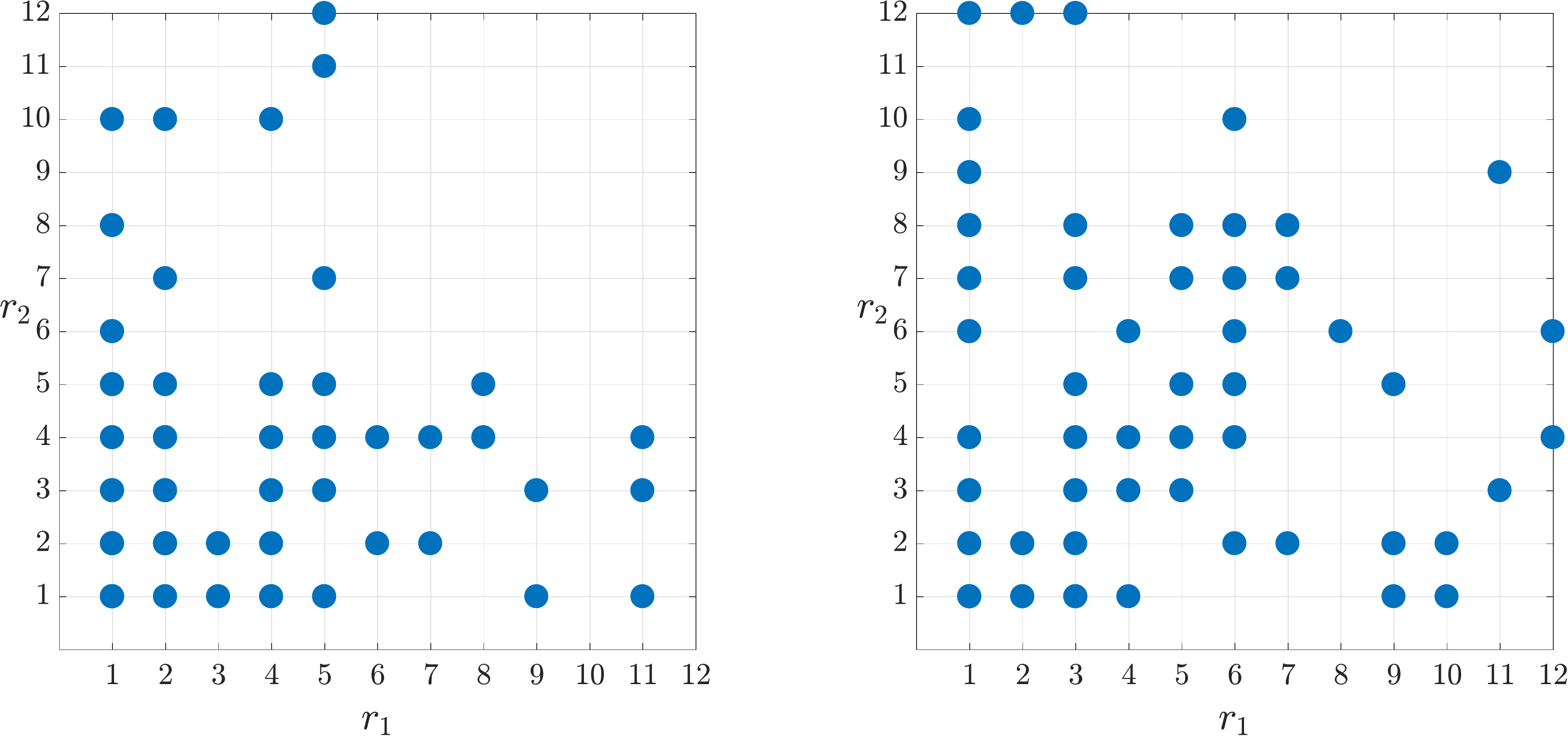}
  \caption{Left: sampling probability $\bm{\pi} = \bm{\nu}/\|\bm{\nu}\|_1$ in the 2D case.
      Center and right: two sets of frequencies $\bm{r} = (r_1,r_2)$ randomly sampled using such sampling probability with $m = 80$.
      Sampling is done with replacement in order to have i.i.d.\ samples (this is why the central and the right panels in the figure show less than 80 dots each).
      The reduced test space employed in \cossiga is the span of the sine functions \eqref{eq:sin_basis_fun_with_vector_index_rr}
      corresponding to the selected frequencies.\label{fig:sampling_rows_with_nu}}
\end{figure}
In Figure~\ref{fig:sampling_rows_with_nu}, we show the sampling probability $\bm{\pi} = \bm{\nu}/\|\bm{\nu}\|_1$ in the 2D case
  and two corresponding random samples of frequencies $\bm{r}$. The probability distribution employed selects lower frequencies with higher probability.
  This is in line with standard compressive sensing results, where sampling distributions concentrated
  on lower frequencies are known to recover multiscale coefficients of real-world signals (e.g., natural images)
  significantly better than the uniform distribution (see, e.g., \cite{BreakingBarrier,krahmer2014stable}).

\begin{remark}
The choice of the upper bound $\bm{\nu}$ made in \eqref{eq:def_nu} is not fully
justified from the theoretical perspective, but it has to be considered heuristic.
This is due to two main reasons:
(i) the set of trial functions $\Psi_{p,l_0,L}$ is a dictionary and not a Riesz basis, as
assumed in the theoretical framework of \cite{brugiapaglia2018wavelet}
and \cite{brugiapaglia2018theoretical};
(ii) the estimate \cite[Equation (56)]{brugiapaglia2018wavelet} used to derive \eqref{eq:def_nu} holds
in the case of B-spline wavelet (with $p=1$) tested against Fourier
functions defined over a periodic tensor product domain.
Consequently, the influence of the degree $p$ and the geometry of $\Omega$ are not
taken into account by \eqref{eq:def_nu}. Deriving a rigorous upper bound $\bm{\nu}$ to $\bm{\mu}$ for \cossiga is beyond the scope of
this paper and is left to future work.
%Figure \ref{fig:mu_vs_nu} shows the effectiveness of computing $\bm{\pi}$ based on $\bm{\nu}$ rather than on $\bm{\mu}$
%  in our context, considering two different geometries and $p=2$;
%  to read this plot, remember that the index $q$ is the lexicographic order given to the trigonometric basis functions $\sin_{\bm{r}}$ defined in \eqref{eq:sin_basis_fun_with_vector_index_rr} with respect to the frequency multi-index $\bm{r}$.
%  First of all, one can immediately see that in this case $\bm{\nu}$ is a reasonable estimate for $\bm{\mu}$
%  for values of $q$ corresponding to low frequencies $\bm{r}$,
%  while is a very pessimistic/conservative estimate for values of $q$ corresponding to high frequencies $\bm{r}$. We also note that,
%  although the trend is similar, the value of $\bm{\pi}$ based on true $a$-coherence shows a little dependence on the geometry,
%  while the bound (\ref{eq:def_nu}) does not. The fact that the bound is pessimistic for high frequencies $\bm{r}$ means that
%  we will sample more than necessary the rows of the matrix corresponding to those frequencies.
%  This can be seen in Figure \ref{fig:sampling_rows_with_mu_and_nu}, where we show three samples of frequencies obtained using
%  the  probability distribution $\mu_q/\|\bm{\mu}\|_1$ (top row) and three samples obtained using
%  $\nu_q/\|\bm{\nu}\|_1$ (bottom row). While both sampling methods concentrate on low frequencies, the samples obtained with 
%  $\nu_q/\|\bm{\nu}\|_1$ show a larger dispersion towards high frequencies.
\end{remark}

\begin{remark}
The convergence theory for CORSING proposed in \cite{brugiapaglia2018theoretical} does not require $\bm{\nu}$ to be a sharp upper bound to $\bm{\mu}$, but only an upper bound to $\bm{\mu}$ such that $\|\bm{\nu}\|_1 \ll N_{\mathrm{dof}}$. In fact, it can be shown that when both the trial and test functions are Riesz bases, drawing $m$ test functions with $m \geq c s \|\bm{\nu}\|_1 (s\ln(\mathrm{e}N_{\mathrm{dof}}/s) + \ln(2s))$ (where $c>0$ is a universal constant) using to the probability distribution  $\bm{\pi}= \bm{\nu}/\|\bm{\nu}\|_1$ is sufficient to achieve a recovery error proportional to the best $s$-term approximation error (see \cite[Theorem 3.15]{brugiapaglia2018theoretical}). 
Although this theory provides sufficient conditions for sparse recovery, there are not results about the optimality of this sampling strategy.
\end{remark}

\begin{figure}[t]
  \centering
  \includegraphics[width=0.24\linewidth]{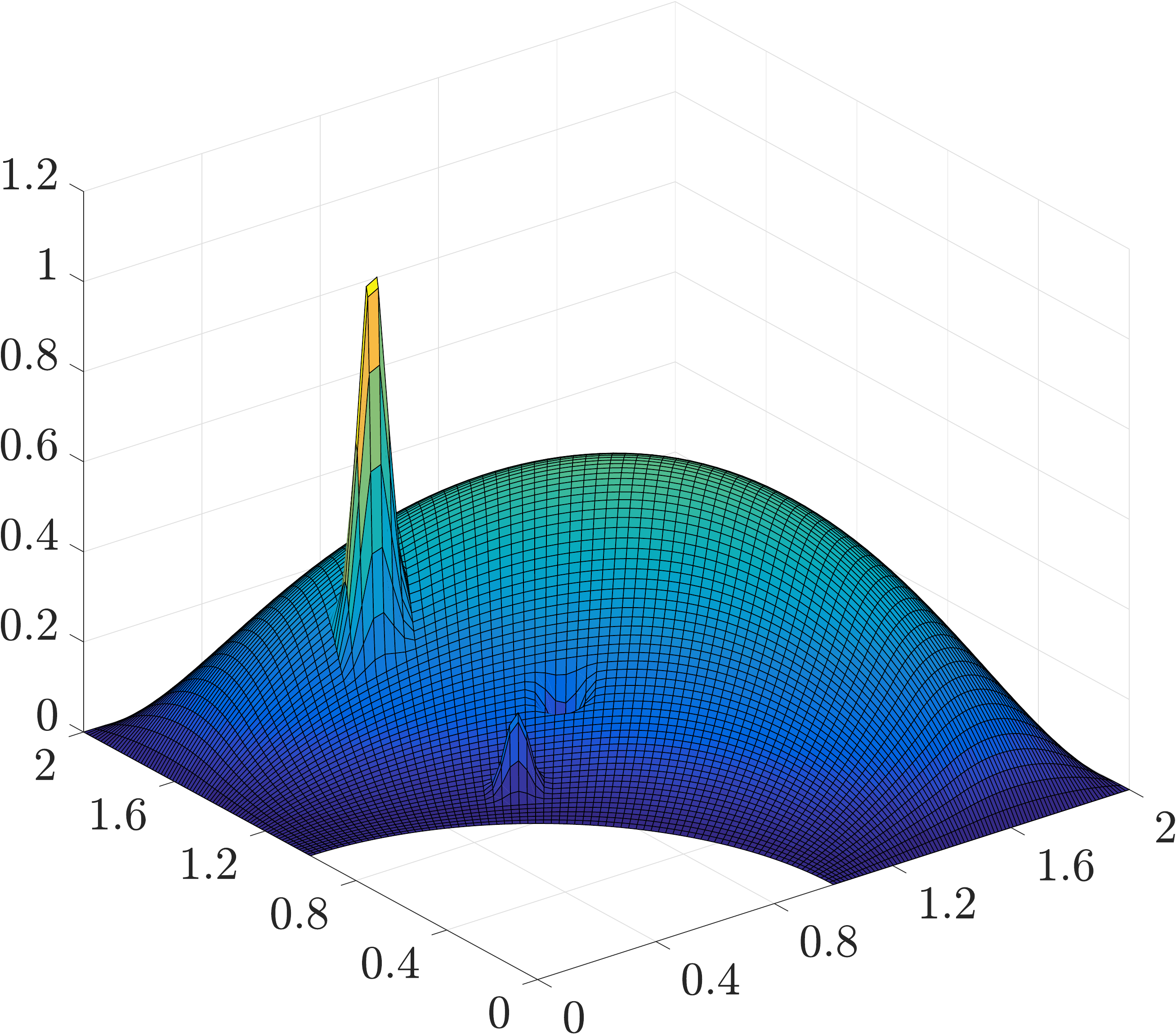}
  \includegraphics[width=0.24\linewidth]{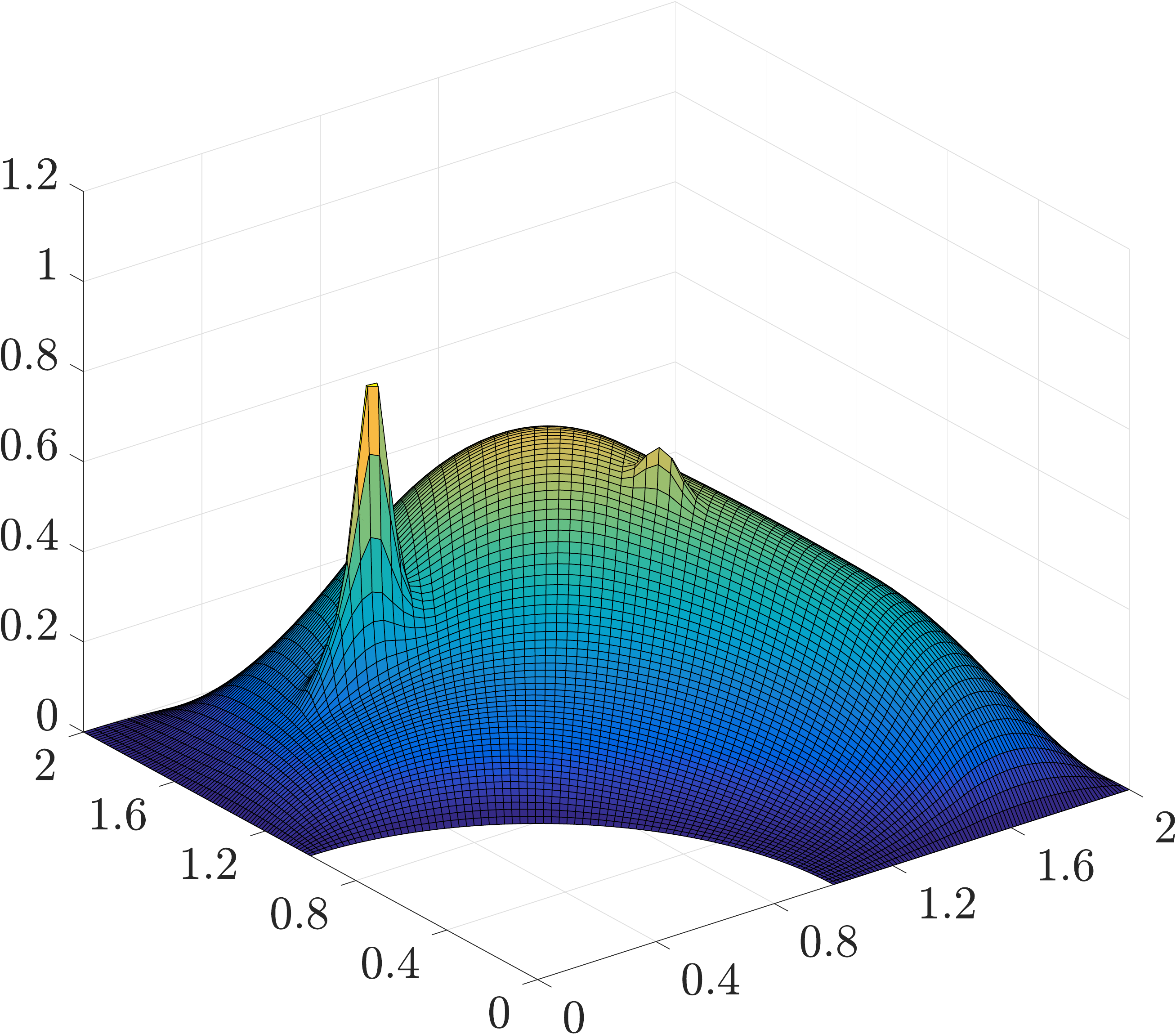} 
  \includegraphics[width=0.24\linewidth]{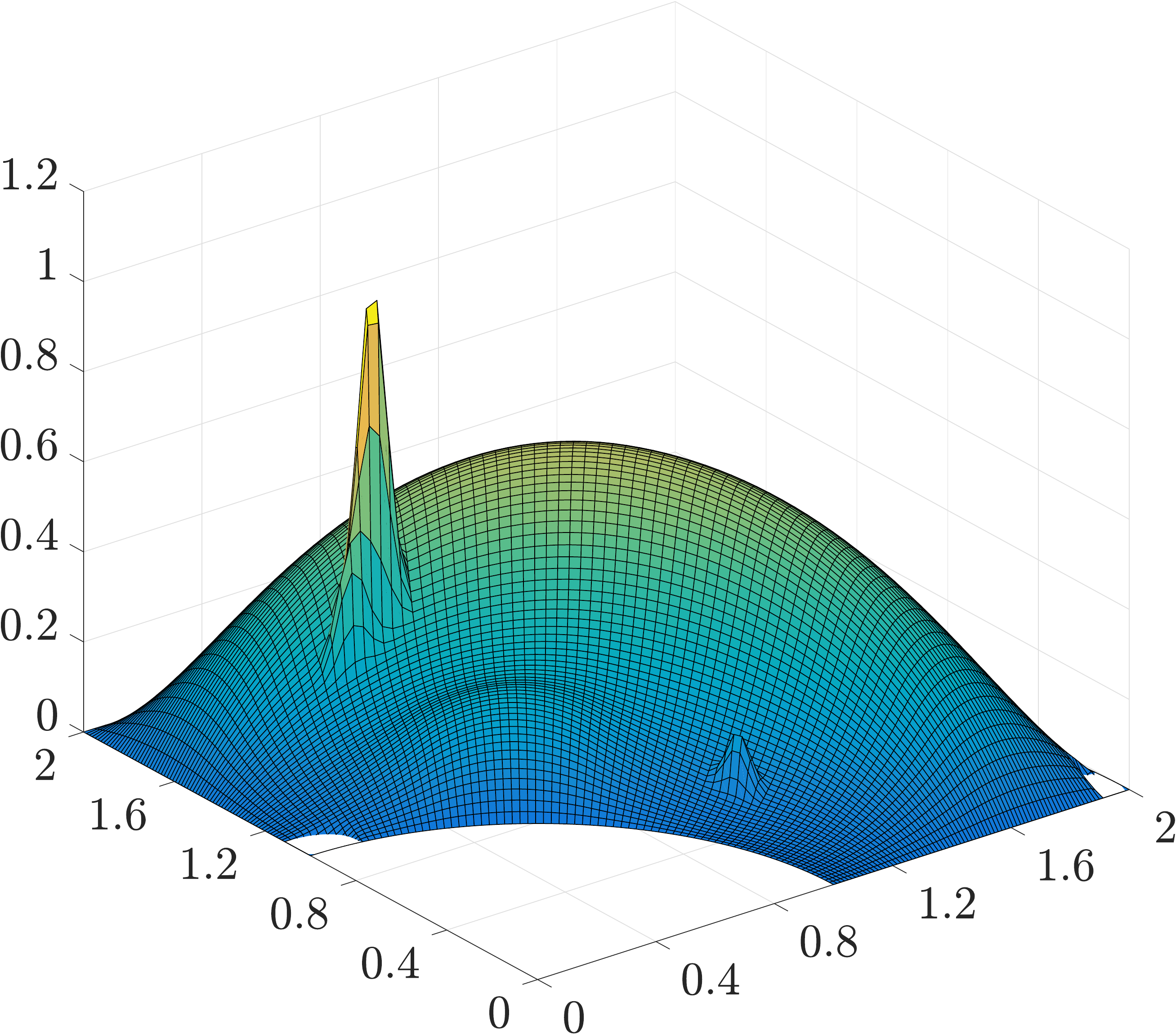} 
  \includegraphics[width=0.24\linewidth]{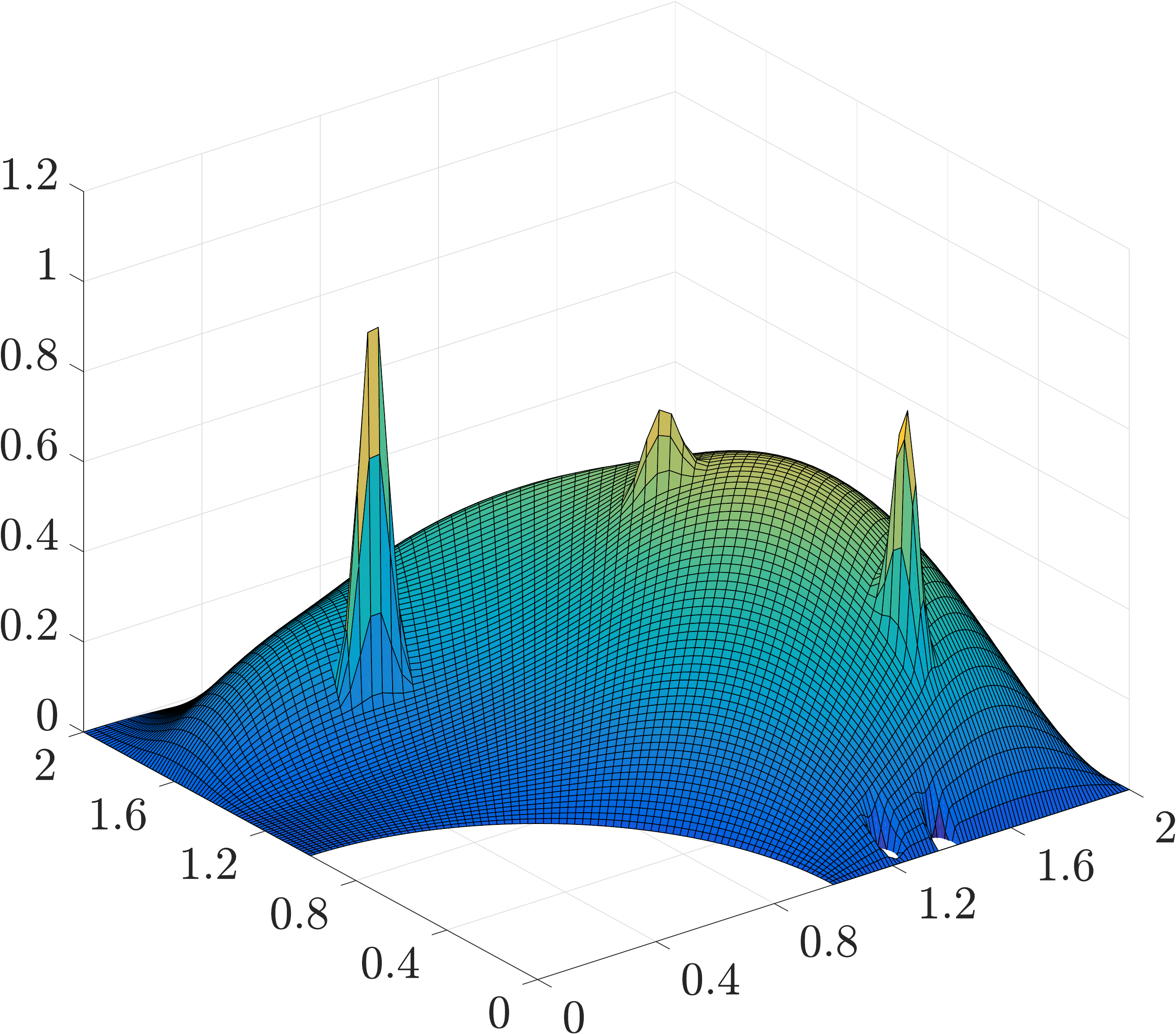} 
  \caption{Four different realizations of \cossiga for the problem of Figure 2, with $m=76$. The subsampling rate is $m/N_{\mathrm{dof}} = 76/1024 = 7.4\%$.}\label{fig:a_bunch_of_cossigas}
  \vspace{20pt}
  \includegraphics[width=0.24\linewidth]{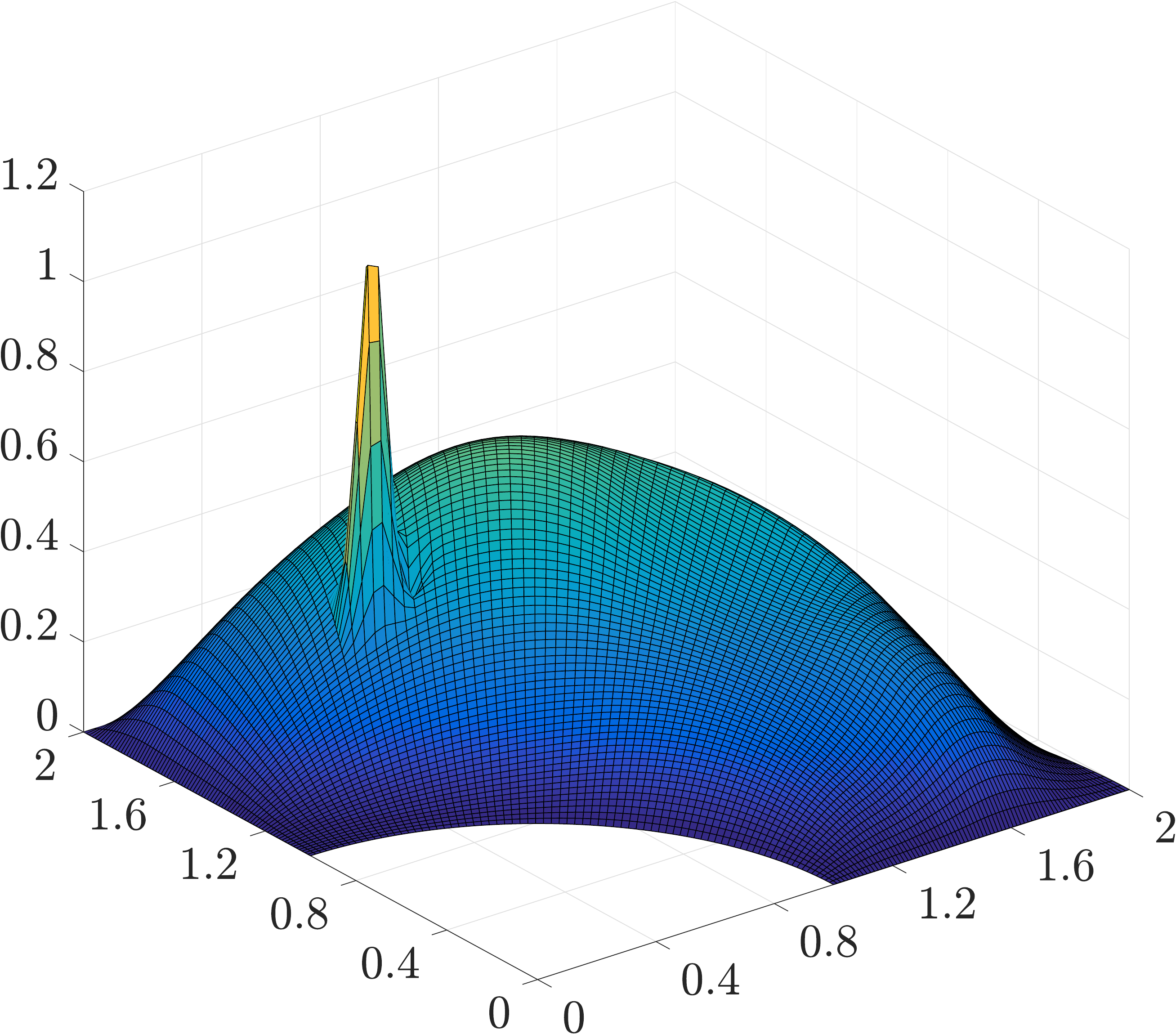} 
  \includegraphics[width=0.24\linewidth]{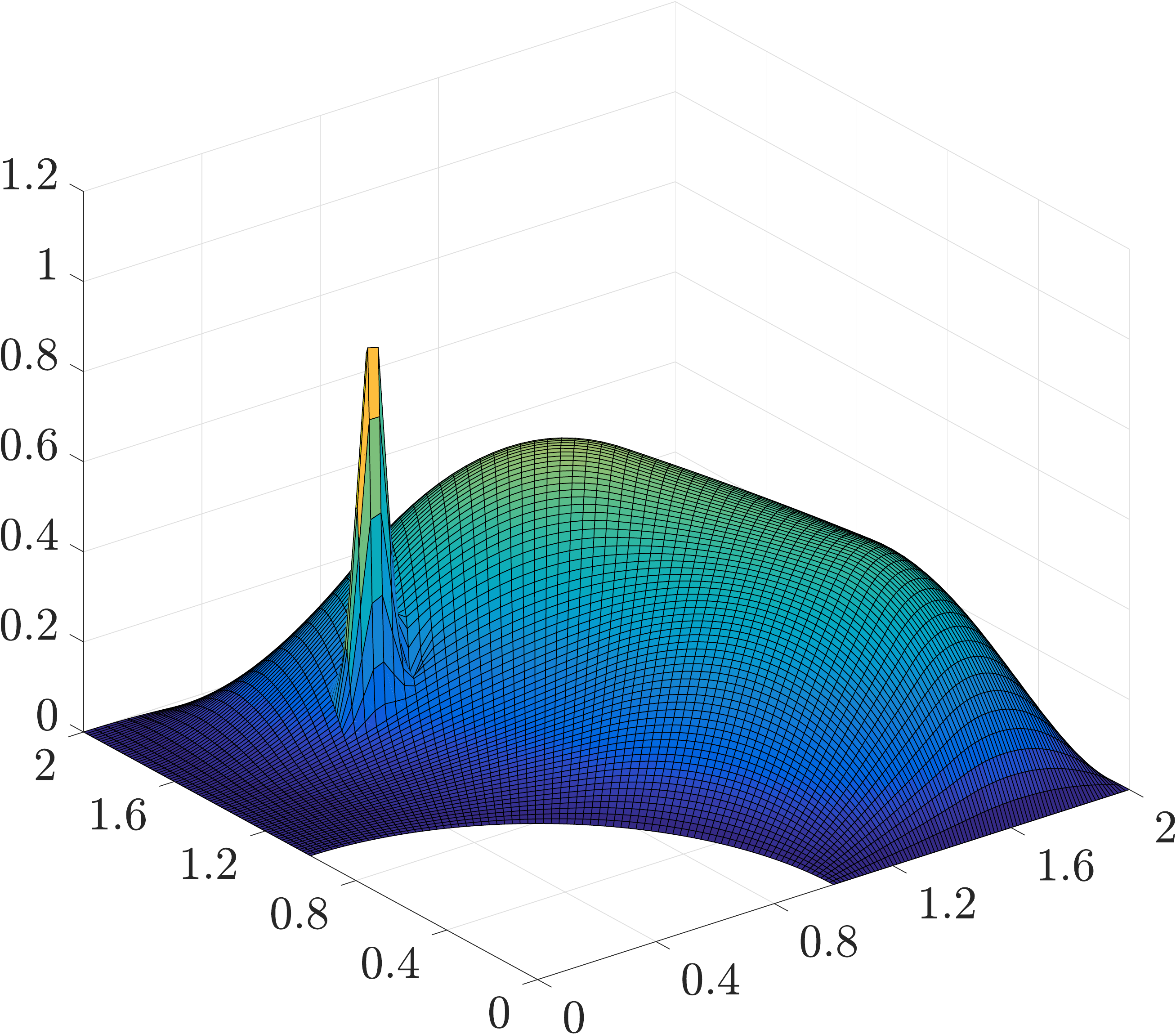}
  \includegraphics[width=0.24\linewidth]{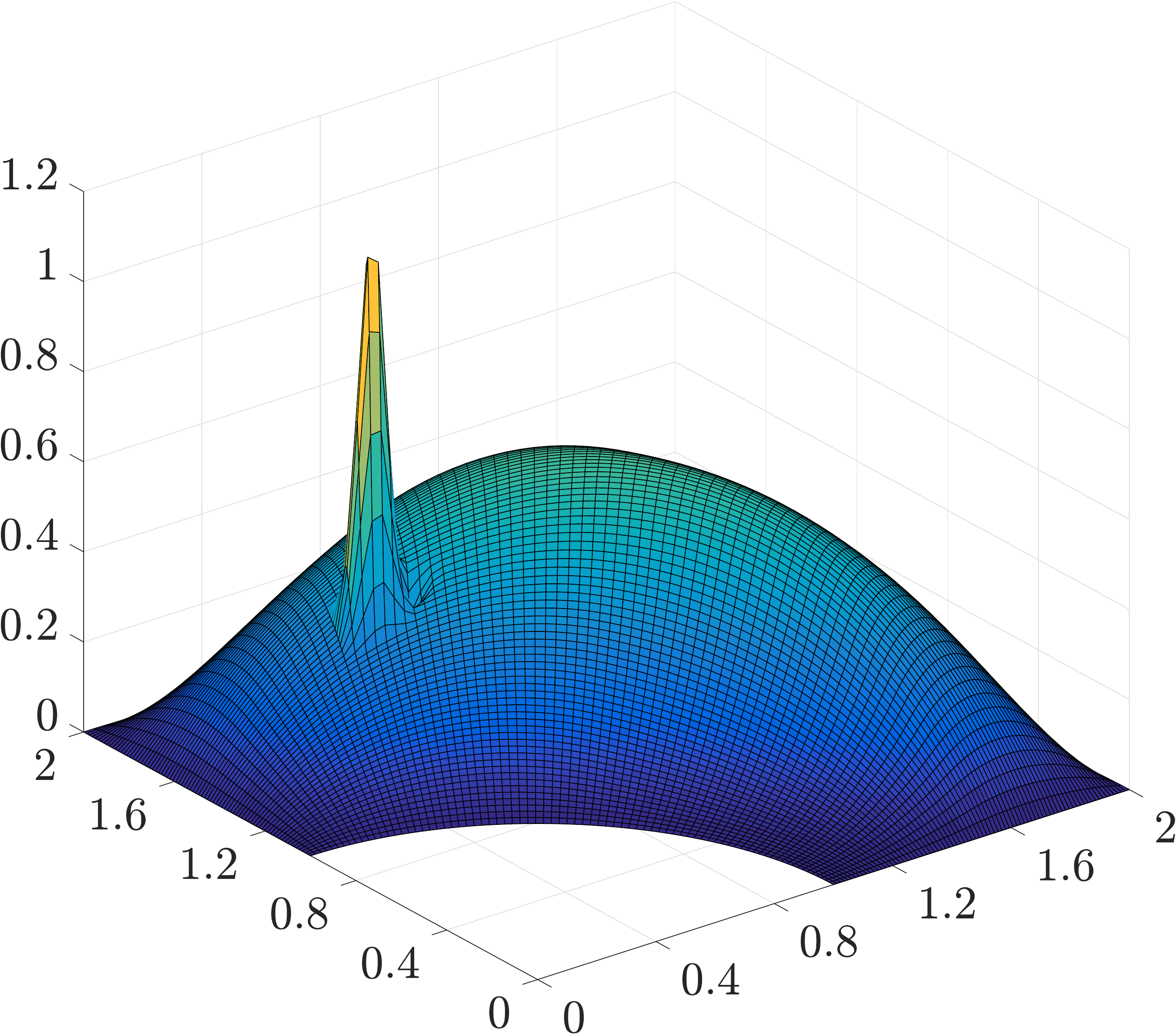} 
  \includegraphics[width=0.24\linewidth]{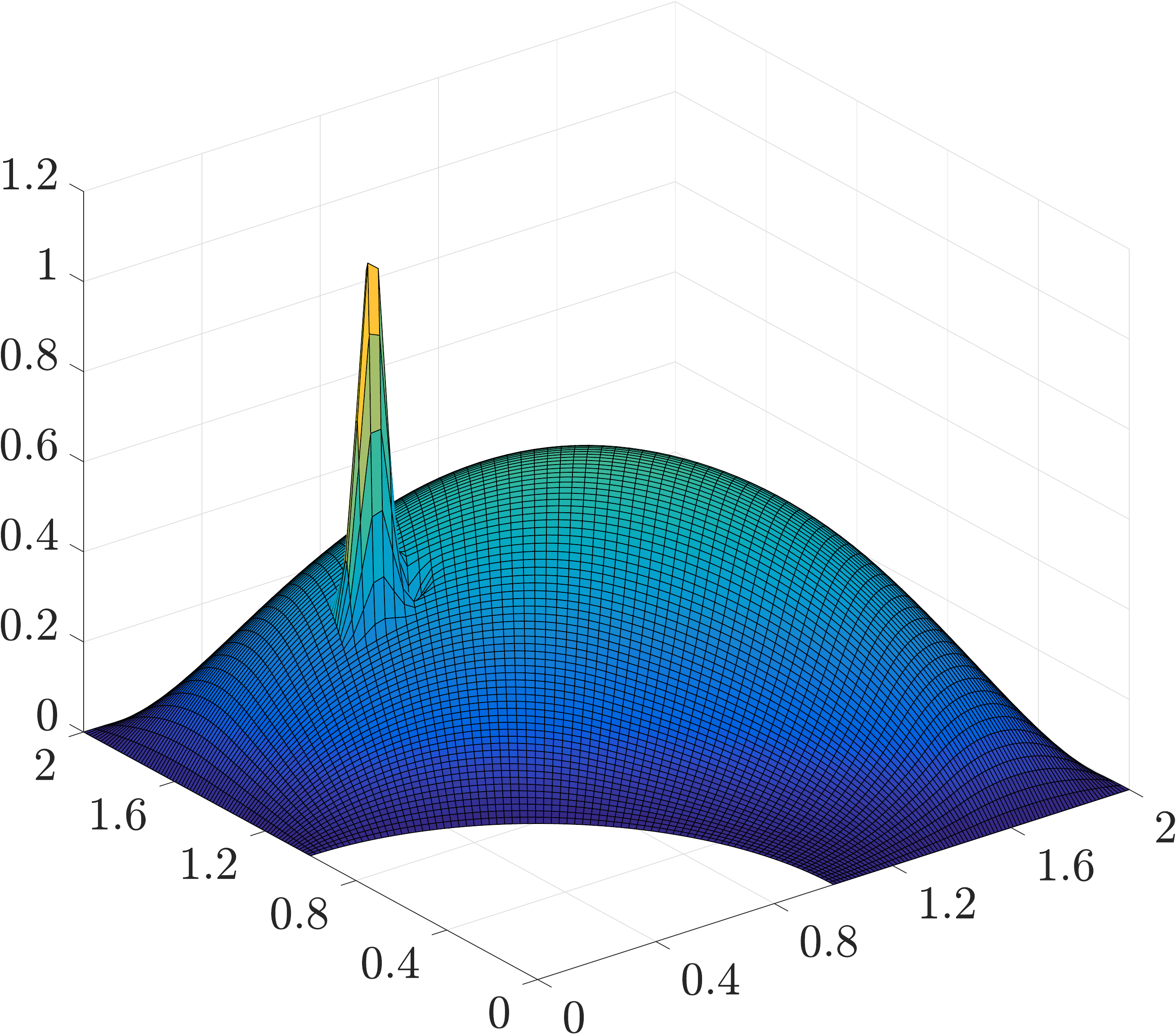} 
  \caption{Four different realizations of \cossiga for the problem of Figure 2, with $m=304$. The subsampling rate is $m/N_{\mathrm{dof}} = 304/1024 = 29.7\%$.}\label{fig:a_bunch_of_overkill_cossigas}
\end{figure}

The \cossiga approach is summarized in Algorithm \ref{alg:CossIGA}.
To convince the reader of the effectiveness of the method,
Figure \ref{fig:a_bunch_of_cossigas} shows four different realizations of CossIGA on the problem discussed in Figure 2,
obtained using a very small value of $m$, i.e. a very small subsampling rate.
Figure \ref{fig:a_bunch_of_overkill_cossigas} shows the same results, obtained using a larger value of $m$.
As expected, increasing $m$ improves the chances of a good recovery of the exact solution.

%Figure \ref{fig:cossiga_show} shows two different \cossiga solutions for different values of $m$, in the same setting as in Figure \ref{fig:IGAvsOMP}; as expected, the larger $m$ the better the chances of a good recovery of the true solution. Figures \ref{fig:a_bunch_of_cossigas} and \ref{fig:a_bunch_of_overkill_cossigas} show four additional realization of \cossiga solutions, to appreciate visually how the ``likelihood'' of a good recovery depends on $m$.

\begin{algorithm}
\caption{\label{alg:CossIGA} \cossiga (COmpreSSive IsoGeometric Analysis)}
%\normalsize
\textbf{Inputs:}
\begin{itemize}
\item $p$: B-spline degree.
\item $\mathrm{reg}$: the regularity of the B-splines;
\item $L$: maximum hierarchical level;
\item $s$: target sparsity level;
\item $m$: number of random test functions.
\end{itemize}
\textbf{Procedure:} 

\vspace{0.2cm}

$\widetilde{u} = \cossiga(p, \mathrm{reg}, L, s, m)$
\begin{enumerate}
\item Let $R$ as in \eqref{eq:choice_R} and $N_{\mathrm{test}} = R^d$.
\item Randomly draw $\tau_1,\ldots, \tau_m \in [N_{\mathrm{test}}]$ i.i.d.\ according to $\bm{\pi}$, defined as in \eqref{eq:prob_density}-\eqref{eq:def_nu}.
\item Build $A\in \mathbb{R}^{m \times N_{\mathrm{dict}}}$, $y \in \mathbb{R}^{m}$, and $E \in \mathbb{R}^{m \times m}$ defined as in \eqref{eq:def_A_b} and \eqref{eq:def_E}.
\item \label{alg:Cossiga:OMPstep} Compute $\widetilde{\bm{x}} \in \mathbb{R}^{N_{\mathrm{dict}}}$ by applying $s$ iterations of OMP to \eqref{eq:sparse_recovery_cossiga}. 
\item Let $\displaystyle\widetilde{u} = \sum_{j = 1}^{N_{\mathrm{dict}}} \widetilde{x}_j \psi_j$.
\end{enumerate}
\textbf{Output:}
\begin{itemize}
\item $\widetilde{u}$: $s$-sparse approximation of the solution $u$ to \eqref{eq:weak_problem}.
\end{itemize}

\end{algorithm}

\subsection{Practical setup for an effective use of \cossiga}\label{subsection:effective_CossIGA}

Algorithm~\eqref{alg:CossIGA} depends on five input parameters.
In analogy with classical IGA-Galerkin method,
we let the user choose the B-spline degree $p$, the regularity of the B-splines (either $p-1$ or $0$) 
and the maximum hierarchical level $L$ (or, equivalently, the mesh size $h = 2^{-L}$).
Therefore,  we are left with identifying two more parameters,
i.e., the target sparsity $s$ and the number of random test functions $m$.
Of course, in an ideal setting one would have at disposal some \emph{a priori} estimates that give an
indication on the optimal choices of $s$ and $m$. In this paper, these choices will instead be
made numerically, based on a calibration procedure. We leave the \emph{a priori} analysis for future work.
Let us explain the philosophy behind the numerical calibration, leaving the technical details to the next section.

We assume that for given values of $p$ and $L$, a good portion of the full accuracy
(say, e.g., no more than twice
the best approximation error of $u$ in the B-spline space $S_p^{\mathrm{int}}(\bm{\Xi}_L,\Omega)$)
can be reached using a certain sparsity value $s^*(p,L)$,
which we assume to linearly depend on $N_{\mathrm{dof}}$, i.e.
\begin{equation}\label{eq:s-star=C*N}
s^*(p,L) = C(p) N_{\mathrm{dof}}, 
\end{equation}
where $N_{\mathrm{dof}} = N_{\mathrm{dof}}(p,L)$ is defined as in \eqref{eq:NDOFs}. More generally, one might conjecture a nonlinear dependence between $s^*(p,L)$ and $N_{\mathrm{dof}}$, of the form $s^*(p,L) = C(p) N_{\mathrm{dof}}^{\alpha(p)}$ for some $\alpha(p) > 0$. We choose $\alpha(p)=1$ for the sake of simplicity.
Of course, we expect $s^*(p,L)$ and $C(p)$ to depend  heavily on the specific solution and,
in particular, on its compressibility with respect to the dictionary $\Psi_{p,l_0,L}$.

Furthermore, for each value of $s$ we need to identify the minimum value of $m$ such that \cossiga reaches a good portion of the accuracy corresponding to the best $s$-term approximation error of the solution with respect to the dictionary $\Psi_{p,l_0,L}$, defined as in \eqref{eq:best_s-term}. In principle, this $m$ might also depend on $L$ and $p$. Hence, we assume a dependence of the form
\begin{equation}
  \label{eq:m=D*s^beta}
m(p,L,s) = D(p,L) s. 
\end{equation}
The constant $D$ measures the compression capabilities of \cossiga and, contrary to the previous constant $C$,
  we expect a mild dependence of $D$ on the exact solution.
According to the compressive sensing theory \cite{foucart2013mathematical}, a sufficient condition to recover $s$-sparse vectors in $\mathbb{R}^N$ is $m \geq c \,s \polylog(N)$, where $c>0$ is a universal constant and $\polylog(N)$ is a polylogarithmic factor depending on the particular sampling scheme employed. These two factors are implicitly included in the constant $D(p,L)$ in \eqref{eq:m=D*s^beta}. Moreover, the constant $D(p,L)$ depends on the bilinear form associated with the PDE considered (see \cite{brugiapaglia2018theoretical}). For this parameter setting, the resulting subsampling rate of the  \cossiga method is
\begin{equation}
\label{eq:compression_ratio}
\frac{m(p,L,s^*(p,L))}{N_{\mathrm{dof}}} = C(p) D(p,L).
\end{equation}
Therefore, successful compression is achieved when $C(p) D(p,L) \ll 1$.

As already mentioned, determining \emph{a priori} the constants $C(p)$ and $D(p,L)$ deserves a careful theoretical
investigation and exceeds the scope of this paper. We will infer them from numerical experiments as discussed in the next section.

\section{Numerical tests}\label{section:tests}

\begin{figure}[t]
\centering
\includegraphics[width = 0.32\linewidth]{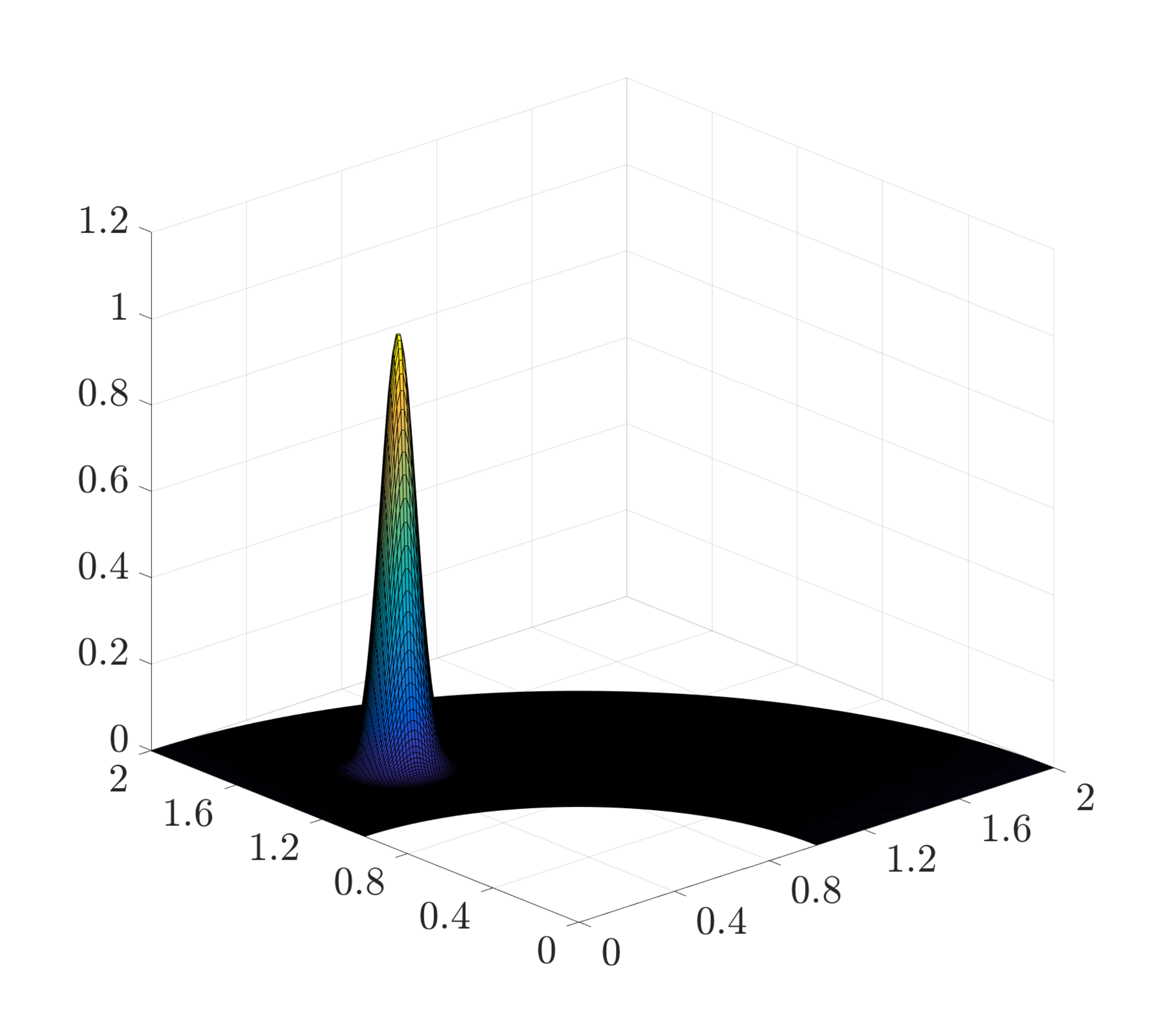}
\includegraphics[width = 0.32\linewidth]{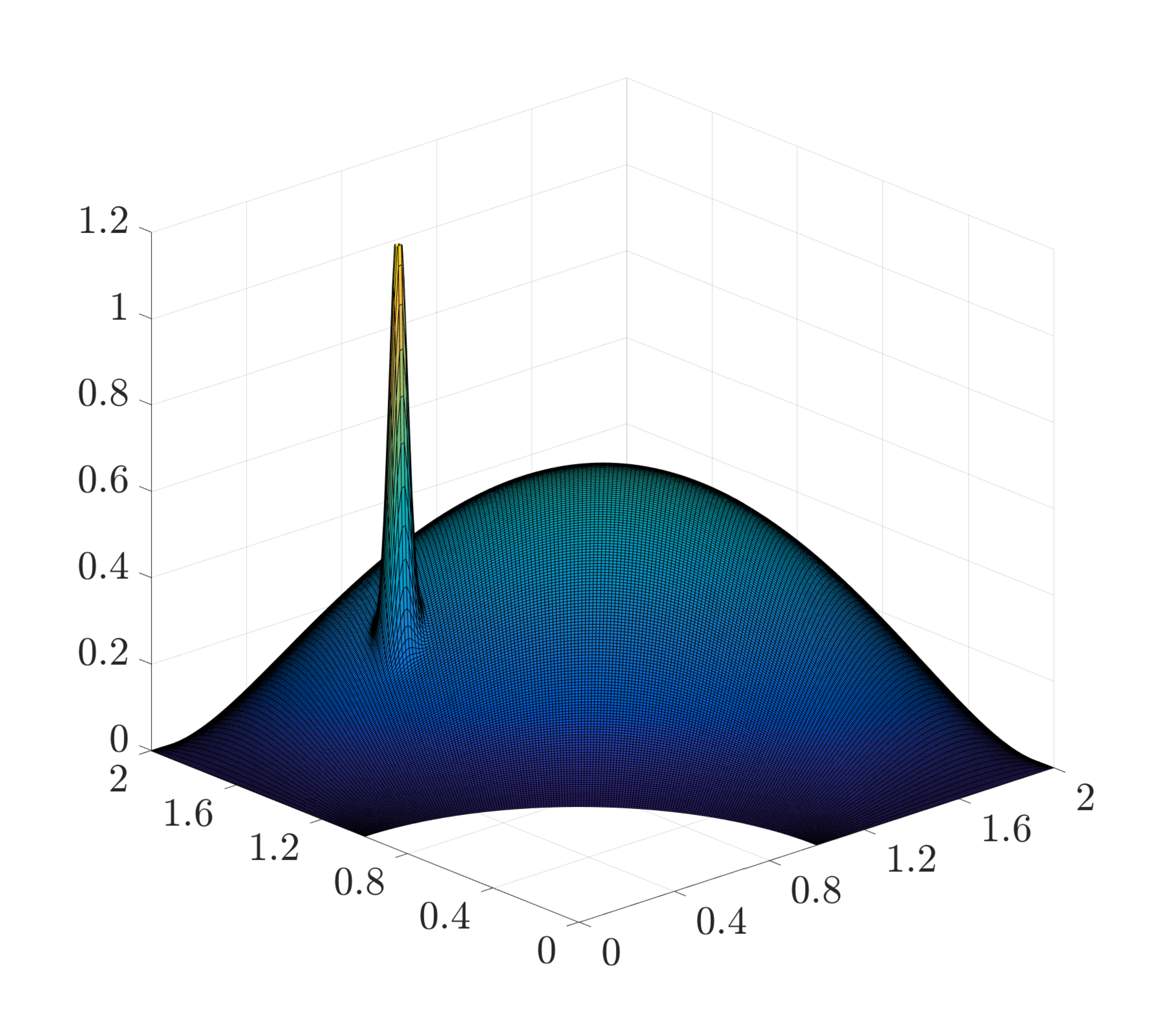}
\includegraphics[width = 0.32\linewidth]{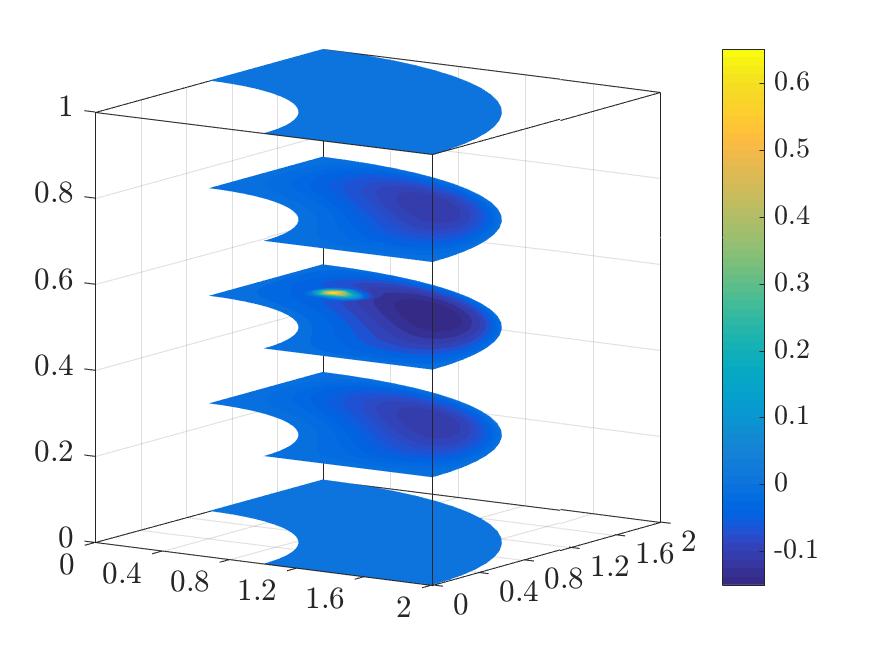}
\caption{Exact solution $u_{\mathrm{Gauss}}$ (left), $u_{\mathrm{polyGauss}}$ (center) and $u_{\mathrm{polyGauss\,3D}}$ 
  defined by \eqref{eq:u_Gauss}, \eqref{eq:u_polyGauss}, and (\ref{eq:u_polyGauss_3D})
  respectively, over the two-dimensional quarter of annulus $\Omega_{\mathrm{ring}}$ and the three-dimensional
  $\Omega_{\mathrm{thick\,ring}}$.} \label{figure:exact_sol}
\end{figure}

In the following, we test the performance of \cossiga for the numerical solution of \eqref{eq:weak_problem}
and for  four test cases. 

\paragraph{Case study I (Gauss 2D)} As a first physical domain, we consider the quarter of ring
$$
\Omega_{\mathrm{ring}} = \{(x_1,x_2): 1 \leq x_1^2 + x_2^2 \leq 4, \; x_1\geq0, \; x_2\geq 0\}.
$$
In order to study the recovery error of the method, we choose an analytical solution to \eqref{eq:weak_problem} defined over
$\Omega_{\mathrm{ring}}$. We call this case study ``Gauss'', corresponding to the exact solution 
\begin{equation}
\label{eq:u_Gauss}
u_{\mathrm{Gauss}}(x_1,x_2) = \exp\left(-\frac{(x_1-0.5)^2 + (x_2-1.4)^2}{(0.08)^2} \right).
\end{equation}
The function $u_{\mathrm{Gauss}}$ is very close to zero on most part of the domain and has a local
feature around $(0.5,1.4)$ (see Figure \ref{figure:exact_sol} (left)).
Although $u_{\mathrm{Gauss}}$ is not mathematically zero at the boundary, the homogeneous boundary conditions
are satisfied within the machine precision accuracy range.
In this case, the solution is expected to be sparse thanks to its small support in the physical domain.
We use this first case study also to
detail the calibration procedure proposed to estimate the constants $C$ and $D$
in \eqref{eq:s-star=C*N} and \eqref{eq:m=D*s^beta} respectively.
We then study the error of the \cossiga solution as a function of $L$
and for fixed $p$, when $s$ and $m$ are chosen according to \eqref{eq:s-star=C*N} and \eqref{eq:m=D*s^beta}
and using the values of $C$ and $D$ obtained via calibration. The accuracy achieved by \cossiga
is compared with the accuracy of the full solution of the PG system and with
the accuracy of the solution obtained by using OMP to approximate the solution to \eqref{eq:PG_system}
with $s^*$ coefficients, where $s^*$ is again chosen as in Equation (\ref{eq:s-star=C*N}).

\paragraph{Case study II (polyGauss 2D)} The second test case is called ``polyGauss'' and corresponds to the exact solution
\begin{equation}
\label{eq:u_polyGauss}
u_{\mathrm{polyGauss}}(x_1,x_2)
= \frac{1}{5}  x_1 x_2 (x_1^2 + x_2^2-1) (4-x_1^2-x_2^2)
+ \exp\left(-\frac{(x_1-0.5)^2 + (x_2-1.4)^2}{(0.04)^2} \right).
\end{equation}
In this case, the solution has a global support and a local feature and exhibits a multiscale behaviour
(see Figure~\ref{figure:exact_sol} (center)).
We use this test to investigate the sensitivity of \cossiga with respect to the parameters $C$ and $D$ in Section~\ref{subsection:polygauss}.
We achieve this goal by comparing the convergence results obtained using the constants
$C$ and $D$ calibrated on $u_{\mathrm{Gauss}}$ with the results of the same test performed after recalibrating $C$ and $D$ on $u_{\mathrm{polyGauss}}$.

\paragraph{Case study III ($C^0$ vs. $C^{p-1}$ splines 2D)} In the third test case, we consider again the solution $u_{\mathrm{polyGauss}}$ and 
investigate the impact of the smoothness of the spline dictionary considered,  i.e.,
whether it is advantageous to use $C^0$ splines instead of $C^{p-1}$ splines as in the first two
test cases (see Section~\ref{subsection:C0}). 

\paragraph{Case study IV (polyGauss 3D)} Finally, in the fourth test case we consider the three-dimensional generalization of the polyGauss
test case, i.e., we let 
\begin{align}
  & \Omega_{\mathrm{thick\,ring}} = \{(x_1,x_2,x_3): 1 \leq x_1^2 + x_2^2 \leq 4, \; x_1\geq0, \; x_2\geq 0,  \; 0 \leq x_3 \leq 1\}, \nonumber \\
  & u_{\mathrm{polyGauss\,3D}}(x_1,x_2,x_3)
    = \frac{1}{5}  x_1 x_2 x_3 (x_1^2 + x_2^2-1) (4-x_1^2-x_2^2) (x_3-1) \nonumber \\
  & \qquad \qquad \qquad  + \exp\left(-\frac{(x_1-0.5)^2 + (x_2-1.4)^2 + (x_3-0.5).^2}{(0.04)^2} \right), \label{eq:u_polyGauss_3D}
\end{align}
which is again zero within machine precision on the boundary of $\Omega_{\mathrm{thick\,ring}}$.
This function is displayed in Figure~\ref{figure:exact_sol}~(right).
For this problem, we perform the same test as in the first case study, i.e.\ we compare the different methods in terms of accuracy when $L$ varies.

\subsection{Case study I: Gauss 2D}\label{subsection:testgauss}

In this test, we consider B-splines of maximal regularity $C^{p-1}$. We will use the following terminology when referring to different approximate solutions to \eqref{eq:weak_problem}:
\begin{description}
%\item [IGA:] Solution computed via IGA Galerkin method with B-splines of degree $p$ and mesh size $h = 2^{-L}$. The corresponding approximation is denoted as $\tilde{u}_{\mathrm{IGA}}$.
\item [PG-OMP$(s)$:] Solution computed applying $s$ iterations of OMP to the PG system \eqref{eq:PG_system},
  i.e., to compute $s$ coefficients that approximately solve the minimization problem \eqref{eq:sparse_recovery}.
  The corresponding approximation is denoted as $\tilde{u}_{\mathrm{PG-OMP}(s)}$.
  The error associated with $\tilde{u}_{\mathrm{PG-OMP}(s)}$ is the best accuracy
  that we can expect from \cossiga$(p,L,s,m)$.
  In particular, in this approach there is no random compression of the test space,
  hence the recovery  error is only due to the PG approximation and to the $s$-sparse approximation computed via OMP.
  Note also that $\|\tilde{u}_{\mathrm{PG-OMP}(s)}-u\|_{H^1(\Omega)}$ is an upper bound to the best $s$-term
  approximation error of $u$ with respect to the dictionary $\Psi_{p,l_0,L}$ (defined in \eqref{eq:best_s-term}) and can be thought as a proxy for it.

\item [PG-BS:] Least-squares solution to the PG system $B \bm{z} = \bm{c}$,
 defined in \eqref{eq:def_B_c}.
  The corresponding approximation is denoted as $\tilde{u}_{\mathrm{PG-BS}}$. 
  Note that here the error of the solution is only due to the PG
  approximation, since we are neither compressing the test space nor
  sparsifying the solution. In particular, we expect the accuracy of
  PG-BS to be the best possible accuracy achievable by PG-OMP$(s)$ for any $s$ since $s \leq N_{\mathrm{dict}}$.
\end{description}

We also performed numerical tests with a standard IGA Galerkin discretization on the finest hierarchical level, and
the error obtained with this approach is almost identical to the one obtained using PG-BS.
This suggests that the empirical choice \eqref{eq:choice_R} is sufficient to achieve the discrete inf-sup stability condition \eqref{eq:discrete_infsup}.
Since also the number of degrees of freedom of the IGA discretization is comparable to the size of the full dictionary we use for PG-BS, only the results of the latter approach are shown in the following.
We begin the discussion of this test case by detailing the procedure to estimate the constants $C$ and $D$.

\subsubsection{$C$-calibration test}
\label{subsection:C-calibration}

For sake of explanation, we fix $p=2$ (we have tested also $p=1,4$ obtaining analogous results,
not shown for brevity). 
\begin{figure}[t]
\centering
\includegraphics[width = 0.45\linewidth]{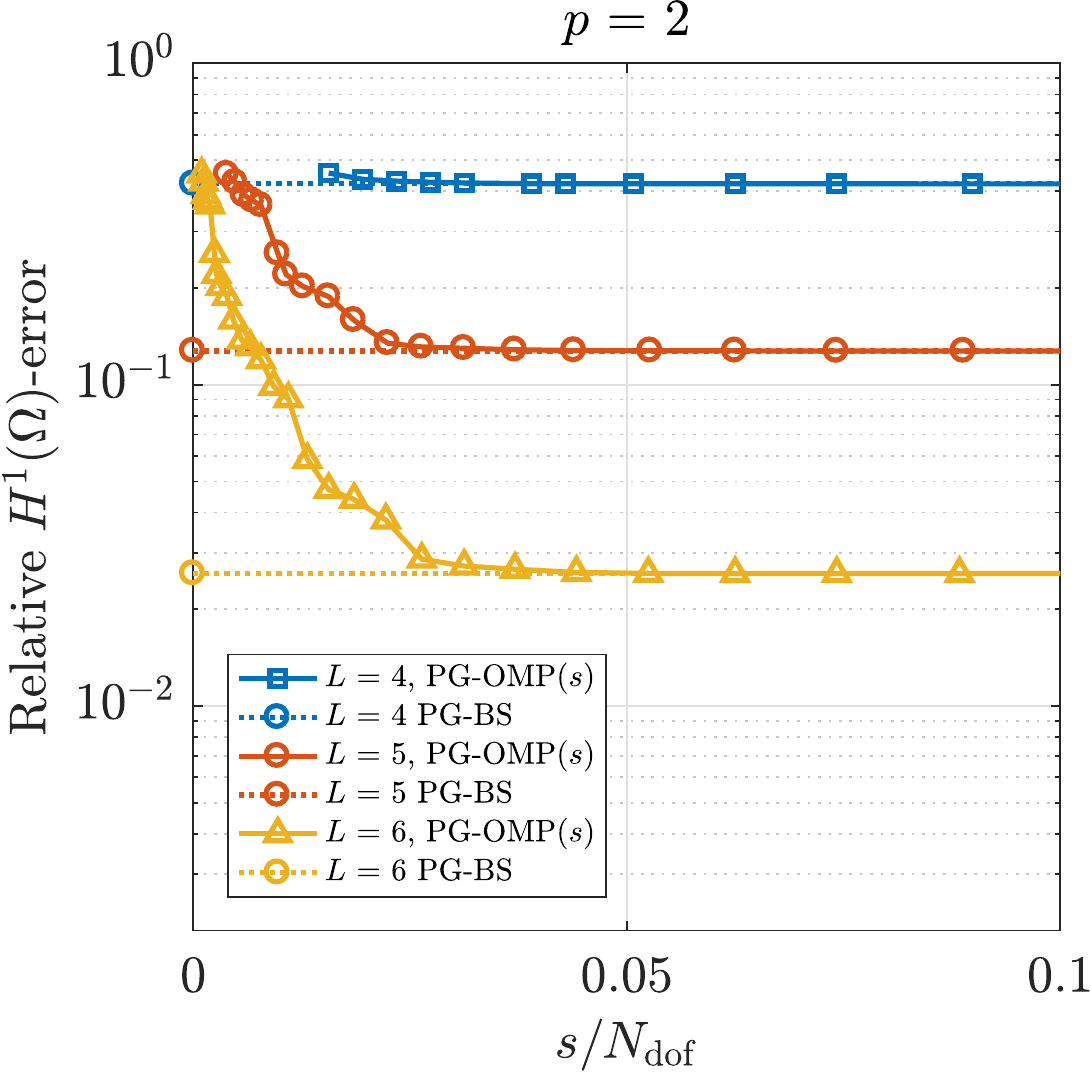}
\includegraphics[width = 0.45\linewidth]{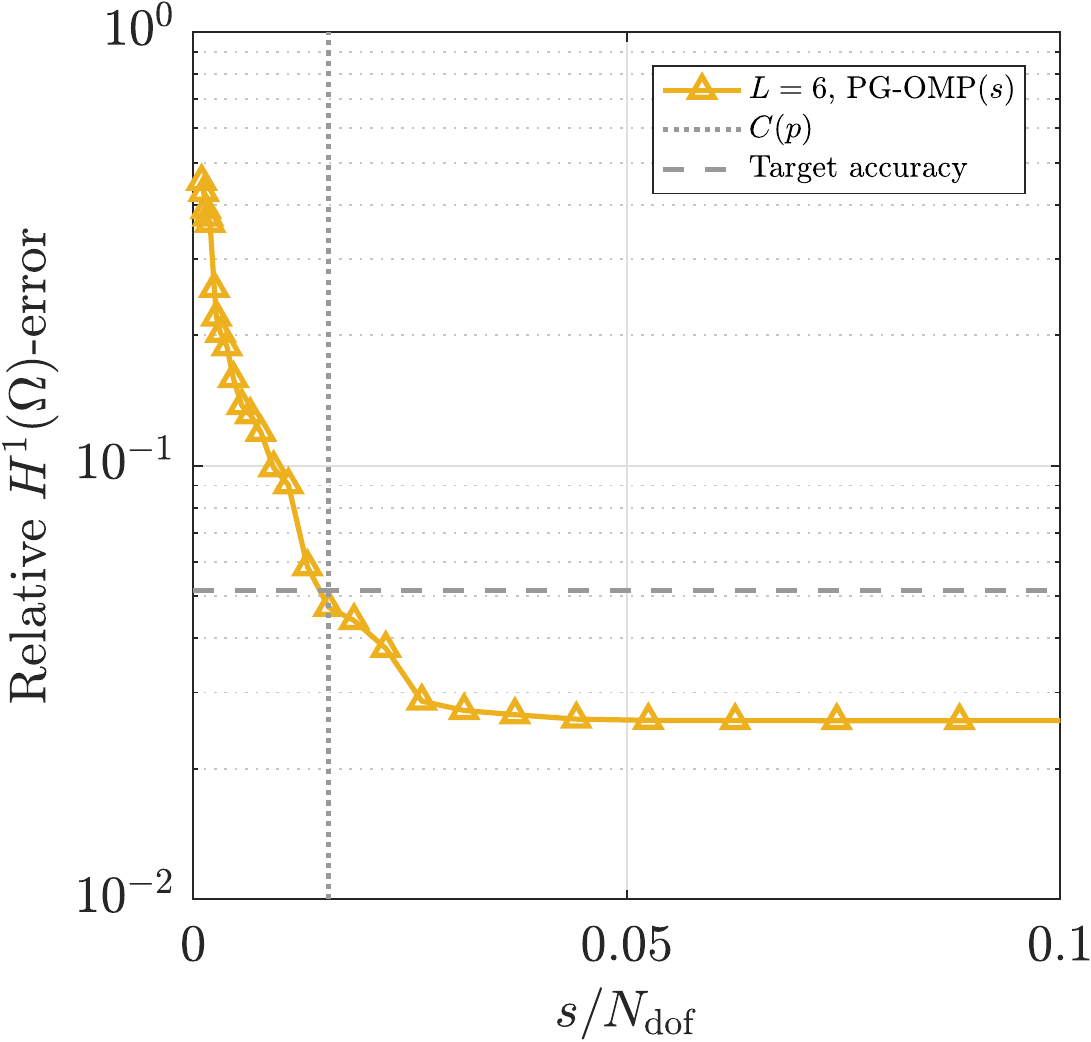}
\caption{Case study I (Gauss 2D). Left: error of $u_{\textsc{OMP}(h,p,s)}$ with respect to $u$ as a function of $s$, 
for increasing mesh refinement $L=4,5,6$ and degree $p=2$. Right: result of the $C$-calibration process for $L=6$.} \label{figure:C_calibration}
\end{figure}
In Figure~\ref{figure:C_calibration} (left), we show the relative $H^1(\Omega)$-error computed with PG-OMP$(s)$ for increasing values of $s$
  (normalized as $s/N_{\mathrm{dof}}$), which, as already discussed, is a proxy (upper bound) to the best $s$-term
  approximation error of $u$ with respect to the dictionary $\Psi_{p,l_0,L}$. We
observe that PG-OMP$(s)$ quickly reaches a plateau for every choice of $L$,
at the level of the PG-BS error. This means that the same accuracy of PG-BS
(or standard IGA Galerkin, as mentioned above) can be reached by activating only a small portion of the available coefficients in the dictionary,
i.e., the function $u_{\mathrm{Gauss}}$ is compressible in $\Psi_{p,l_0,L}$, as claimed before.
For a fixed accuracy, increasing the level $L$ leads to a reduction in 
the sparsity percentage needed to reach such target accuracy. However, the approximation
is globally less sparse for large values of $L$, i.e.\ it takes a larger
percentage of dofs to reach the accuracy plateau.
This means, in particular, that our assumption in Equation (\ref{eq:s-star=C*N}) only holds approximately (it would hold exactly if the ``elbows'' of the three convergence lines occurred at the same abscissa). Yet, in order to keep the complexity of our numerical illustration moderate,
we do not want to complicate the model (\ref{eq:s-star=C*N}) and we choose a conservative value for $C(p)$
(i.e., the one for the largest value of $L$ tested).

Let us now illustrate in detail how to perform the $C$-calibration test (see Figure~\ref{figure:C_calibration} (right)).
Our goal is to choose a suitable sparsity level $s^* = s^*(p,L)$ such that the relative error achieved via PG-OMP$(s^*)$
is comparable with the relative error of PG-BS, i.e.,
$\| u - \tilde{u}_{\textsc{OMP}(s)}\|_{H^1(\Omega)} \approx \mu \cdot \|u - \tilde{u}_{\mathrm{BS}}\|_{H^1(\Omega)}$,
where $\mu$ is a small constant larger than $1$; in particular, we choose $\mu=2$.
We then  look for the value of $s$ among those tested
that renders the error $\| u - \tilde{u}_{\textsc{OMP}(s)}\|_{H^1(\Omega)}$
as close as possible to $\mu \cdot \|u - \tilde{u}_{\mathrm{BS}}\|_{H^1(\Omega)}$. More precisely, we numerically compute $s^*(p,L)$ as
\begin{equation}
\label{eq:s^*_choice}
s^*(p,L) = \arg  \min_{s \in S_{\mathrm{tested}}}
\left| \| u - \tilde{u}_{\textsc{OMP}(s)}\|_{H^1(\Omega)} -  \mu \cdot \|u - \tilde{u}_{\mathrm{BS}}\|_{H^1(\Omega)}\right |, \;
S_{\mathrm{tested}} = \lceil 2^{[2:0.25:11]} \rceil, \; \mu = 2,
\end{equation}
where we used Matlab notation to define $S_{\mathrm{tested}}$.
 Once $s^*(p,L)$ is computed, one simply has $C(p) = s^*(p,L)/N_{\mathrm{dof}}$. 
The $C$-calibration process is repeated for multiple values of $p$, leading to the values of $C(p)$ given in Table~\ref{table:C-calibration_Gauss}.
\begin{table}[t]
\centering
\begin{tabular}{c|ccc}
$p$ & 1 & 2 & 4 \\%& 5 \\
\hline
$C(p)$ &  $8.0 \cdot 10^{-3}$ &  $1.6\cdot 10^{-2}$ & $3.5\cdot 10^{-2}$ %& $8.0\cdot 10^{-2}$
\end{tabular}
\caption{\label{table:C-calibration_Gauss}Case study I (Gauss 2D). Numerical estimate of $C(p)$ computed via the $C$-calibration test.}
\end{table}
In Figure~\ref{figure:C_calibration} (right), the value of $C(p)$ with $p=2$ (corresponding to the ratio $s^*/N_{\mathrm{dof}}$) is marked with a vertical dashed line
and the target accuracy 
\begin{equation}
\label{eq:target_accuracy_Gauss}
\mu \cdot \frac{\|u - \tilde{u}_{\textsc{BS}}\|_{H^1(\Omega)}}{\|u\|_{H^1(\Omega)}}, \qquad \mu = 2,
\end{equation}
with a horizontal dashed line.
We note that $C(p)$ is monotonically increasing with respect to $p$. 
%This can be perhaps explained by the fact
%that achieving the best accuracy is more and more challenging as $p$ increases. 

%\begin{remark}[Choice of $R$]
%\label{rem:choice_R}
%We note in passing that the accuracy of IGA and PG-BS are practically identical in Figure~\ref{figure:C_calibration}.
%This suggests that the empirical choice \eqref{eq:choice_R} is sufficient to achieve the discrete inf-sup stability condition \eqref{eq:discrete_infsup}.
%\end{remark}

% \begin{figure}[t]
% \centering
% \includegraphics[width = 0.7\linewidth]{figures/plain_IGA_error}
% \caption{error of $u_{IGA}$ as a function of $L$ for degree $p=1\ldots,6$.} \label{figure:plain_IGA}
% \end{figure}

% The plots shown so far focused on the sparsity properties of the Petrov-Galerkin
% solution only. We now present a sequence of tests in which we analyze the 
% performance of \cossiga, which is influenced by the parameters $L$, $s$, $m$ and $p$.
% Analyzing how the interplay among these parameters affects the quality of the 
% \cossiga solution is the goal of the next sections.

\subsubsection{$D$-calibration test}
\label{subsection:D-calibration}

This experiment aims at estimating the constant $D(p,L)$ in Equation (\ref{eq:m=D*s^beta}).
We fix the degree $p$ and the maximum hierarchical
level $L$;
in particular, similarly to the previous experiments, we consider $p=1,2,4$ and $L=4,5,6$.
For each combination of the parameters $L$ and $p$ we further consider different values of
$s$, with $2 \leq s \leq s^*(p,L)$, $s^*(p,L)$ obtained via \eqref{eq:s-star=C*N},
where $C(p)$ estimated by the $C$-calibration procedure just detailed.
For each value of $s$, we then let $m$ vary in the interval
$s \leq m \leq N_{\mathrm{dof}}$ and we perform $n_{\mathrm{runs}}=25$ random runs of
\cossiga for each combination of $p$, $L$, $s$, and $m$.
In particular, we let $m\geq s$ because sampling less than $m=s$ rows
does not allow to compute an $s$-sparse approximation via OMP and we
choose $m \leq N_{\mathrm{dof}}$ because we want to achieve
compression.

\begin{figure}[t]
\centering
\includegraphics[width = 0.374\linewidth]{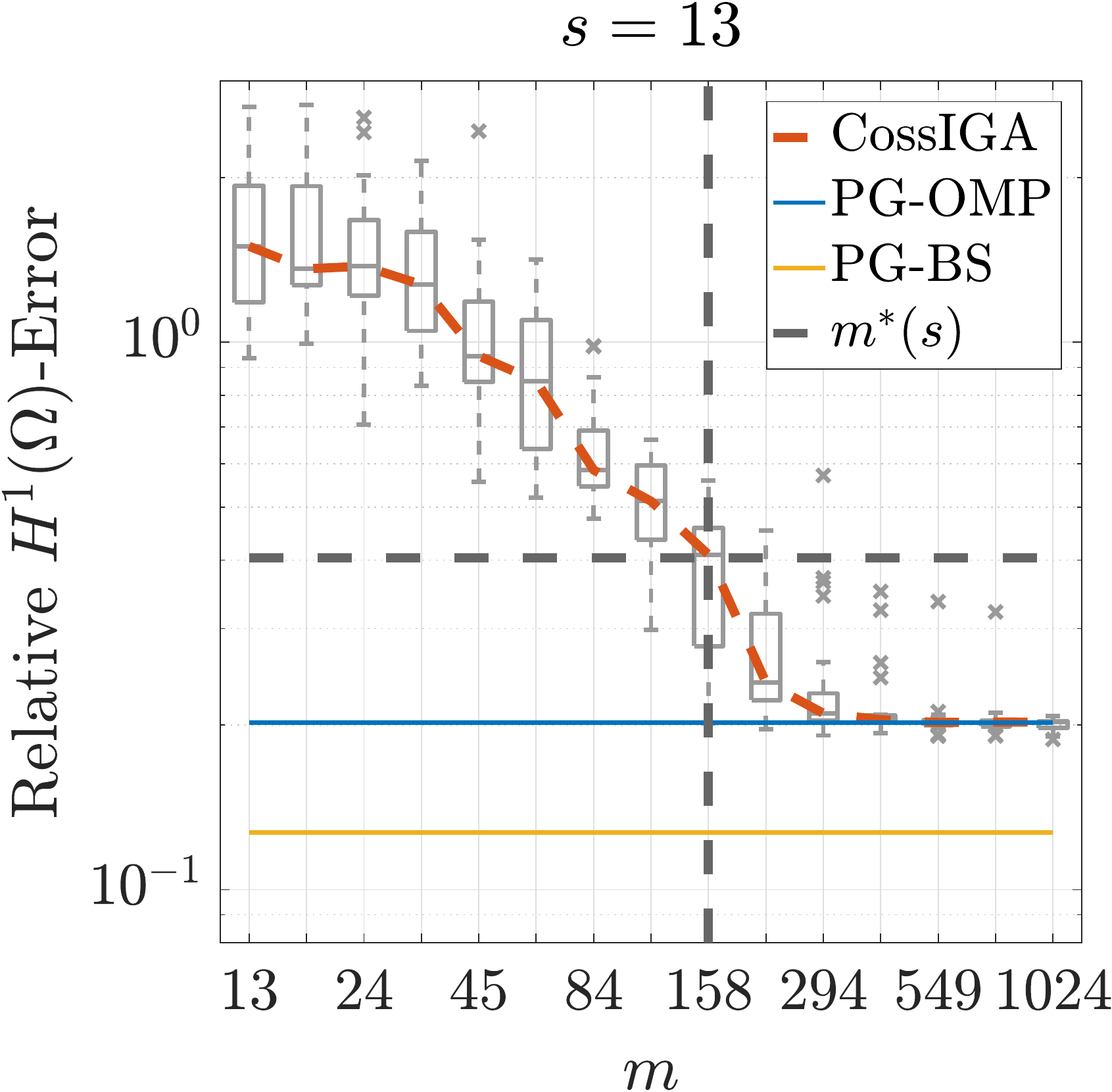} 
\includegraphics[width = 0.3\linewidth]{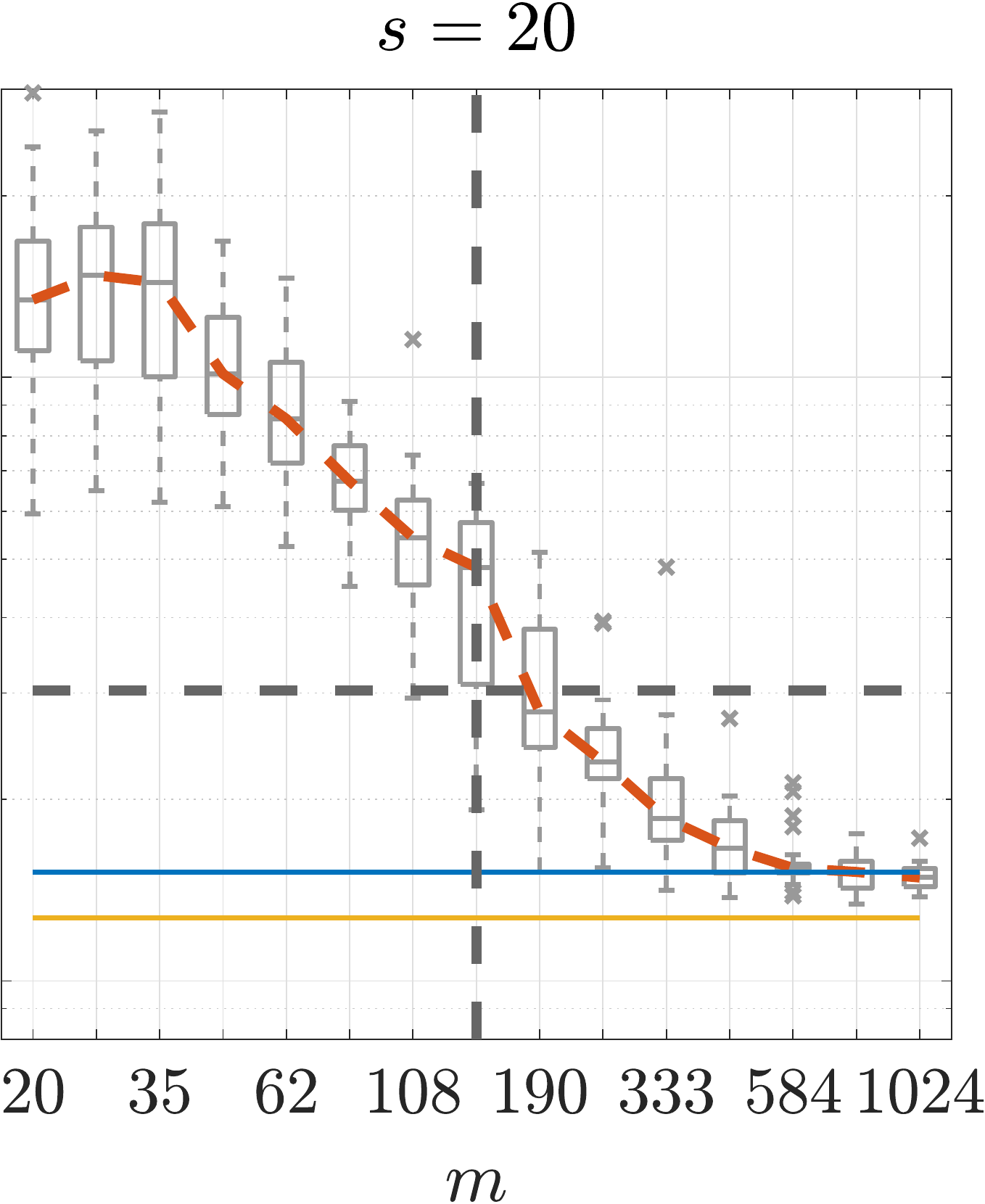} 
\includegraphics[width = 0.3\linewidth]{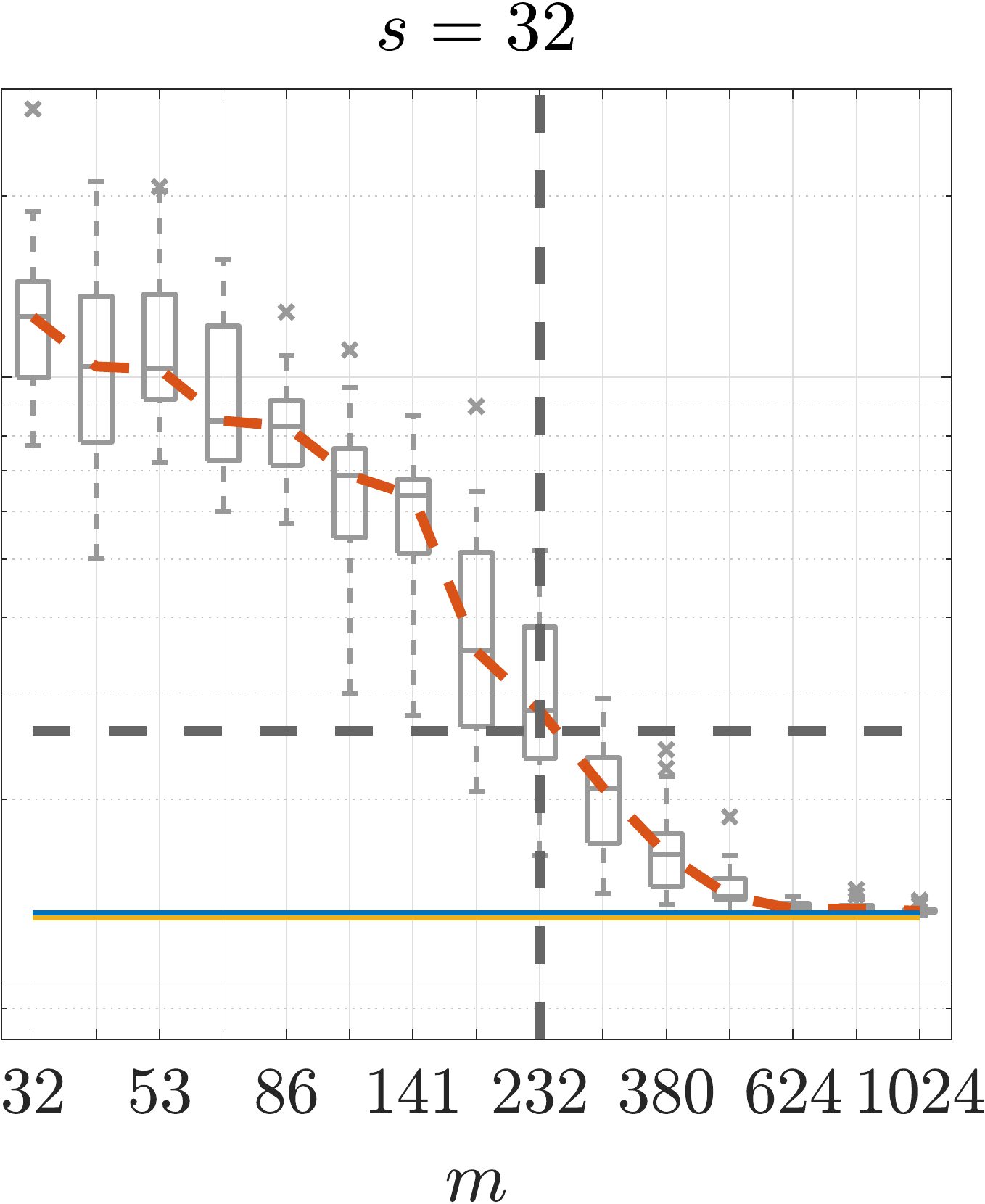}
\caption{\label{fig:m_vs_err_Gauss}Case study I (Gauss 2D). $D$-calibration test for $p = 2$ and $L = 5$; we show error vs.\ $m$ for different values of $s$.}
\end{figure}

We plot the \cossiga relative $H^1(\Omega)$-error as a
function of $m$ for $s = 13, 20, 32$ and $L=5$ in Figure~\ref{fig:m_vs_err_Gauss}.
The variability among the $n_{\mathrm{runs}}$ runs for each value of $s$ and $m$ is illustrated by using so-called box plots,
which are classical tools used in statistics to represent the variability of an ensemble of values.
More specifically, the rectangle extends from the $25$-th to the $75$-th percentile of the computed $n_{\mathrm{runs}}$ values;
the median ($50$-th percentile) is marked by a horizontal line inside the rectangle; whiskers (horiziontal
ticks connected to the rectangle by a line) mark the smallest and largest value out of the $n_{\mathrm{runs}}$ values
that are considered not to be outliers (in our case, the whiskers mark the $2.7$-th percentile
and $99.3$-th percentiles), and values exceeding these bounds are marked by ``cross'' markers.
In each plot, we add lines connecting the median values of the box plots to ease the visualization
of the convergence of \cossiga, and horizontal lines that mark the accuracy obtained by PG-OMP$(s)$ and PG-BS.
We can make several observations: 
\begin{itemize}
\item As $m$ increases, the accuracy of \cossiga eventually reaches the accuracy of
  PG-OMP$(s)$.
  Note that when the convergence curve of \cossiga approaches this bound, 
  it exhibits an ``elbow'', marking the point where the decay of the error with respect to $m$ slows significantly.

\item As $s$ increases (plots from  left to  right),
  the PG-OMP$(s)$ solution error decreases (as predicted already by Figure \ref{figure:C_calibration}),
  and reaches the PG-BS accuracy for $s$ large enough ($s = 32$ in Figure~\ref{fig:m_vs_err_Gauss}).
\item As we increase $s$, we need a larger $m$ to reach full accuracy.

\end{itemize}

\begin{figure}[t]
\centering
\begin{minipage}{0.45\textwidth}
\centering
\includegraphics[width = 0.95\linewidth]{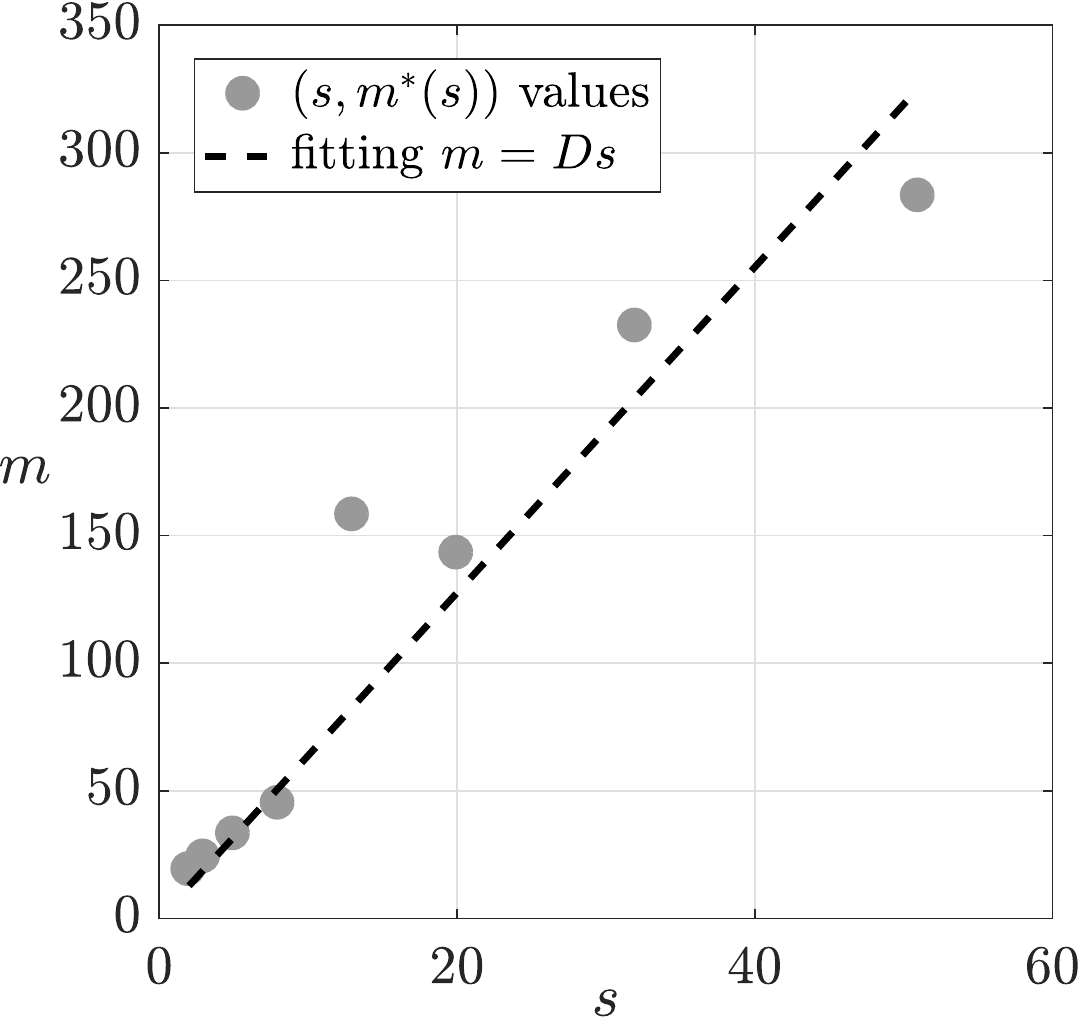}
\caption{\label{fig:D-calibration_Gauss}Case study I (Gauss 2D). $D$-calibration test for $p = 2$, $L = 5$.}
\end{minipage}
\begin{minipage}{0.05\textwidth}
\hspace{0.05\textwidth}
\end{minipage}
\begin{minipage}{0.4\textwidth}
\centering
\captionsetup{type=table} %% tell latex to change to table
\begin{tabular}{c|ccc}
$L$ $\backslash$ $p$ & 1 & 2 & 4 \\%& 5 \\
\hline
$4$ & 1.34 & 3.19 & 3.33 \\% & 2.74 \\
$5$ & 5.31 & 6.38 & 7.29 \\% & 4.97 \\
$6$ & 5.54 & 6.42 & 11.6 \\% & 17.2 \\
\end{tabular}  
\caption{\label{table:D-calibration_Gauss}Case study I (Gauss 2D). Numerical estimate of $D(p,L)$ obtained via $D$-calibration.}

\vspace{0.5cm}

\centering
\captionsetup{type=table}
\begin{tabular}{c|ccc}
$L$ $\backslash$ $p$ & 1 & 2 & 4 \\%& 5 \\
\hline
$4$ &  1.07\%  &  5.1\% & 11.7\% \\%&  21.9\%\\
$5$ &  4.25\%  & 10.2\% & 25.5\% \\%&  39.8\%\\
$6$ &  4.43\% & 10.3\% & 40.6\% \\%& 	[137\%]
\end{tabular}  
\caption{\label{table:m/N}Case study I (Gauss 2D). subsampling rate $m/N_{\mathrm{dof}}$ for different values of $p$ and $L$.}

\end{minipage}
\end{figure}

Now, let us explain how to perform $D$-calibration given the data computed in the setting above.
For every value of $s$,  we select the value of $m$ closest to the elbow of the convergence curve
(up to a prescribed relative tolerance) similarly to the case of $C$-calibration. In particular, we choose
\begin{equation}
\label{eq:m_choice}
m^*(s) = \arg\min_{s \leq m \leq N_{\mathrm{dof}}} 
\left|\|u-\tilde{u}_{\cossiga(s,m)}\|_{H^1(\Omega)} -  \mu \cdot \|u-\tilde{u}_{\textsc{OMP}(s)}\|_{H^1(\Omega)}\right|,  \quad \mu=2.
\end{equation}
With this choice, $\|u-\tilde{u}_{\cossiga(s,m^*(s))}\|_{H^1(\Omega)} \approx 2 \|u-\tilde{u}_{\textsc{OMP}(s)}\|_{H^1(\Omega)}$.
Recalling \eqref{eq:m=D*s^beta}, we can now find $D(p,L)$ by computing the best curve of the form $m = Ds$
fitting in the least-squares sense the data $(s,m^*(s))$ for the considered values of $s$.
Note that we are deliberately not considering the zero-order term in the equation, i.e., we are not fitting an affine model $m = Ds + m_0$,
since we are looking for a linear law of the form \eqref{eq:m=D*s^beta}.
Figure~\ref{fig:D-calibration_Gauss} illustrates this process. The resulting values of $D(p,L)$ are listed in Table~\ref{table:D-calibration_Gauss}.
As expected, $D$ increases overall with $L$ and $p$. 
%Some exceptions to this general rule appear 
%In fact, in a few cases $D$ decreases when moving from $p=4$ to $p=5$ for low values of $L$: this is due to the fact that for such values of $L$ we are not yet in the asymptotic regime, i.e., the PG-BS
%error of the approximations with $p=4$ is actually smaller that for $p=5$ (not shown for brevity).

\subsubsection{Convergence test}

We are now in a position to study the convergence of \cossiga with respect to the hierarchical level $L$
(or, equivalently, to the mesh size $h=2^{-L}$) for fixed $p$.
We consider different values of $L = 4,5,6$ and we study the recovery error as a function of the hierarchical level $L$.
For each value of $L$ and $p$ (we recall that in this test we consider B-splines of maximal regularity $C^{p-1}$),
we consider a the following \cossiga approximation (see Algorithm~\ref{alg:CossIGA}):
$$
\cossiga(p, p-1, L, s^*, m^*),
$$
with $s^* = C(p)N_{\mathrm{dof}}$ and $m^* = D(p,L) s^* = D(p,L)C(p)N_{\mathrm{dof}}$.
Table~\ref{table:m/N} illustrates the subsampling rate $m / N_{\mathrm{dof}} = C(p)D(p,L)$ (recall \eqref{eq:compression_ratio})
of \cossiga for different values of $p$ and $L$, using the constants $C(p)$ and $D(p,L)$
estimated in Sections~\ref{subsection:C-calibration} and \ref{subsection:D-calibration}.
In particular, with these choices we have %where  $s$ and $m$ are chosen such that 
$$
\|u-\tilde{u}_{\cossiga(s,m)}\|_{H^1(\Omega)} 
\approx 2 \|u-\tilde{u}_{\textsc{OMP}(s)}\|_{H^1(\Omega)}
\approx 4 \|u-\tilde{u}_{\textsc{BS}}\|_{H^1(\Omega)},
$$
i.e., we are losing a factor 4 from the best accuracy available at resolution $h$.
All the values in Table~\ref{table:m/N} are below $100\%$, hence corresponding
to a successful subsampling. %, apart from the case $p=5$ and $L = 6$.
% This means that the user is forced to choose $m$ and/or $s$ smaller than the values prescribed by the calibration tests in order to achieve
% any compression when $p$ and $L$ are both large (of course, this will lead to a further loss of accuracy).

Given that most of the values in Table~\ref{table:m/N} are significantly smaller than $100\%$, in the convergence test
we also investigate the effects of taking the constants $C$ and $D$ larger than prescribed by the calibration tests,
i.e., multiplying both $C$ and $D$ by a factor $\lambda \geq 1$: a choice of $\lambda$
strictly greater than 1 is expected to decrease the compression but improve the accuracy and robustness of \cossiga.
We also set a upper threshold, so that an minimal compression is always enforced. Specifically, we set
$$
m / N_{\mathrm{dof}} = \min \{ \lambda^2 C(p)D(p,L), 80\% \}.
$$ 
We repeat again $n_{\mathrm{runs}}=25$ tests of \cossiga for each value of $L$ using the above recipe.

\begin{figure}[t]
\centering
\includegraphics[width = 0.33\textwidth]{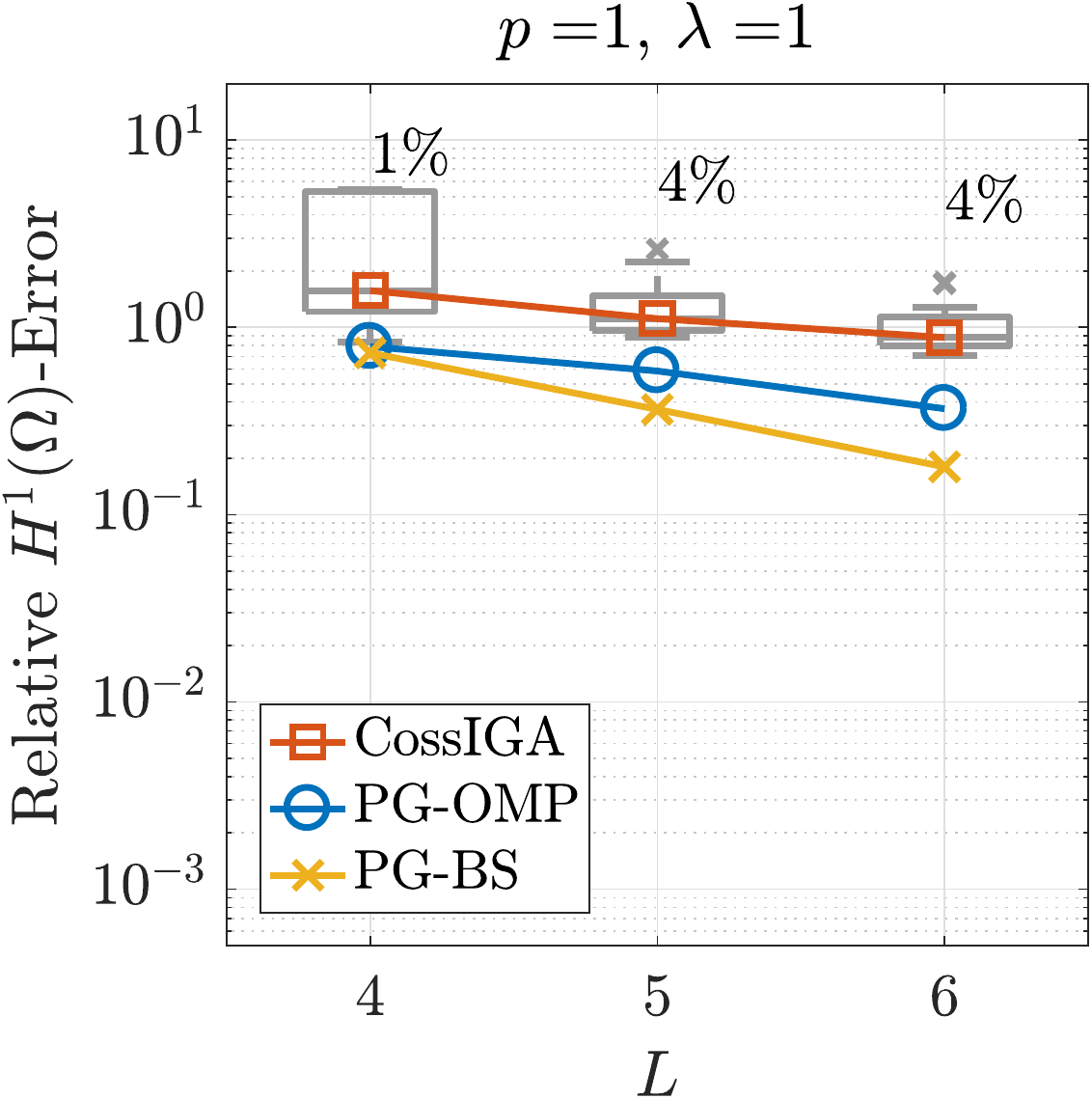}
\includegraphics[width = 0.262\linewidth]{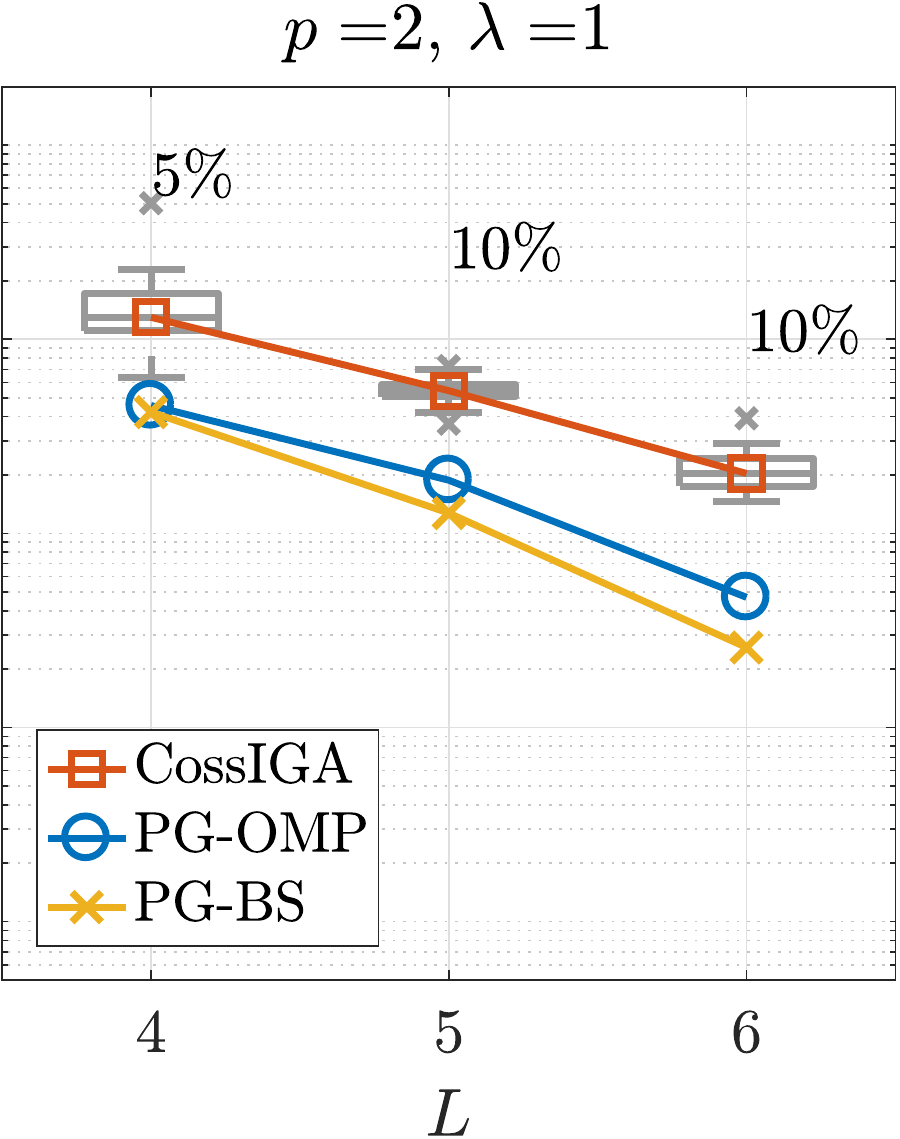}
\includegraphics[width = 0.262\linewidth]{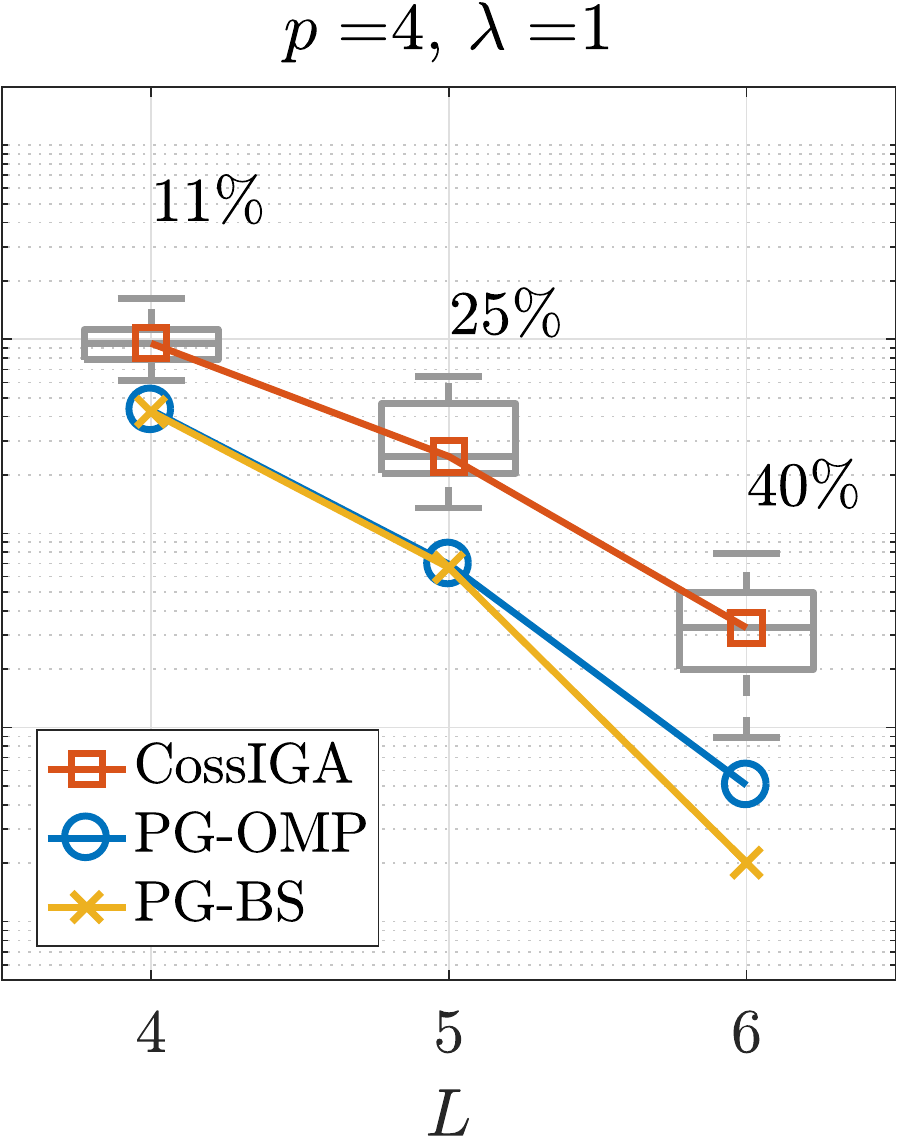}
\caption{\label{fig:convergence_Gauss}Case study I (Gauss 2D). Convergence analysis of \cossiga for the ``Gauss'' case study
for $p=1,2,4$ (left to right).
  The percentage above each box is the subsampling rate $m/N_{\mathrm{dof}}$, with $N_{\text{dof}}$ defined as in \eqref{eq:NDOFs}.}
\end{figure}

Results are reported in Figure~\ref{fig:convergence_Gauss} in terms of relative $H^1(\Omega)$-error vs.\ $L$.
We show the $n_{\text{runs}}$ 
values using box plots, and we add convergence curves for the
PG-BS solution and the PG-OMP$(s)$ solution with $s = s^*(p,L)$.
This Figure shows results for increasing $p$ and $\lambda=1$ (i.e., using $C$ and $D$ as calibrated).
We first note that, as expected, PG-OMP$(s)$ converges at a lower rate than PG-BS.
The convergence of \cossiga has an even lower rate, but the loss of accuracy of \cossiga
with respect to PG-OMP$(s)$ is moderate,
especially for lower degrees $p=1,2$, and considered the quite small subsampling rates imposed (reported in the plots with numbers above each box).

\newcommand{\convfigheight}{4cm}

%\clearpage

\begin{figure}[t]
  \centering
  \includegraphics[width = 0.33\textwidth]{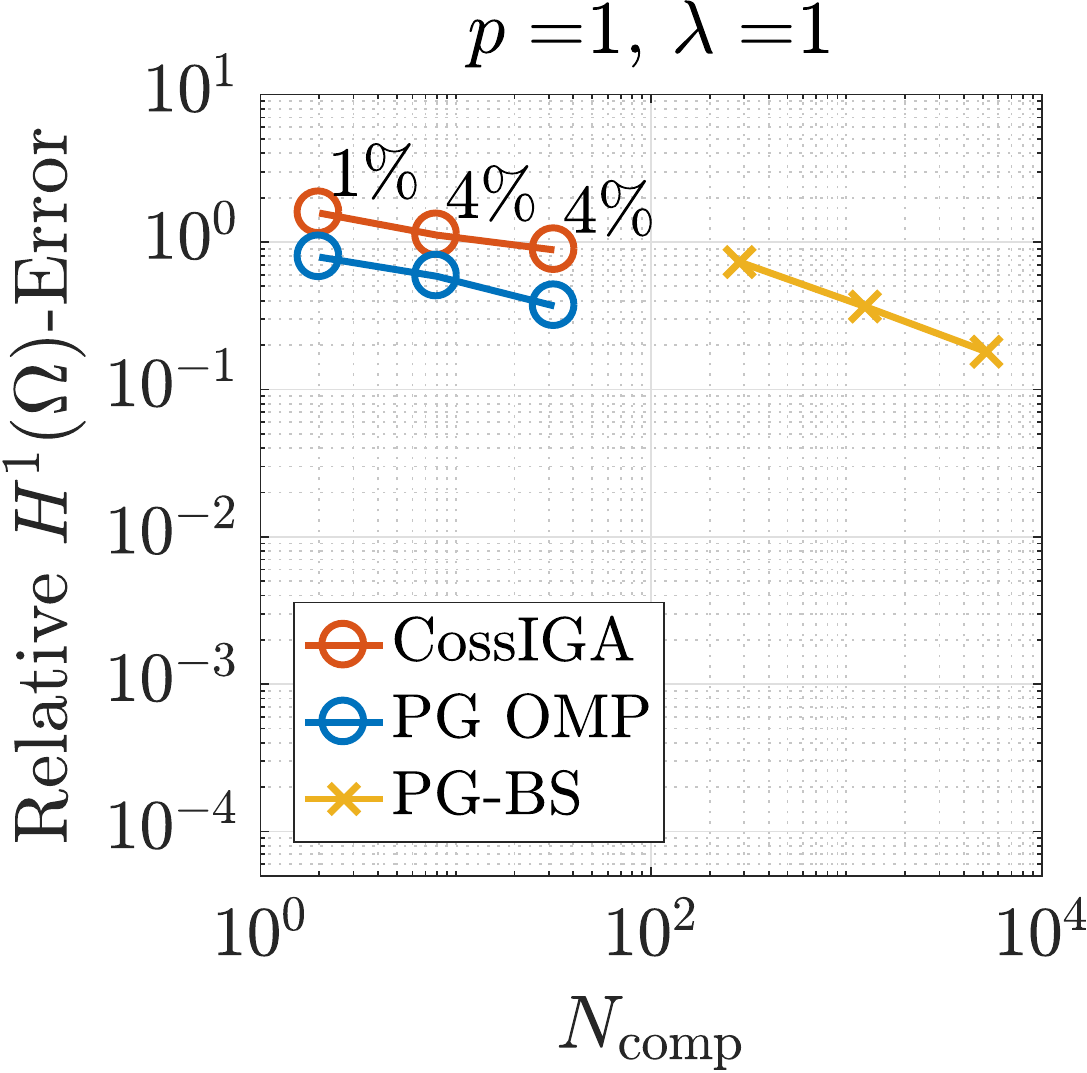}
  \includegraphics[width = 0.262\textwidth]{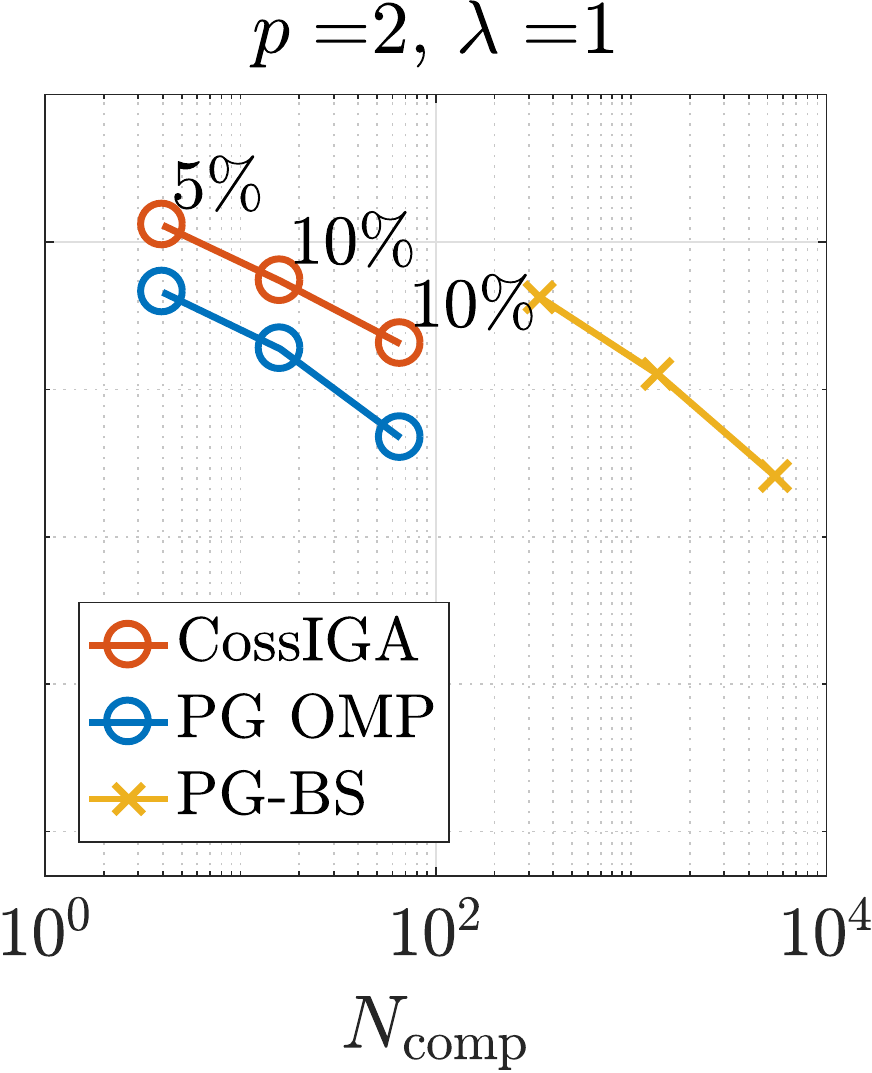}
  \includegraphics[width = 0.262\textwidth]{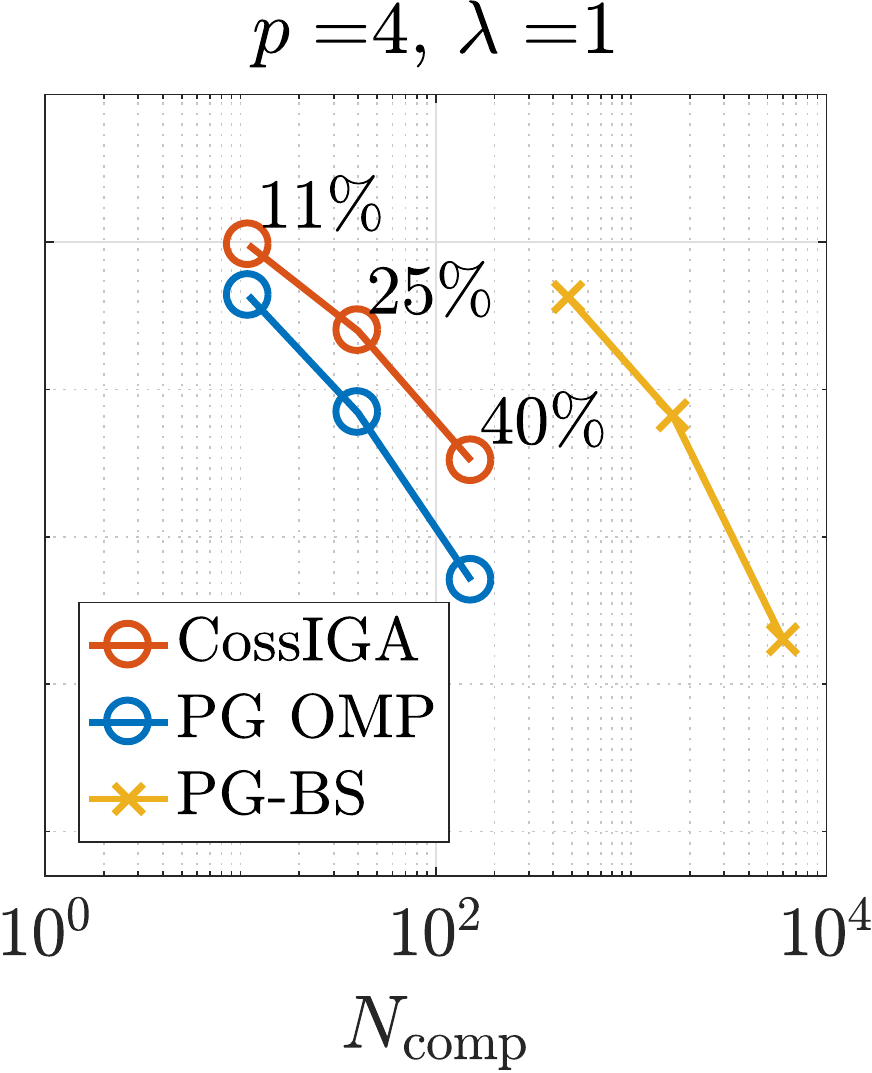} \\[10pt]
  \includegraphics[width = 0.33\textwidth]{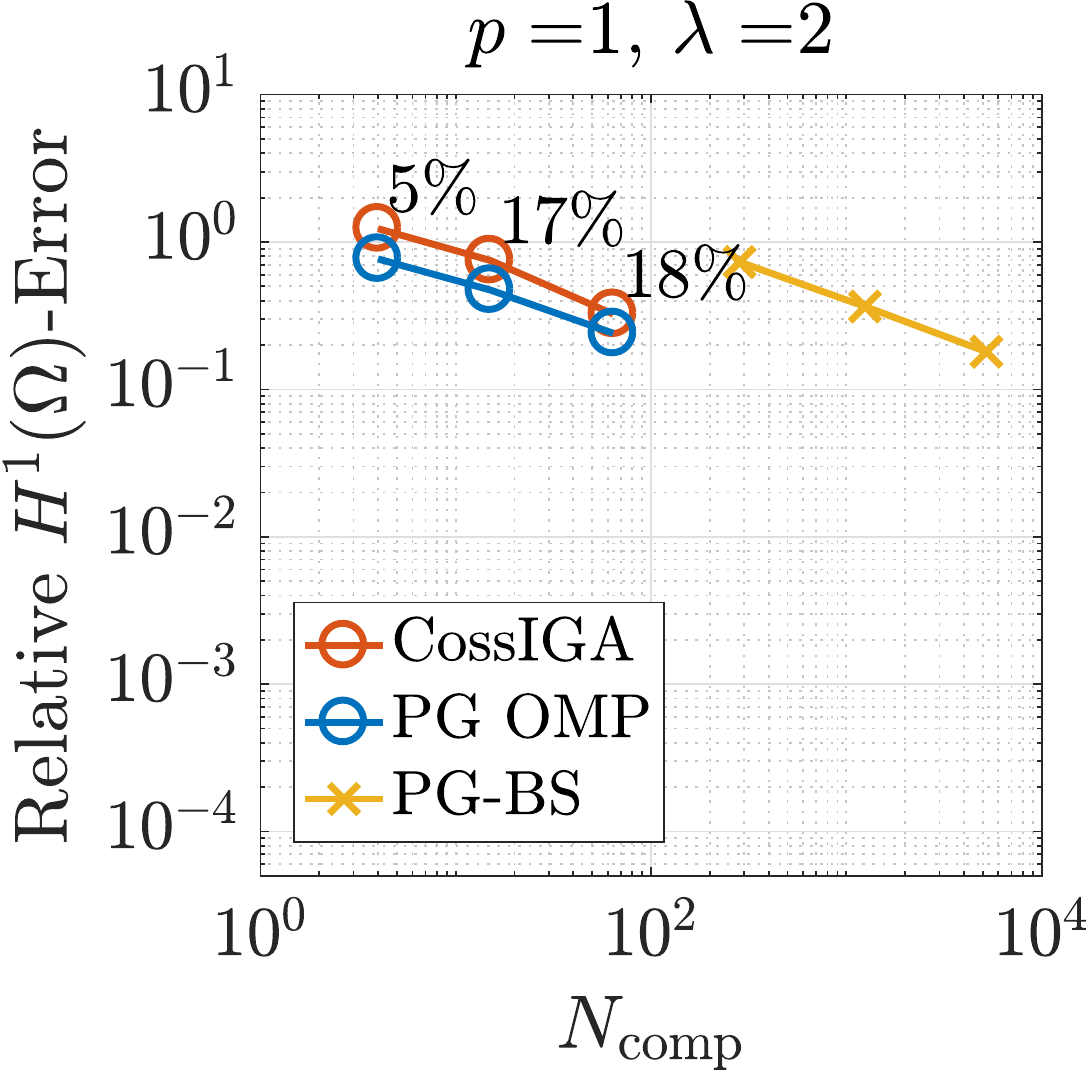}
  \includegraphics[width = 0.262\textwidth]{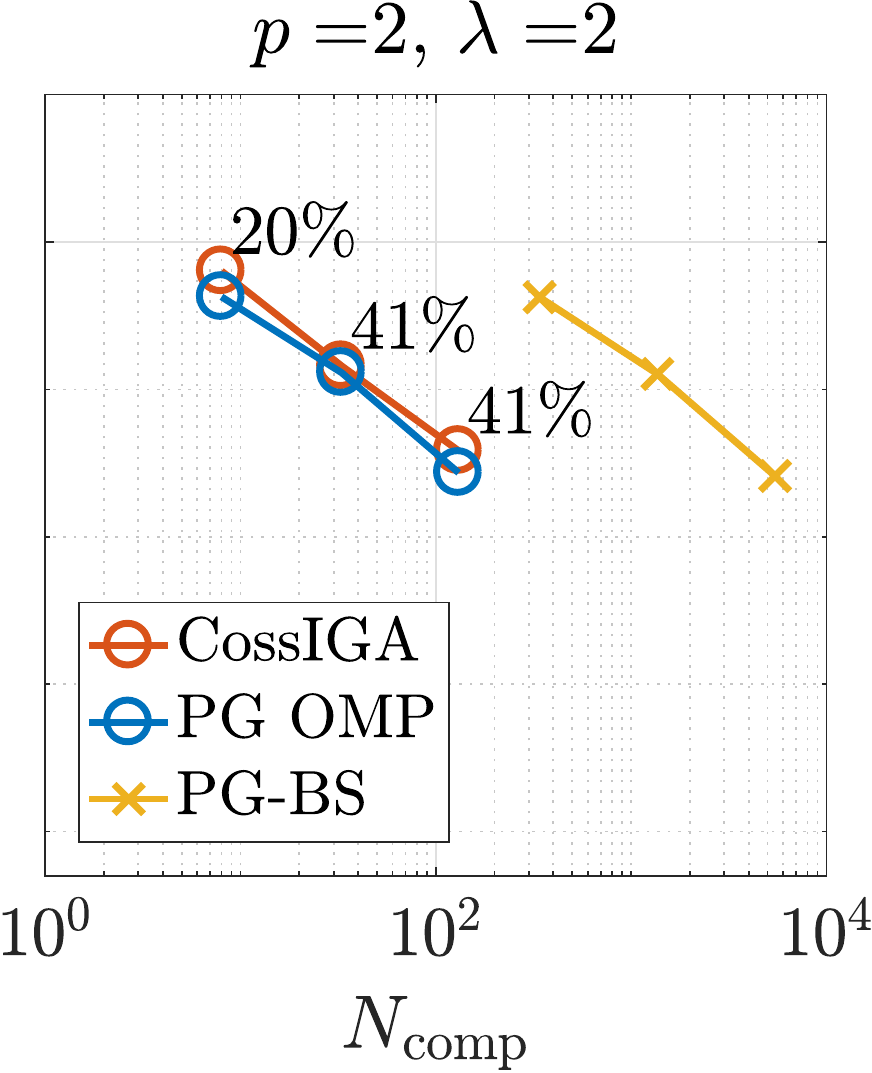}
  \includegraphics[width = 0.262\textwidth]{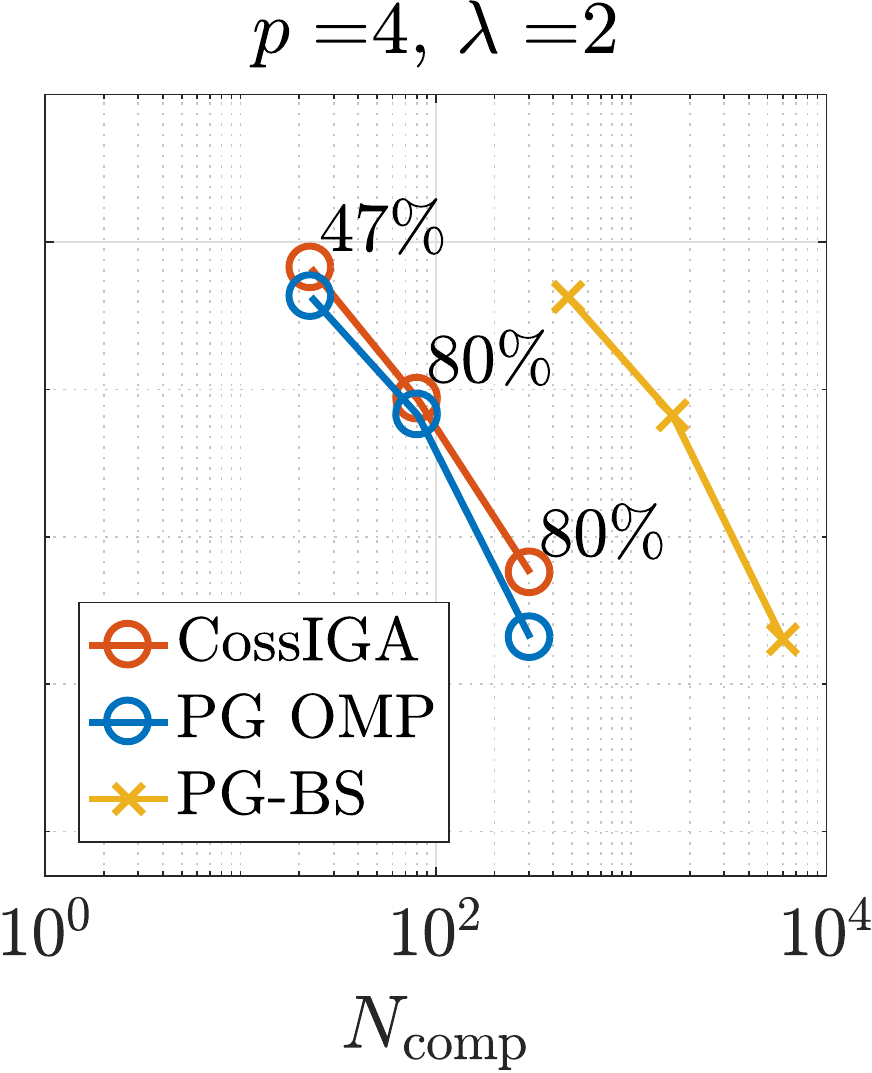}
  \caption{\label{fig:convergence_Gauss_dof}Case study I (Gauss 2D). Convergence analysis of \cossiga for the ``Gauss'' case study.
    In this figure we report error versus number of degrees of freedom.
    Top line: results for $p=1,2,4$ and $\lambda=1$ (left to right). Bottom line: same thing with $\lambda=2$.
    The percentage above each marker is the subsampling rate $m/N_{\mathrm{dof}}$, with $N_{\text{dof}}$ defined as in \eqref{eq:NDOFs}.}
\end{figure}

%Figure~\ref{fig:convergence_Gauss} was only showing plots of convergence with respect to $L$,
Convergence with respect to $L$, shown in Figure~\ref{fig:convergence_Gauss}, 
is not really a representative quantity for the computation cost of \cossiga. We therefore
also compare the methods by plotting the relative $H^1(\Omega)$-error as a function of the number of computed coefficients, defined by 
\begin{equation}
\label{eq:def_Ncomp}
N_{\mathrm{comp}}:= 
\begin{cases}
s, & \text{for PG-OMP$(s)$ and \cossiga},\\
N_{\mathrm{dict}}, & \text{for PG-BS},
\end{cases}
\end{equation}
(recall that the size $N_{\text{dict}}$ of the spline dictionary is comparable to the size $N_{\text{dof}}$ of its last hierarchical level -- see also the discussion after Equation~\eqref{eq:NDOFs}). This quantity is also not entirely representative of the actual computational cost.
It rather represents the optimal cost that can be achieved with ideal algorithms and implementation
(for comparison, the class of sublinear-time algorithms known as ``sparse Fourier transforms'' 
are able to recover an $s$-sparse signal of $\mathbb{R}^N$ with $O(s\polylog(N))$ flops from compressive Fourier measurements \cite{gilbert2014recent}).
Results are reported in Figure \ref{fig:convergence_Gauss_dof}
and show the rather significant improvement in convergence that could be potentially reached with
a careful implementation of \cossiga. In this Figure (and in most of the remaining ones of this paper),
we show only the median convergence of \cossiga instead of the box plots. In detail, the top half shows the results with $\lambda=1$, while the bottom one show
results for $\lambda=2$. As expected, the effect of setting $\lambda=2$ are: (i) PG-OMP is closer to PG-BS (since the number $s$
of coefficients of the solution that we are computing is doubled) (ii) \cossiga is closer to PG-OMP (since we are doubling
the number of rows that we sample from the PG matrix) (iii) the subsampling rate is four times larger, and in particular
for $p=4$ and large $L$ the threshold of 80\% is enforced.
% (specifically, for $p=5$ we observe that the convergence
% line of \cossiga departs from the line of PG-OMP, because we would need a subsampling rate above the threshold to stay close to it).

%\newpage

\subsection{Case study II: polyGauss 2D}\label{subsection:polygauss}

In this second test, we consider the Poisson problem with exact solution $u_{\mathrm{polyGauss}}$ defined in (\ref{eq:u_polyGauss}).
This solution has a clear multilevel structure since it is composed of a ``coarse component''
(the polynomial in (\ref{eq:u_polyGauss})) and a ``fine detail'' (the Gaussian peak),
while the coarse component was missing in the previous test.
In this test, as in the previous one, we consider B-splines of maximal regularity $C^{p-1}$.

We take advantage of this test also to verify the robustness of the method with respect to the calibration of the constants $C$ and $D$,
which is an expensive procedure (that one would rather do a limited number of times in advance, if not once --
or possibly skip altogether if theoretical estimates of $C$ and $D$ were available).
The results for the new calibrations are reported in Tables \ref{table:C-calibration_polyGauss} and \ref{table:D-calibration_polyGauss}
for $C$ and $D$, respectively, and they should be compared with the corresponding Tables~\ref{table:C-calibration_Gauss} and 
\ref{table:D-calibration_Gauss} obtained for the Gauss test case.

This comparison shows that the constant $C$ changes slightly (order of magnitude are identical though), as expected.
In particular, the values of $C$ for the polyGauss test are smaller than for Gauss,
which means that $u_{\mathrm{polyGauss}}$ is more compressible than $u_{\mathrm{Gauss}}$.
Perhaps more surprising (and against our assumptions) is that also the constants $D$ change, albeit being again in the same
range of magnitude. There is however no clear trend, i.e., sometimes the $D$ associated with the Gauss test is larger than the $D$ associated with the polyGauss test, and vice versa. We emphasize, however, that the calibration is a numerical procedure that can be sensitive to many tuning parameters
(e.g., choice of the values of $s,m$, number of runs per test $n_{\mathrm{runs}}$, tolerance factors, sampling probability distribution $\bm{\pi}$ for the test functions).  
Consequently, a conclusive statement on whether the value of $D$ is independent of the solution $u$ 
or not is hard to make and is postponed to further and deeper analyses.

\begin{figure}[t]
\centering
\includegraphics[width = 0.33\textwidth]{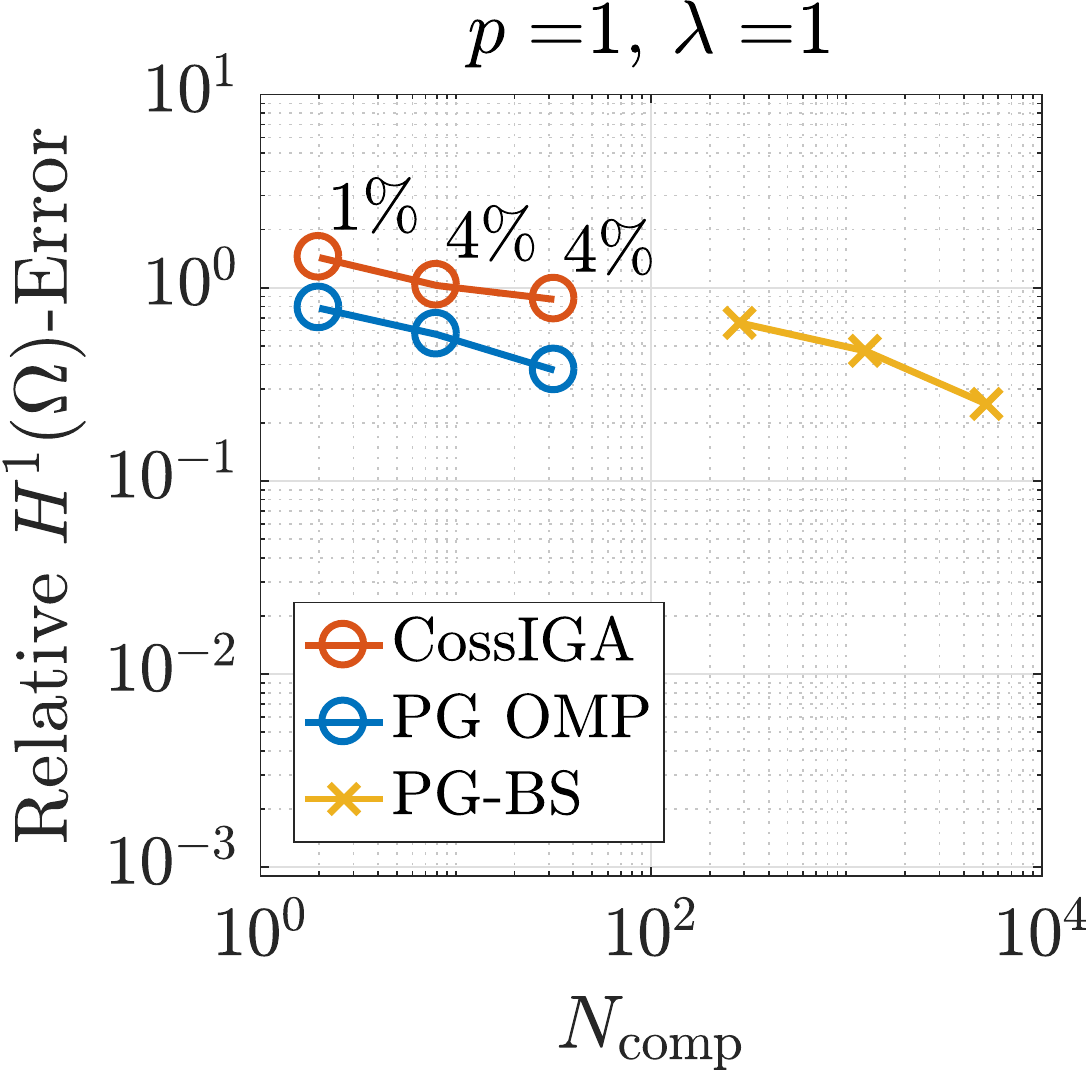}
\includegraphics[width = 0.264\textwidth]{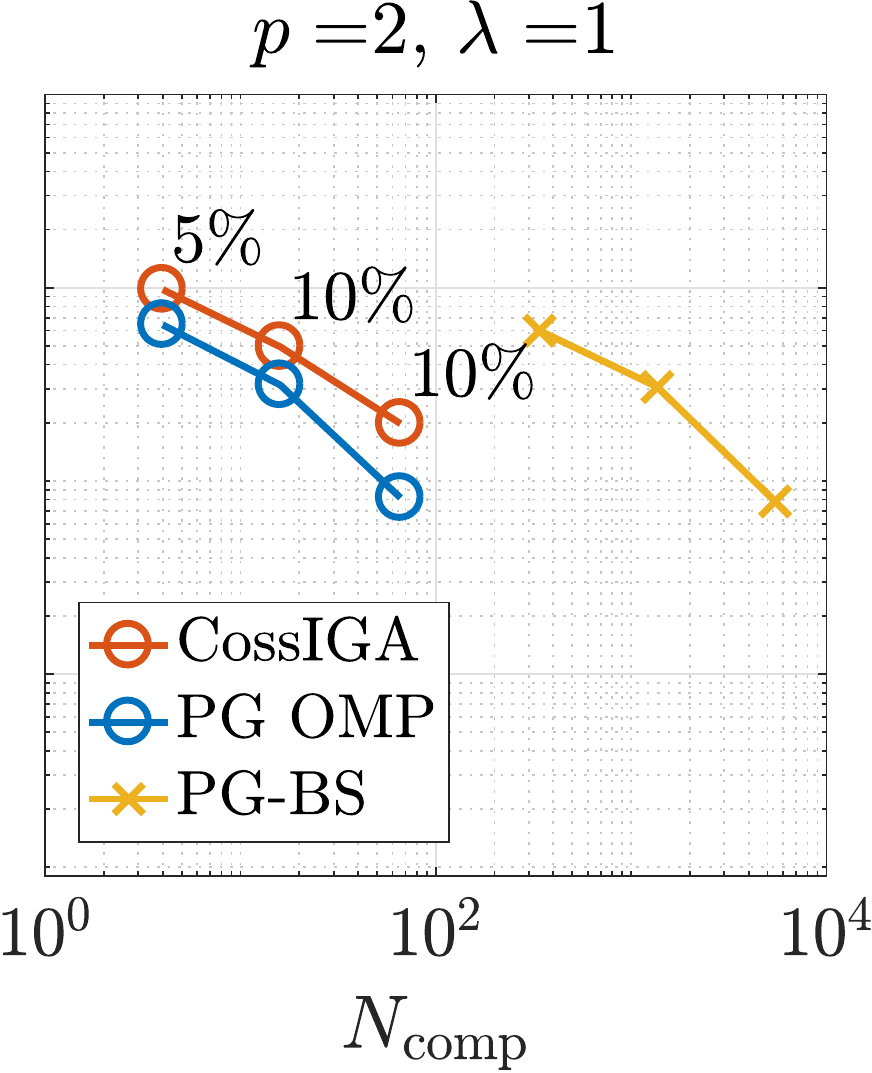}
\includegraphics[width = 0.264\textwidth]{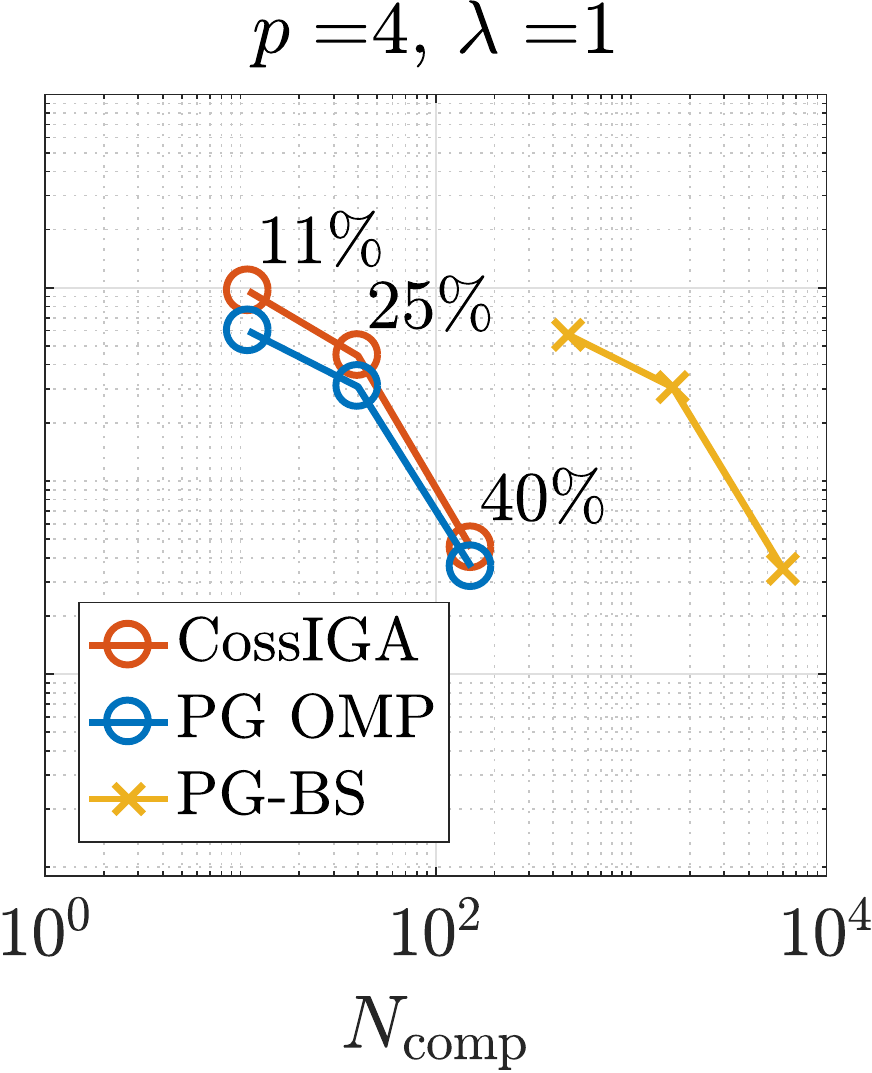}
\caption{\label{fig:convergence_polyGauss_C_Gauss}Case study II (polyGauss 2D). Convergence analysis of \cossiga
  for the polyGauss case study without performing $C$-calibration on $u_{\mathrm{polyGauss}}$.
  The percentage above each marker is the subsampling rate $m/N_{\mathrm{dof}}$, with $N_{\text{dof}}$ defined as in \eqref{eq:NDOFs}.}
\end{figure}

\begin{figure}[t]
  \centering
  \includegraphics[width = 0.33\textwidth]{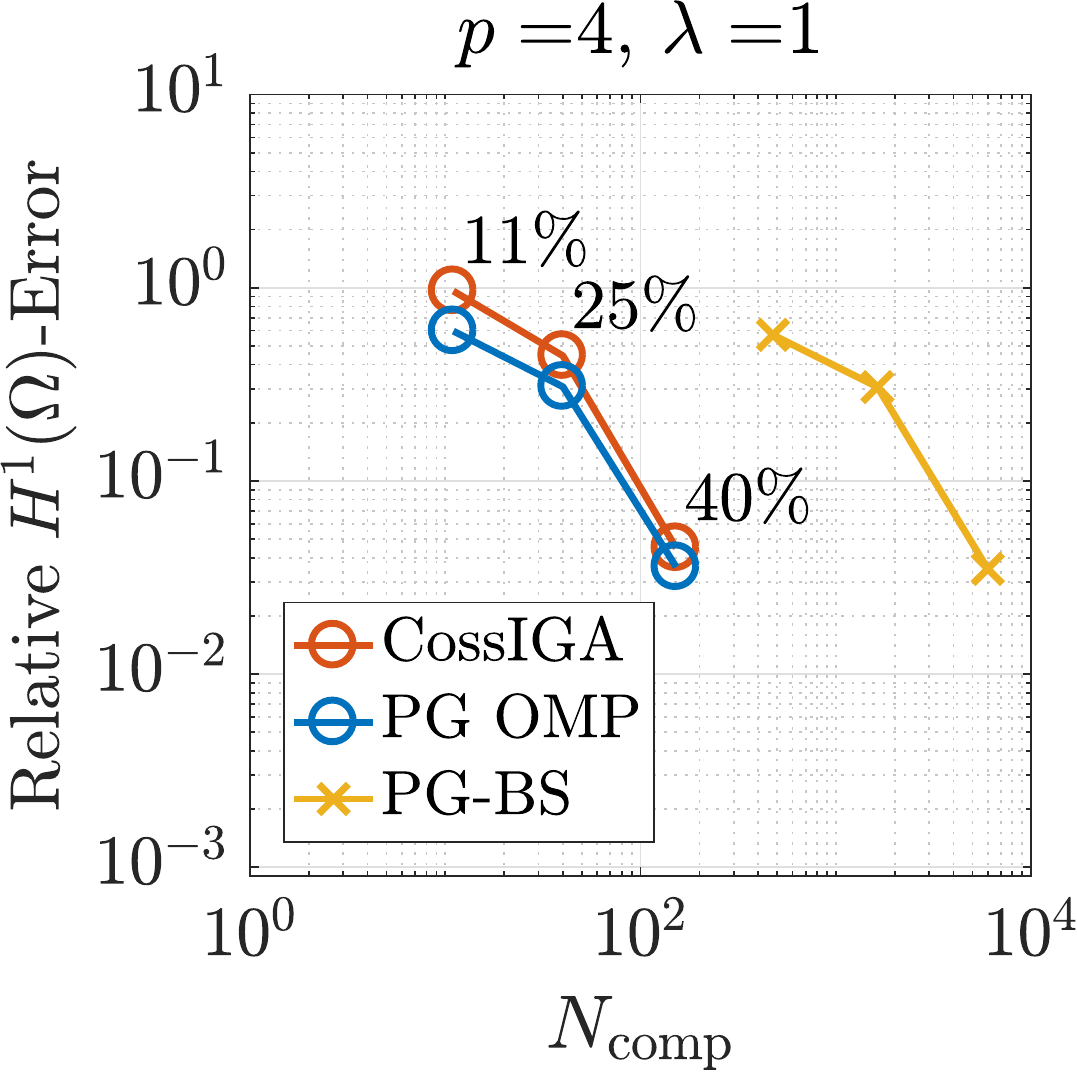}
  \includegraphics[width = 0.268\textwidth]{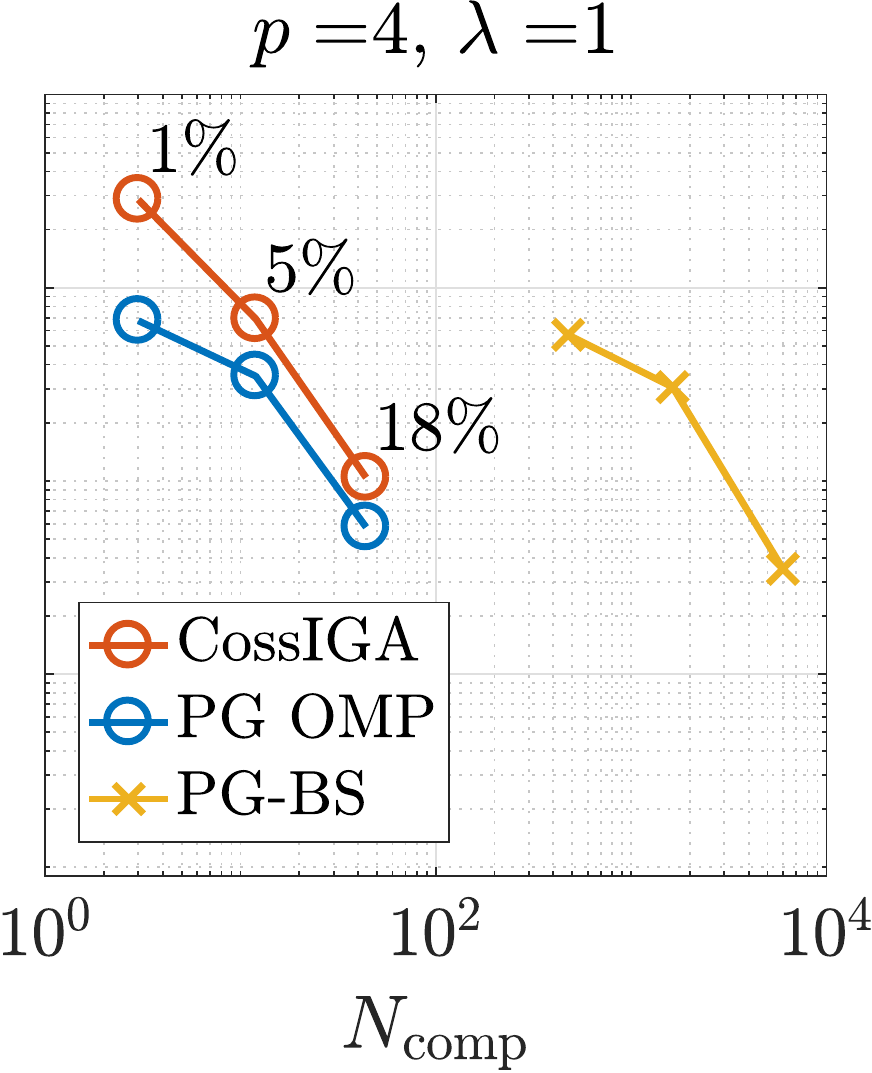}
\caption{\label{fig:convergence_polyGauss_C_polyGauss}Case study II (polyGauss 2D). Convergence analysis of \cossiga
  for the polyGauss case study with $C$-calibration on $u_{\mathrm{polyGauss}}$.
  Here we fix $p=4$ and we compare the following convergence plots.
  Left: uncalibrated $C,D,\lambda=1$ (this plot is also in Figure \ref{fig:convergence_polyGauss_C_Gauss});  
  Right: calibrated $C,D,\lambda=1$.
  The percentage above each marker is the subsampling rate $m/N_{\mathrm{dof}}$, with $N_{\text{dof}}$ defined as in \eqref{eq:NDOFs}.} 
\end{figure}

We compare the convergence plots when both $C$ and $D$ uncalibrated (i.e., using the constants for the
Gauss test in the polyGauss one) and calibrated for this test, with the aim of studying the sensitivity of
\cossiga with respect to the choice of these parameters.  The results are reported in
Figures \ref{fig:convergence_polyGauss_C_Gauss} and \ref{fig:convergence_polyGauss_C_polyGauss}.
The former shows results for $p=1,2,4$ and $\lambda=1$ for uncalibrated $C,D$ and the  latter compares
the results with calibrated and uncalibrated $C,D$ for selected values of $p$ and $\lambda$.
Figure \ref{fig:convergence_polyGauss_C_Gauss} shows that \cossiga is more
effective than in the Gauss test since the convergence of \cossiga
is closer to the convergence of PG-OMP, and the error reached by PG-OMP is closer to the error reached by PG-BS
than in the previous test (cf.\ Figure \ref{fig:convergence_Gauss_dof}).
This is due to the higher compressibility of the solution at hand: since we ``froze'' the subsampling rate but
  the solution considered in this test is more compressible, PG-OMP gets closer in error to PG-BS for the given number of coefficients $s$,
  and \cossiga gets closer to PG-OMP for the given number of sampled rows.

%The ``dual'' set of results appears instead if we compare the results of the calibrated and uncalibrated procedure in Figure \ref{fig:convergence_polyGauss_C_polyGauss} (case $p=4$ only for brevity; other values of $p$ give similar results).
In Figure \ref{fig:convergence_polyGauss_C_polyGauss} we compare the results of the calibrated
and uncalibrated procedure for $p=4$ (other values of $p$ give similar results).
We see that in the calibrated case, the convergence of PG-OMP is actually further from the PG-BS results,
which is to be expected since $C$ is substantially smaller after recalibration, so less terms are computed;
and similarly, the convergence of \cossiga is further from PG-OMP because less rows are now computed.
In other words, as one would expect, the results obtained with the uncalibrated constants are
(in this case) suboptimal, in the sense that the same target accuracy relative to PG-BS
can be obtained with a smaller subsampling rate (i.e., with smaller values of $C$ and $D$). 

\begin{table}[t]
\begin{minipage}[t]{0.53\linewidth}
\centering
\begin{tabular}{c|ccc}
$p$ & 1 & 2 & 4 \\%& 5 \\
\hline
$C(p)$ &  $4.0\cdot 10^{-3}$ &  $4.6\cdot 10^{-3}$ &  $1.0\cdot 10^{-2}$ %&    $7.0\cdot 10^{-3}$\\
\end{tabular}
\caption{\label{table:C-calibration_polyGauss}Case study II (polyGauss 2D). Numerical estimate of $C(p)$ computed via the $C$-calibration test for the function $u_{\mathrm{polyGauss}}$.}
\end{minipage}
\quad
\begin{minipage}[t]{0.45\linewidth}
\centering
\begin{tabular}{c|cccc}
$L$ $\backslash$ $p$ & 1 & 2 & 4 \\%& 5 \\
\hline
$4$ & 1.80 & 1.00 & 1.00  \\%& 3.69 \\
$5$ & 1.76 & 4.24 & 4.64  \\%& 5.43 \\
$6$ & 4.48 & 11.45& 17.77 \\%& 11.53 \\
\end{tabular}  
\caption{\label{table:D-calibration_polyGauss}Case study II (polyGauss 2D). Numerical estimate of $D(p,L)$ obtained via $D$-calibration for the polyGauss test.} 
\end{minipage}
\end{table}

\subsection{Case study III: $C^0$ vs.\ $C^{p-1}$ splines 2D}\label{subsection:C0}

In this test, we assess the performance of \cossiga when $C^0$ splines are employed instead of $C^{(p-1)}$ splines,
motivated by the fact that $C^0$ splines are supported on one or two elements only (instead of $p+1$ elements as in the case of $C^{(p-1)}$ splines), which might further promote sparsity of solutions with localized features.
On the other hand, it is well-known that $C^{(p-1)}$ splines yield a better accuracy per degree of freedom, see e.g. \cite{BBSV14},
so it is not clear \emph{a priori} what choice should more favorable in terms of error-dof ratio.
We limit ourselves to the polyGauss test and recalibrate once more $C,D$ for this test.

\begin{figure}[t]
  \centering
  \includegraphics[height = 0.32\textwidth]{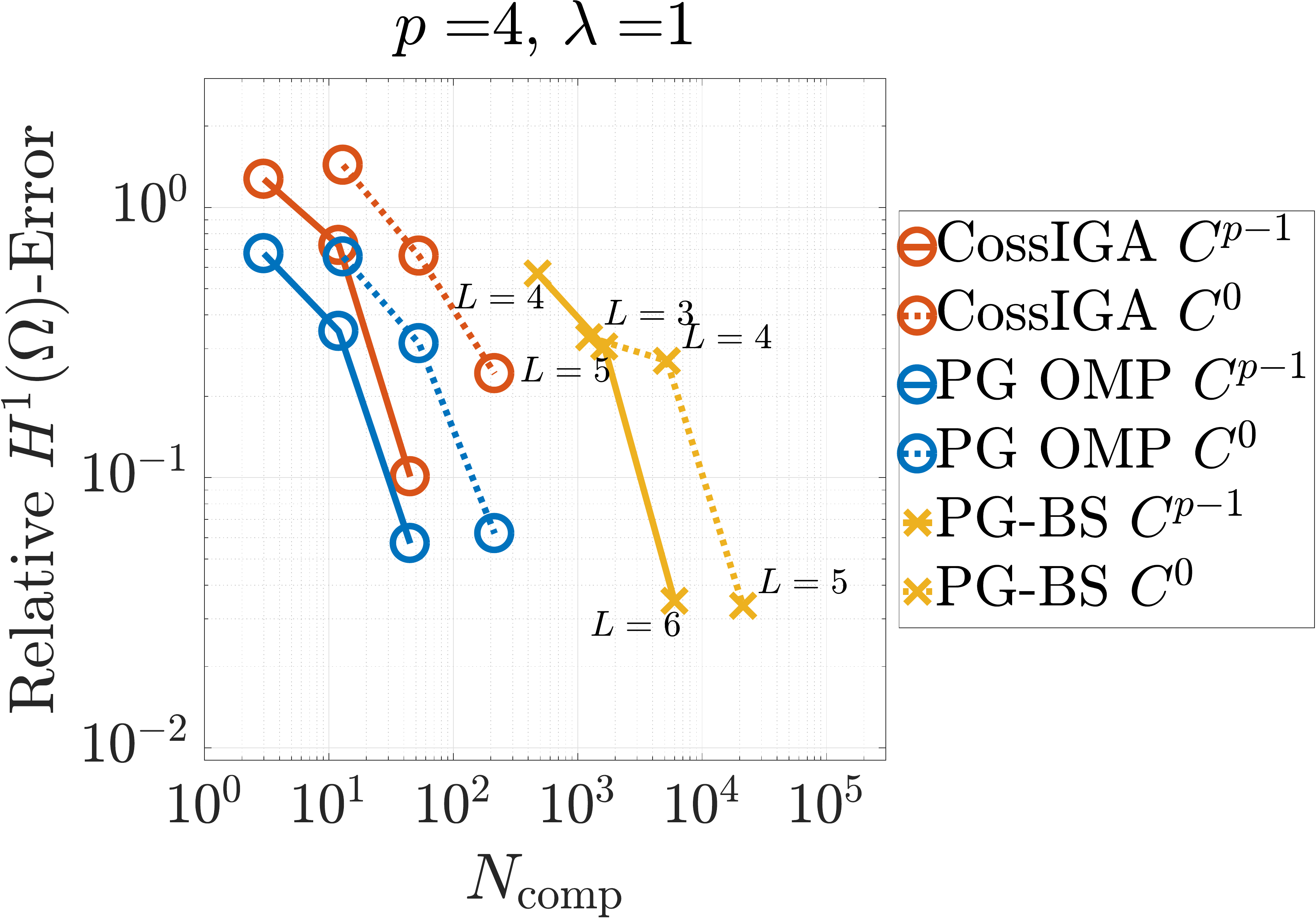}
	\caption{Case study III ($C^0$ vs. $C^{p-1}$ splines 2D). Results obtained with $p = 4$, $\lambda = 1$. We use different discretization levels for the $C^0$ case ($L=3,4,5$) and for the $C^{p-1}$ case ($L=4,5,6$).}
	\label{fig:C0_vs_Cmax}
\end{figure}

In Figure \ref{fig:C0_vs_Cmax}, we show the results obtained with $p=4, \lambda=1$.
In the case of $C^0$ splines we consider $L=4,5,6$ and represent the results with full lines, while in the case of $C^0$ splines we consider $L=3,4,5$ and represent the results with dashed lines.
The use of different discretization levels for the two approaches allows to make a better comparison.
Indeed, as can be seen in the plot, in this case the error we obtain with PG-BS and PG-OMP using $C^{p-1}$ splines for a given discretization level is almost identical to the error obtained using $C^{0}$ splines for the following discretization level.

The $C^{p-1}$ approach, however, yields a lower number of degrees of freedom, showing an advantage over the $C^0$ case.
This is assessed by the distance between the PG-BS lines, and a similar distance (possibly a bit larger) can be observed between the PG-OMP lines is approximately the same. 
The advantage of $C^{p-1}$ splines is then of course inherited by the \cossiga results. 
This suggests that the better error-dof ratio yielded by high regularity holds also in the context of compressed sensing. 

\renewcommand{\convfigheight}{4cm}

\subsection{Case study IV: polyGauss 3D}\label{subsection:3d}

% Here just run an example with good subsampling rate to show that the method works. polyGauss with $C^{(p-1)}$ splines, $p=2$, C and D calibrated as
% \begin{itemize}
% \item C=0.011,
% \item D=2.20 14.5 9.51  for L=2,3,4
% \end{itemize}

In this case study, we consider the three-dimensional version of the Poisson problem with exact solution in Equation (\ref{eq:u_polyGauss_3D}). 
Intuitively, we expect this solution to be even sparser than the 2D equivalent since
the localized feature (i.e., the exponential term in \eqref{eq:u_polyGauss_3D}) is essentially supported on the
horizontal mid-plan of the domain, and is zero in most of the rest of the volume.
Thus, the setup is ideal for \cossiga. In this test, we recalibrate the constants $C,D$ and we fix $p=2$ for simplicity.
An immediate verification of the fact that the solution is much sparser and compressible than before is that now we can choose a much
smaller constant $\mu$ in the calibration of both $C$ and $D$, cf.\ Equations (\ref{eq:s^*_choice}) and (\ref{eq:m_choice}):
specifically, we choose $\mu=1.01$ instead of $\mu=2$  for the $C$-calibration (which means that much fewer coefficients
are significantly nonzero) and $\mu=1.2$ instead of $\mu=2$ for the $D$-calibration (and, hence, that fewer rows of the matrix
are needed to recover a satisfactory approximation of the solution). 

\begin{figure}[t]
\centering
\includegraphics[height = 0.32\linewidth]{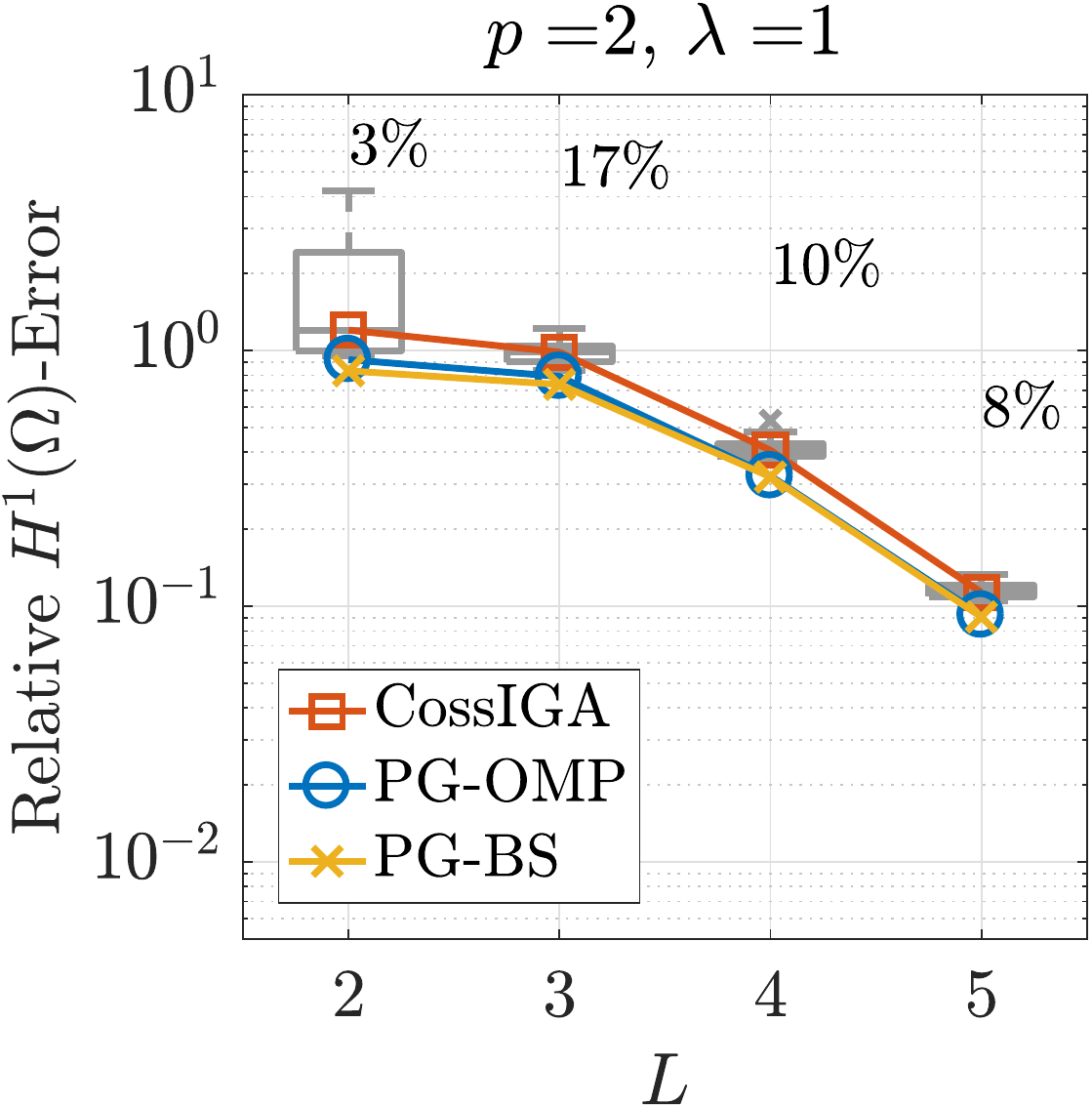} \quad \quad
\includegraphics[height = 0.32\linewidth]{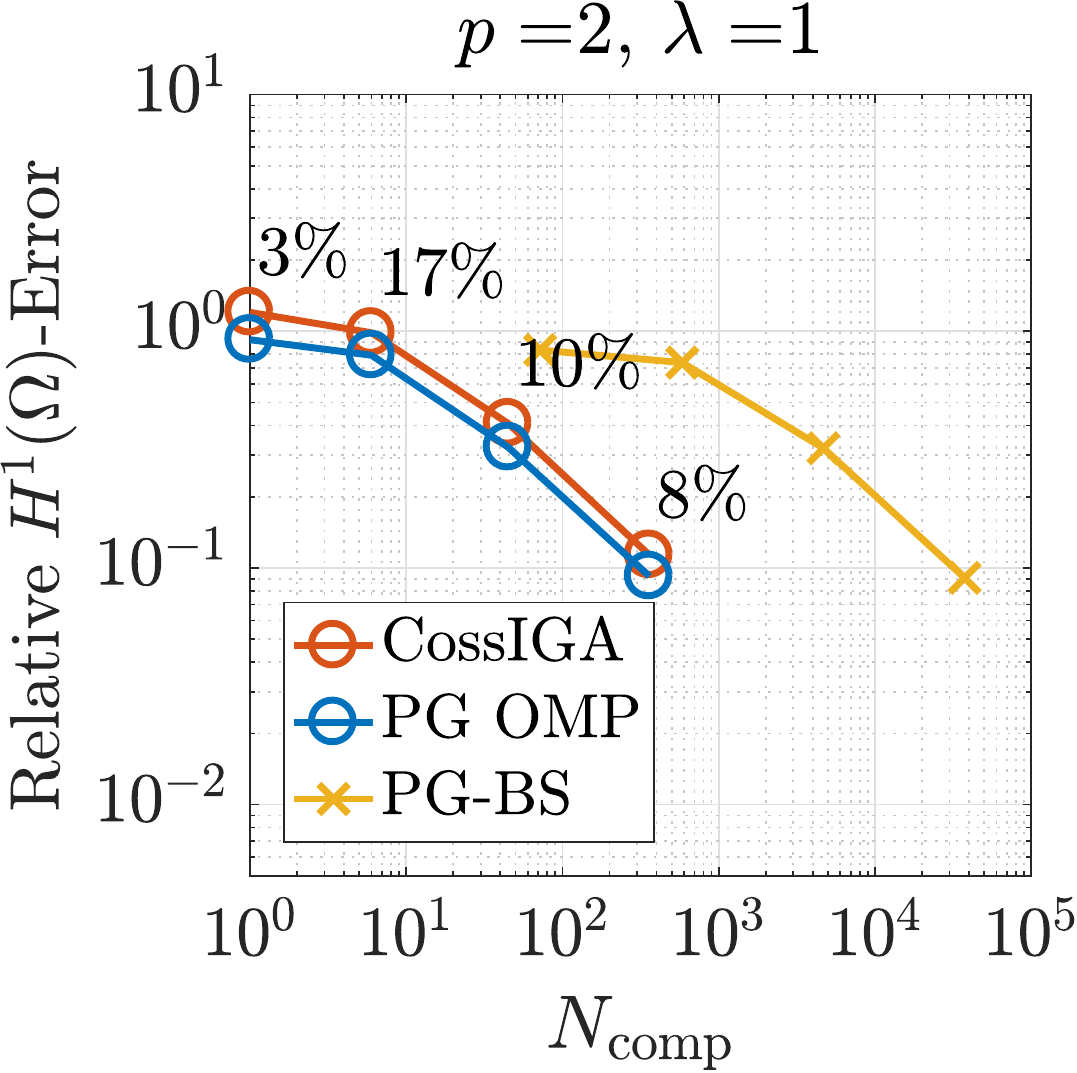}
\caption{Case study IV (polyGauss 3D). Results for $p=2$, $C^{(p-1)}$ splines, $\lambda=1$. Left: error vs $L$; right: error vs $N_{\mathrm{dof}}$. The percentage above each box/marker is the subsampling rate $m/N_{\mathrm{dof}}$, with $N_{\text{dof}}$ defined as in \eqref{eq:NDOFs}. }\label{fig:polyGauss3D}
\end{figure}

Given these premises, we expect a good performance of
\cossiga and indeed this is what can be deduced from the convergence plots reported in Figure \ref{fig:polyGauss3D}, 
which shows the convergence of the error with respect to both the refinement level $L$
and the number of computed coefficients $N_{\mathrm{comp}}$, defined in Equation \eqref{eq:def_Ncomp},
in the left and right panel, respectively.
The boxes in the left panel are very thin. This means that the variability due to randomness is almost negligible.
Moreover, in the left panel the \cossiga median convergence is very close to the one of PG-BS.
This  can be explained by the smaller values of $\mu$ used in the $C$- and $D$-calibration procedures.
In the right panel we see that we can essentially recover the accuracy of the full PG-BS solution with two orders
of magnitude less degrees of freedom and an overall subsampling rate smaller than $20\%$.

% \newpage

\section{Conclusions}\label{section:conclusions}

We have shown that the compressive sensing paradigm can be successfully applied to
solve PDEs on domains with a nontrivial geometry, and that the sparsity
principle can be leveraged to discretize PDEs by means of a
compressive Petrov-Galerkin discretization, leading to the
\cossiga (COmpreSSive IsoGeometric Analysis) method.
This paper is essentially a proof of concept. Its aim is to show
that the proposed method can be an attractive alternative to speed up IGA solvers
whenever the solution can be expressed over a basis (or, possibly, a dictionary)
that enhances its sparsity. 

From a theoretical perspective, many issues remain open, such as
estimating the local $a$-coherence and providing effective \emph{a priori}
estimates for the constants $C$ and $D$, which are two of the
cornerstones to make the method effective. 
These issues should be addressed in order to prove a formal convergence theorem for \cossiga,
which hence seems far from being a trivial task. We also note that the calibration procedure used to estimate  $C$ and $D$ (or, equivalently, $s$ and $m$) requires multiple runs of \cossiga and it is therefore not computationally efficient. Although the optimization of the calibration process is an interesting open problem, the fine tuning of $C$ and $D$ does not seem to play a crucial role in practice. In fact, recalling Figure~\ref{fig:convergence_polyGauss_C_polyGauss}, the method numerically converges in the uncalibrated scenario as well (i.e., for constants $C$ and $D$ calibrated on a different problem).

Concerning the computational efficiency, we remark that our Matlab
implementation of \cossiga is not optimized yet.  As a consequence,
the different algorithms (\cossiga, PG-BS, PG-OMP, and IGA) 
considered in this paper were compared only in terms of number of
computed coefficients, and not in terms of computational
time, even though the former choice does not take into account the
fact that the corresponding matrices have different sparsity patterns.
Furthermore, we remark that even the standard implementation usually employed in
CORSING and compressed sensing, where one assembles only the (subsampled)
matrix $A$ and uses OMP for sparse recovery, would be far from the ideal
computational cost $N_{\text{comp}} = s$ (recall
Equation~\eqref{eq:def_Ncomp}), since such
cost would depend on the number of columns $N_{\text{dict}}$.
As mentioned in Section~\ref{section:tests}, a promising research direction is the use
of sublinear-time algorithms to bridge this gap (note that this type
of algorithms do not require the assembly of $A$, but only fast access
to its entries). However, sublinear-time algorithms are not
available yet for the type of matrices considered in CORSING and in \cossiga.
Recent work in this direction can be found in \cite{choi2020sparse, choi2019sparse},
where algorithms of this kind were applied to high-dimensional function approximation
in the context of random sampling from bounded orthonormal systems.
To be applied to CORSING and \cossiga, they should be generalized to the case of random sampling
from Riesz bases and dictionaries, respectively. These
extensions and their efficient implementation for CORSING and \cossiga
are promising open directions, currently under investigation.

Finally, the connections between the \cossiga approach and the standard
local adaptivity algorithms for IGA deserve further investigations.

\section{Acknowledgements}
SB acknowledges the PIMS Postdoctoral Training Centre in Stochastics, NSERC through grant R611675, and the Faculty of Arts and Science of Concordia University for the financial support. The authors thank Ben Adcock for fruitful discussions about \cossiga and for supporting LT's visits at SFU in 2018 and 2019, partially funded by the PIMS CRG in ``High-dimensional Data Analysis''. The authors also thank Fabio Nobile and John Evans for their feedback on an earlier version of this manuscript.
LT and MT also received support from the Gruppo Nazionale Calcolo Scientifico-Istituto Nazionale di Alta Matematica
``Francesco Severi'' (GNCS-INDAM).

\bibliographystyle{plain}
\bibliography{biblio,IGA_biblio}

\end{document}